\documentclass[11pt]{amsart}

\usepackage{amssymb,amsmath,amsthm}
\usepackage{graphicx}
\usepackage{subfigure}
\usepackage{multicol}
\usepackage{caption}

\makeindex

\usepackage{color}
\usepackage{enumerate,amssymb}
\usepackage[arrow,curve,matrix,tips,2cell]{xy}
  \SelectTips{eu}{10} \UseTips
  \UseAllTwocells

\theoremstyle{plain}
\newtheorem{theorem}{Theorem}[section]

\newtheorem{corollary}[theorem]{Corollary}
\newtheorem{conjecture}[theorem]{Conjecture}

\theoremstyle{definition}

\newtheorem{definition}[theorem]{Definition}
\newtheorem{example}[theorem]{Example}

\newtheorem{remark}[theorem]{Remark}
\newtheorem{notation}[theorem]{Notation}

\newtheorem{convention}[theorem]{Convention}

{\catcode`@=11\global\let\c@equation=\c@theorem}

\newcommand{\id}{\operatorname{id}}
\newcommand{\Tor}{\operatorname{Tor}}
\newcommand{\colim}{\operatorname{colim}}
\newcommand{\Ext}{\operatorname{Ext}}
\newcommand{\supp}{\operatorname{supp}}
\newcommand{\im}{\operatorname{im}}
\newcommand{\alg}{\operatorname{alg}}

\numberwithin{figure}{section} 
\numberwithin{table}{section} 
\makeatletter
\let\c@table\c@theorem
\let\c@figure\c@theorem
\makeatother

\usepackage{hyperref}  
\hypersetup{%
  bookmarksnumbered=true,%
  bookmarks=true,%
  colorlinks=true,%
  linkcolor=blue,%
  citecolor=blue,%
  filecolor=blue,%
  menucolor=blue,%
  pagecolor=blue,%
  urlcolor=blue,%
  pdfnewwindow=true,%
  pdfstartview=FitBH}
  
\usepackage{xcolor}
\usepackage{makecell}
\usepackage{url}
\usepackage{soul}
\usepackage{caption}

\newcommand{\upar}[2]{{\uparrow}_{#1}^{#2}}
\newcommand{\wS}{\mathbb{S}_n}
\newcommand{\wG}{{\mathbb{G}_n}}

\newcommand{\abGal}{\mathrm{Gal}}

\newcommand{\wK}{{K_n}}

\newcommand{\smsh}{\wedge}
\newcommand{\GG}{\mathbb{G}}
\newcommand{\bbS}{\mathbb{S}}
\newcommand{\WW}{\mathbb{W}}
\newcommand{\W}{\mathbb{W}}
\newcommand{\FF}{\mathbb{F}}
\newcommand{\F}{\mathbb{F}}

\newcommand{\N}{\mathbb{N}}
\newcommand{\Z}{\mathbb{Z}}
\newcommand{\D}{\mathbb{D}}
\newcommand{\Q}{\mathbb{Q}}

\newcommand{\C}{\mathbb{C}}
\newcommand{\Pic}{\mathrm{Pic}}
\newcommand{\Gal}{\mathrm{Gal}}
\newcommand{\Mod}{\mathrm{Mod}}
\newcommand{\Hom}{\mathrm{Hom}}
\newcommand{\End}{\mathrm{End}}
\newcommand{\Aut}{\mathrm{Aut}}
\newcommand{\ndiv}{\hspace{-4pt}\not|\hspace{2pt}}
\newcommand{\bF}{\overline{\mathbb{F}}}

\usepackage[nameinlink,capitalise,noabbrev]{cleveref}
\usepackage{stmaryrd}

\newcommand{\bP}{\mathbb{P}}

\newcommand{\cA}{\mathcal{A}}
\newcommand{\cC}{\mathcal{C}}
\newcommand{\cD}{\mathcal{D}}
\newcommand{\cF}{\mathcal{F}}
\newcommand{\cH}{\mathcal{H}}
\newcommand{\cM}{\mathcal{M}}
\newcommand{\cT}{\mathcal{T}}
\newcommand{\cU}{\mathcal{U}}

\newcommand{\CC}{\mathfrak{C}}
\DeclareMathOperator{\cdim}{cd}
\DeclareMathOperator{\vcdim}{vcd}
\DeclareMathOperator{\Loc}{Loc}

\newcommand{\fm}{\mathfrak{m}}
\newcommand{\fM}{\mathfrak{M}}

\newcommand{\Set}{\mathrm{Set}}
\newcommand{\Cat}{\mathrm{Cat}}
\newcommand{\Fun}{\mathrm{Fun}}

\newcommand{\Sp}{\mathrm{Sp}}

\newcommand{\Spc}{\mathrm{Spc}}
\newcommand{\Spec}{\mathrm{Spec}}
\newcommand{\Proj}{\mathrm{Proj}}
\newcommand{\Spf}{\mathrm{Spf}}
\newcommand{\QCoh}{\mathrm{QCoh}}
\newcommand{\IndCoh}{\mathrm{IndCoh}}
\newcommand{\Comod}{\mathrm{Comod}}

\newcommand{\StMod}{\mathrm{StMod}}
\newcommand{\Map}{\mathrm{Map}}

\newcommand{\Open}{\mathrm{Open}}

\newcommand{\Def}{\mathrm{Def}}

\newcommand\noloc{%
  \nobreak
  \mspace{6mu plus 1mu}
  {:}
  \nonscript\mkern-\thinmuskip
  \mathpunct{}
  \mspace{2mu}
}

\makeatletter
\DeclareFontFamily{OMX}{MnSymbolE}{}
\DeclareSymbolFont{MnLargeSymbols}{OMX}{MnSymbolE}{m}{n}
\SetSymbolFont{MnLargeSymbols}{bold}{OMX}{MnSymbolE}{b}{n}
\DeclareFontShape{OMX}{MnSymbolE}{m}{n}{
    <-6>  MnSymbolE5
   <6-7>  MnSymbolE6
   <7-8>  MnSymbolE7
   <8-9>  MnSymbolE8
   <9-10> MnSymbolE9
  <10-12> MnSymbolE10
  <12->   MnSymbolE12
}{}
\DeclareFontShape{OMX}{MnSymbolE}{b}{n}{
    <-6>  MnSymbolE-Bold5
   <6-7>  MnSymbolE-Bold6
   <7-8>  MnSymbolE-Bold7
   <8-9>  MnSymbolE-Bold8
   <9-10> MnSymbolE-Bold9
  <10-12> MnSymbolE-Bold10
  <12->   MnSymbolE-Bold12
}{}

\let\llangle\@undefined
\let\rrangle\@undefined
\DeclareMathDelimiter{\llangle}{\mathopen}%
                     {MnLargeSymbols}{'164}{MnLargeSymbols}{'164}
\DeclareMathDelimiter{\rrangle}{\mathclose}%
                     {MnLargeSymbols}{'171}{MnLargeSymbols}{'171}
\makeatother

\newcommand{\moravamod}[1]{\Mod_{(E_{#1})_*}^{\GG_{#1}} }

\newcommand{\WG}[1]{\WW_n\llangle #1 \rrangle}
\newcommand{\WGT}[1]{\WW_2\llangle #1 \rrangle}

\newcommand{\WFT}[1]{\WW_2\langle #1 \rangle}
\newcommand{\twistgr}[2]{#1\llangle #2 \rrangle}
\newcommand{\EE}{\mathcal{E}}
\newcommand{\EED}{\mathcal{E}^{\cD}}

\newcommand{\EEC}{\mathcal{E}^{\mathcal{C}}}

\newcommand{\GMod}[2]{\Mod_{#1\llangle#2 \rrangle}}

\usepackage{todonotes}

\usepackage{todonotes}

\definecolor{darkspringgreen}{rgb}{0.09, 0.45, 0.27}
\newcommand{\stkout}[1]{\ifmmode\text{\sout{\ensuremath{#1}}}\else\sout{#1}\fi}

\theoremstyle{theorem}
\newtheorem*{chromaticapproach*}{Chromatic Approach}
\newtheorem*{linhypo}{Linearization Hypothesis}
\newtheorem{descript}[theorem]{Description}

\DeclareRobustCommand\circled[1]{\tikz[baseline=(char.base)]{
   \node[shape=circle,draw,inner sep=0pt] (char) {#1};}}

\title{Chromatic structures in stable homotopy theory}
\author{Tobias Barthel and Agn\`es Beaudry}
\date{\today}

\begin{document}

\maketitle

\begin{abstract}
In this survey, we review how the global structure of the stable homotopy category gives rise to the chromatic filtration. We then discuss computational tools used in the study of local chromatic homotopy theory, leading up to recent developments in the field. Along the way, we illustrate the key methods and results with explicit examples. 
\end{abstract}
\tableofcontents

\section{Introduction}

At its core, chromatic homotopy theory provides a natural approach to the computation of the stable homotopy groups of spheres $\pi_*S^0$. Historically, the first few of these groups were computed geometrically through the classification of stably framed manifolds, using the Pontryagin--Thom isomorphism $\pi_*S^0 \cong \Omega_*^{\mathrm{fr}}$. However, beginning with the work of Serre, it soon turned out that algebraic tools were more effective, both for the computation of specific low-degree values as well as for establishing structural results. In particular, Serre proved that $\pi_*S^0$ is a degreewise finitely generated abelian group with $\pi_0S^0 \cong \Z$ and that all higher groups are torsion. 

Serre's method was essentially inductive: starting with the knowledge of the first $n$ groups $\pi_0S^0,\ldots,\pi_{n-1}S^0$, one can in principle compute $\pi_nS^0$. Said differently, Serre worked with the Postnikov filtration of $\pi_*S^0$, in which the $(n+1)$st filtration quotient is given by $\pi_nS^0$. The key insight of chromatic homotopy theory is that $\pi_*S^0$ comes naturally equipped with a completely different filtration---\emph{the chromatic filtration}---which systematically exhibits the large scale symmetries hidden in the stable homotopy category.

Chromatic homotopy theory is the study of the chromatic filtration and the structures that arise from it, both on $\pi_*S^0$ but also on the category of spectra itself. As with many young and active fields, points of views are evolving rapidly and there are few surveys that keep up with the developments. Our goal for this chapter is to present our perspective on the subject and, in the process, to draw one of the possible maps of the field in its current state. 

We would like to emphasize that our exposition is in many ways revisionistic and certainly far from comprehensive, but rather reflects our own understanding of and point of view on the subject. We apologize to those who would have preferred us to present the material from a different point of view, or for us to include important topics we have left untouched. Hopefully, they will take this as a cue to write an expository piece of their own as we feel there is a great need for more background literature in this vibrant field.

\bigskip

In the rest of this short introduction, we give a brief overview of the content of the chapter. 

The goal of Section \ref{sec:landscape} is to introduce and study the chromatic filtration and its consequences from an abstract point of view. More precisely, we will:
\begin{enumerate}[(1)]
	\item Explain that the chromatic filtration arises \emph{canonically} from the global structure of the stable homotopy category. See Section \ref{ssec:toy}.
	\item Describe the geometric origins of the chromatic filtration through the relation with the stack of formal groups. See Section \ref{sec:geomodel}.
	\item Demonstrate that many geometric structures have homotopical manifestations in the chromatic picture that motivate and guide the past and recent developments in the subject. See Section \ref{ssec:chromfiltration1} and Section \ref{ssec:chromfiltration2}.
\end{enumerate}

While Section~\ref{sec:landscape} focuses mostly on the global picture, in Section~\ref{ssec:lcht} we zoom in on $K(n)$-local homotopy theory. In Section~\ref{ssec:lcht}, we introduce Morava $E$-theory $E_n$ and the Morava stabilizer group $\GG_n$, which play a central role in this story because of their relationship to the $K(n)$-local sphere via the equivalence $L_{K(n)}S^0 \simeq E_n^{h\GG_n}$. The resulting descent spectral sequence, whose $E_2$-term is expressed in terms of group cohomology, is one of the most important computational tools in the subject. For this reason, Section~\ref{ssec:gn} is devoted to the study of $\GG_n$ and its homological algebraic properties.

At this point, we go on a hiatus and give an overview of the chromatic story at height $n=1$. This is the content of Section~\ref{sec:cscheight1}, whose goal is to provide the reader with a concrete example to keep in mind for the rest of the chapter. 

The most technical part of this overview of chromatic homotopy theory is Section~\ref{sec:fin}, which presents the theory of finite resolutions. These are finite sequences of spectra that approximate the $K(n)$-local sphere by spectra of the form $E_n^{hF}$ for $F$ finite subgroups of $\GG_n$. The advantage of this approach is that the spectra $E_n^{hF}$ are computationally tractable. Finite resolutions have been one of the most important tools in computations at height $n=2$
and we gives detailed examples in this case in Section~\ref{sec:resheight2}.

In the last part, Section~\ref{sec:thms}, we provide an overview of three topics in chromatic homotopy theory that have seen recent breakthroughs:
\begin{enumerate}[(1)]
\item In Section~\ref{sec:cscmore}, we discuss chromatic reassembly, which describes the passage from the $K(n)$-local to the $p$-local picture. The main open problem is the \emph{chromatic splitting conjecture} and we give an overview of the current state of affairs on this question.
\item In Section~\ref{sec:invdual}, we turn to the problem of computing the group of invertible objects in the symmetric monoidal category of $K(n)$-local spectra. We also touch upon the closely related topic of $K(n)$-local dualities.
\item In Section~\ref{sec:assymptotics} we talk about the asymptotic behavior of local chromatic homotopy theory when $p\to \infty$.
\end{enumerate}
These developments demonstrate how chromatic homotopy theory uncovers structures in the stable homotopy category that reveal the many interactions between homotopy theory and other areas of mathematics.

\subsection*{Conventions and prerequisites}

We will assume that the reader is familiar with basic stable homotopy theory and category theory, as for example contained in the appendices to Ravenel \cite{RavNil}. Throughout this chapter, $\Sp$ will denote a good symmetric monoidal model for the category of spectra, as for example $S$-modules \cite{ekmm}, symmetric spectra \cite{hss_symm}, orthogonal spectra \cite{mmss_orth}, or the $\infty$-category of spectra \cite{ha}. Note that all of these categories model the stable homotopy category, i.e., their associated homotopy categories are equivalent to the stable homotopy category, so the homotopical constructions in this chapter will be model-independent. In fact, Schwede's rigidity theorem \cite{schwede_rigid} justifies that we may work in a model-independent fashion. 

In particular, we freely use the theory of ring spectra in $\Sp$ and module spectra over them, formal groups, and spectral sequences. A full triangulated subcategory of a triangulated category $\cT$ is called thick if it is closed under suspensions and desuspensions, fiber sequences, and retracts. If $\cT$ is cocomplete, then a thick subcategory is called localizing if it closed under all set-indexed direct sums, and we write $\Loc(S)$ for the smallest thick subcategory of $\cT$ containing a given set $S$ of objects in $\cT$. Further, recall that an object $C \in \cT$ is said to be compact (or small) if $\Hom_{\cT}(C,-)$ commutes with arbitrary direct sums in $\cT$; we will write $\cT^{\omega}$ for the full subcategory spanned by the compact objects in $\cT$. If $\cC$ denotes a model (i.e., a stable model category or stable $\infty$-category) for $\cT$, then the corresponding notions for $\cC$ are defined analogously.

\subsection*{Acknowledgements}

We would like to thank Paul Goerss, Hans-Werner Henn, Mike Hopkins and Vesna Stojanoska for clarifying conversations, and are grateful to Dustin Clausen and especially Haynes Miller for useful comments on an earlier version of this document. This material is based upon work supported by Danish National Research Foundation Grant DNRF92, the European Unions Horizon 2020 research and innovation programme under the Marie Sklodowska-Curie grant agreement No. 751794, and the National Science Foundation under Grant No.~DMS-1725563. The authors would like to thank the Isaac Newton Institute for Mathematical Sciences for support and hospitality during the program Homotopy Harnessing Higher Structures when work on this chapter was undertaken. This work was supported by EPSRC Grant Number EP/R014604/1.

\section{A panoramic view of the chromatic landscape}\label{sec:landscape}

The goal of this section is to give an overview of the global structure of the stable homotopy category from the chromatic perspective. Motivated by the analogy with abelian groups and the geometry of the moduli stack of formal groups, we will explain how the solution of the Ravenel Conjectures by Devinatz, Hopkins, Ravenel, and Smith leads to a canonical filtration in stable homotopy theory. The construction as well as the coarse properties of the resulting chromatic filtration are then summarized in the remainder of this section, which prepares for the in-depth study of the local filtration quotients in Section~\ref{ssec:lcht}.  

\begin{remark}
The global point of view taken in this section goes back to Hopkins' original account~\cite{hopkins_global} of his work with Devinatz and Smith on the nilpotence conjectures. It has subsequently led to the study of the global structure of more general tensor-triangulated categories and the systematic development of tt-geometry by Balmer and his coauthors. We refer to Balmer's chapter in this handbook for background and a plethora of further examples.  
\end{remark}

\subsection{From abelian groups to spectra}\label{ssec:toy}

As expressed in Waldhausen's vision of \emph{brave new algebra}, the category $\Sp$ of spectra should be thought of as a homotopical enrichment of the derived category $\cD_{\Z}$ of abelian groups. Consequently, before beginning with our analysis of the global structure of the stable homotopy category, we may consider the case of abelian groups as a toy example. The starting point is the \emph{Hasse square}\index{Hasse square} for the integers, displayed as the pullback square on the left:
\begin{equation}\label{eq:hassesquare}
\xymatrix{\Z \ar[r] \ar[d] & \prod_{p}\Z_p \ar[d] & & M \ar[r] \ar[d] & \prod_pM_p^{\wedge} \ar[d] \\
\Q \ar[r] & \Q \otimes \prod_{p}\Z_p & & \Q \otimes M \ar[r] & \Q \otimes \prod_pM_p^{\wedge}.}
\end{equation}
This is a special case of a local-to-global principle for any chain complex $M \in \cD_{\Z}$, expressed by the homotopy pullback square on the right, in which $M_p^{\wedge}$ denotes the derived $p$-completion of $M$. While the remaining terms in this square seem to be more complicated than $M$ itself, they are often easier from a structural point of view. This is the reason that problems in arithmetic geometry---for example finding integer valued solutions to a set of polynomial equations---can often be divided into two steps: First solve the usually simpler question at individual primes $p$, and then attempt to globalize the solutions. 
 
This approach is tied closely to the global structure of the category $\cD_{\Z}$. Let $\cD_{\Q}$ be the derived category of $\Q$-vector spaces and write $(\cD_{\Z})_p^{\wedge}$ for the category of derived $p$-complete complexes of abelian groups. (Recall that a complex $C$ of abelian groups is derived $p$-complete if it is $p$-local and $\Ext^i(\Q,C) = 0$ for $i=0,1$ or, equivalently, if $C$ is in the image of the zeroth left derived functor of $p$-completion on $\cD_{\Z}$.) We highlight three fundamental properties of these subcategories of $\cD_{\Z}$:
\begin{enumerate}[(1)]
	\item The category $(\cD_{\Z})_p^{\wedge}$ is compactly generated by $\Z/p$. In particular, an object $X \in (\cD_{\Z})_p^{\wedge}$ is trivial if and only if $X \otimes \Z/p$ is trivial. 
	\item The only proper localizing subcategory of $(\cD_{\Z})_p^{\wedge}$ is $(0)$, i.e., if $X$ is any non-trivial object in $(\cD_{\Z})_p^{\wedge}$, then $\Loc(X) = (\cD_{\Z})_p^{\wedge}$, i.e., the smallest full triangulated subcategory of $(\cD_{\Z})_p^{\wedge}$ closed under shifts and colimits which contains $X$ is $(\cD_{\Z})_p^{\wedge}$ itself. 
	\item Any object $M \in \cD_{\Z}$ can be reassembled from its derived $p$-completions $M_p^{\wedge} \in (\cD_{\Z})_p^{\wedge}$, its rationalization $\Q \otimes M \in \cD_{\Q}$, together with the gluing information specified in the pullback square displayed on the right of \eqref{eq:hassesquare}.
\end{enumerate}
Therefore, we may think of $(\cD_{\Z})_p^{\wedge}$ as an irreducible building block of $\cD_{\Z}$. In fact, we can promote these observations to a natural bijection between the residue fields of $\Z$, which are parametrized by the points of $\Spec(\Z)$, and the irreducible subcategories of $\cD_{\Z}$ they detect:
\begin{equation}\label{eq:zbijection}
\begin{Bmatrix}
\text{Prime fields} \\
\Q \text{ and } \F_p \text{ for } p \text{ prime}
\end{Bmatrix} 
\xymatrix@C=1.7pc{ \ar@{<->}[r]^-{\sim} &}
\begin{Bmatrix}
\text{Minimal localizing subcategories}  \\
\cD_{\Q} \text{ and } (\cD_{\Z})_p^{\wedge} \text{ for } p \text{ prime}
\end{Bmatrix}
\end{equation}

A convenient language and framework for describing the global structure of categories like $\cD_{\Z}$ and $\Sp$ is provided by Balmer's \emph{tensor triangular geometry}. Roughly speaking, the Balmer spectrum $\Spc(\cT)$ of a tensor triangulated category $\cT$ has as points the thick $\otimes$-ideal of $\cT^{\omega}$ (where $\cT^{\omega}$ denotes the subcategory of compact objects), equipped with a topology that encodes the inclusions among these subcategories. Whenever $\cT$ is compactly generated by its $\otimes$-unit, as is the case for example for $\Sp$, thick $\otimes$-ideals coincide with thick subcategories of $\cT^{\otimes}$. We refer to Balmer's~\cite{balmer_chapter} for precise definitions and many examples.

With this terminology at hand, we are now ready to make the slogan at the beginning of this section more precise. First note that we can truncate the homotopy groups $S^0$ above degree $0$ to obtain a ring map $\phi\colon S^0 \to \tau_{\le 0}S^0 \simeq H\Z$, which is the Hurewicz map for integral homology.
Base-change along $\phi$ then provides a functor 
\[
\xymatrix{\Sp \simeq \Mod_{S^0}(\Sp) \ar[r]^-{\phi^*} & \Mod_{H\Z}(\Sp) \simeq \cD_{\Z}}
\]
which represents the passage from higher algebra to classical algebra; here, the second equivalence was established by Shipley in \cite{shipley_hz}. Moreover, identifying $\Z \cong [S^0,S^0]$, Balmer constructs a canonical comparison map $\rho$ from the Balmer spectrum of $\Sp$ to the Zariski spectrum of $\Z$. The bijection \eqref{eq:zbijection} implies that the composite
\[
\xymatrix{\Spc(\cD_{\Z}) \ar[rr]^-{\Spc(\phi^*)} & & \Spc(\Sp) \ar[rr]^{\rho} & & \Spec(\Z)}
\]
is an isomorphism, so $\Spc(\Sp)$ contains $\Spec(\Z)$ as a retract. This leads to the following natural question: For $p \in \Spec(\Z)$, what is the fiber $\rho^{-1}(p)$ in $\Spc(\Sp)$? We will see in Theorem~\ref{thm:hopsmith} below that, for each prime ideal $(p) \in \Spec(\Z)$, there is an infinite family of points in $\Spc(\Sp)$ that interpolates between $(p)$ and $(0) \in \Spec(\Z)$, the so-called \emph{chromatic primes}. In other words, the global structure of the stable homotopy category refines the global structure of $\cD_{\Z}$; see \cite[Theorem 1.3.3]{balmer_chapter} for a picture. 

Let $\Sp_{(p)}$ be the category of $p$-local spectra,  i.e., those spectra whose homotopy groups are $p$-local abelian groups. It turns out that $\rho^{-1}(p)$ is determined by $\Spc(\Sp_{(p)})$. We will address the following two problems:
\begin{enumerate}[(1)]
	\item Classify the thick subcategories of $\Sp_{(p)}^{\omega}$.
	\item\label{it2} Find the analogues of prime fields of $\Sp_{(p)}$.  
\end{enumerate}
As we will see, the classification of thick subcategories is a consequence of the answer to Problem \ref{it2}, but before we can get there, we will exhibit a geometric model that serves as a good approximation to stable homotopy theory. 

\begin{convention}\label{conv:plocal}
From here onwards, we fix a prime $p$ and only consider the category of $p$-local spectra. We write $\Sp = \Sp_{(p)}$ and assume without further mention that our spectra have been localized at $p$.
\end{convention}

\subsection{A geometric model for stable homotopy theory}\label{sec:geomodel}

In order to prepare for the resolution of the questions above, we first exhibit a geometric model for the stable homotopy category whose main structural features will turn out to reflect that of $\Sp$ rather closely. 

Recall that the mod $p$ singular cohomology $H^*(X,\F_p)$ of any space or spectrum $X$ is endowed with an action of cohomology operations $\pi_*\Hom(H\F_p,H\F_p)$, which form the mod $p$ Steenrod algebra $\cA_p$. In other words, singular cohomology naturally factors through the functor that forgets the $\cA_p$-module structure and only remembers the underlying $\F_p$-vector space of $H^*(X,\F_p)$:
\[
\xymatrixcolsep{3pc}
\xymatrix{& & \Mod_{\cA_p}^{\mathrm{graded}} \ar[d]^{\text{forget}} \\
\Sp^{\mathrm{op}} \ar[rr]_{H^*(-,\F_p)} \ar@{-->}[rru]^{H^*(-,\F_p)} & & \Mod_{\F_p}^{\mathrm{graded}}}
\]
The Adams spectral sequence, first introduced in \cite{adams_ss}, can then be interpreted as a device that attempts to go back, or at least recover partial information about $X$: There is a spectral sequence
\[
E_2^{s,t} \cong \Ext_{\cA_p}^{s}(H^*(Y),H^*(X))_t \Longrightarrow [X,Y_{p}^{\wedge}]_{t-s},
\] 
which converges whenever $X$ and $Y$ are spectra of finite type with $X$ finite, see for a general study of the convergence properties of (generalized) Adams spectral sequences~\cite{bousfield_locspectra}. Here, finite type means that the mod $p$ cohomology is finitely generated in each degree, and $Y_{p}^{\wedge}$ denotes the $p$-completion of $Y$. The subscript $t$ on the $\Ext$-indicates the internal grading, arising from the grading of cohomology groups involved. Informally speaking, this spectral sequence measures to what extent $\Mod_{\cA_p}$ deviates from being a perfect model for $\Sp$. 

\begin{remark}
Paraphrasing, the \emph{Mahowald uncertainty principle}\index{Mahowald uncertainty principle} asserts that any spectral sequence that computes the stable homotopy groups of a finite spectrum with a machine computable $E_2$-term will be infinitely far from the actual answer. In practical terms this means that the Adams spectral sequence for $X=S^0$ and $Y$ a finite spectrum contains many differentials that require additional input to be determined.
\end{remark}

Building on the work of Novikov~\cite{novikov_mu} and Quillen~\cite{quillen_mu}, Morava~\cite{morava} realized that replacing $H\F_p$ by the Brown--Peterson spectrum $BP$ gives rise to a geometric model for $\Sp$ that resembles its global structure more closely. To describe it, recall that $BP$ is an irreducible additive summand in the $p$-localized complex cobordism spectrum $MU$ with coefficients $BP_* = \Z_{(p)}[v_1,v_2,\ldots]$ and $\deg(v_n) = 2p^n-2$. The generator $v_{i+1}$ is uniquely determined only modulo the ideal $(p,v_1,\ldots,v_{i})$ and there are different choices available, for example the Araki or Hazewinkel generators. See, for example, \cite[A2.2]{ravgreen}. 
The corresponding Hopf algebroid $(BP_*,BP_*BP)$ is a presentation of the moduli stack\index{moduli stack of formal groups} of ($p$-typical) formal groups $\cM_{fg}$ and the category of evenly graded comodules over $(BP_*,BP_*BP)$ is equivalent to the category of quasi-coherent sheaves over $\cM_{fg}$, see for example~\cite{naumann_stacks} for a general treatment. Miller~\cite{miller_left} explains how this equivalence can be extended to all graded comodules by replacing $\cM_{fg}$ by a moduli stack of \emph{spin formal groups}, see also \cite{goerss_mfg}. Taking $BP$-homology induces a functor
\[
\xymatrix{\Sp \ar[r] &\Comod_{BP_*BP}^{\mathrm{even}} \simeq \QCoh(\cM_{fg})  , & X \mapsto BP_*(X),}
\]
where $\Comod_{BP_*BP}^{\mathrm{even}}$ denotes the abelian category of evenly graded $BP_*BP$-comodules. The associated Adams--Novikov spectral sequence\index{Adams--Novikov spectral sequence} has signature
\[
E_2^{s,t} \cong H^s(\cM_{fg};BP_*(X))_t \cong \Ext_{BP_*BP}(BP_*,BP_*(X)) \Longrightarrow \pi_{t-s}X.
\]
The structure of this spectral sequence, whose computational exploitation was a major impetus in the development of chromatic homotopy theory (see \cite{MRW}), is governed by the particularly simple geometric structure of $\cM_{fg}$, which we describe next:

As explained in great detail in \cite{goerss_mfg}, the height filtration of formal groups manifests itself in a descending filtration by closed substacks
\begin{equation}\label{eq:heightfiltration}
\cM_{fg} \supset \cM(1) \supset \cM(2) \supset \cdots 
\end{equation}
where $\cM(n)$ is cut out locally by the ideal defined by the regular sequence $(p,v_1,v_2,\ldots,v_{n-1})$. Note that this filtration is not separated, as the additive formal group has height $\infty$. Write 
\begin{itemize}
\item $\cM_{fg}^{\le n}$ for the open complement of $\cM(n+1)$ representing formal groups of height at most $n$ with $i_n\colon \cM_{fg}^{\le n} \to \cM_{fg}$ the inclusion,
\item $\cH(n) = \cM(n) \cap \cM_{fg}^{\le n}$ for the locally closed substack of formal groups of height exactly $n$, and
\item $\widehat{\cH}(n)$ for its formal completion. 
\end{itemize}
If $\Gamma$ is any formal group of height $n$ over $\F_p$, then $\cH(n)$ is equivalent as a stack to $B\Aut_{\overline{\F}_p}(\Gamma)$, so the filtration quotients of the height filtration \eqref{eq:heightfiltration} contain a single geometric point. Furthermore, there is a (pro-)Galois extension 
\begin{equation}\label{eq:galgeometric}
\Def(\overline{\F}_p,\Gamma) \longrightarrow  \widehat{\cH}(n)
\end{equation}
with Galois group $\Gal(\overline{\F}_p/\F_p) \ltimes \Aut_{\overline{\F}_p}(\Gamma)$, with $\Def(\overline{\F}_p,\Gamma)$ being the Lubin--Tate deformation space. See Remark~\ref{rem:connextensions} below.

In light of \eqref{eq:heightfiltration}, any quasi-coherent sheaf $\cF \in \QCoh(\cM_{fg})$ can be approximated by its restrictions to the open substacks $\cM_{fg}^{\le n}$, so the geometric filtration on $\cM_{fg}$ gives rise to a filtration of $\QCoh(\cM_{fg})$. It follows that the computation of the cohomology of a quasi-coherent sheaf $\cF$ on $\cM_{fg}$ can be restricted to the computation of the cohomology of $\cF$ reduced to the strata $\cH(n)$ together with the gluing data between different strata. The insight of Bousfield, Morava, and Ravenel was that the resulting structure on the $E_2$-term of the Adams--Novikov spectral sequence is in fact manifest in $\pi_*S^0$ and $\Sp$ as well, as we shall see in the next sections. 

\begin{remark}\label{rem:landweber}
An early hint there is such a close relation between $\Sp$ and $\QCoh(\cM_{fg})$ is the Landweber exact functor theorem\index{Landweber exact functor theorem}, which shows that any flat map $f\colon \Spec(R) \to \cM_{fg}$ can be lifted to a complex orientable ring spectrum with formal group classified by $f$. We refer to \cite{behrens_chapter} for more details. 
\end{remark}

\subsection{The chromatic filtration: construction}\label{ssec:chromfiltration1}

The goal of this section is to answer the questions raised at the end of Section~\ref{ssec:toy} and to construct the chromatic filtration. We continue to work in the category $\Sp$ of $p$-local spectra for a fixed prime $p$ as in Convention~\ref{conv:plocal}.

In loose analogy with algebra, a ring spectrum $K \in \Sp$ is said to be a field\index{field object} if every $K$-module splits into a wedge of shifted copies of $K$. In particular, for any spectra $X$ and $Y$, there is a K\"unneth isomorphism
\[
K_*(X \smsh Y) \cong K_*(X) \otimes_{K_*}K_*(Y).
\]
There exists a family of distinct field spectra $K(n)$ for $0 \le n \le \infty$ called the \emph{Morava $K$-theories}\index{Morava $K$-theory}, whose construction will be reviewed in  Section~\ref{sec:moravaEthy}.

As a result of the seminal nilpotence theorem\index{nilpotence theorem} proven by Devinatz, Hopkins, and Smith~\cite{nilpotence1,nilpotence2}, we obtain a classification of fields in $\Sp$.

\begin{theorem}[Hopkins--Smith]\label{thm:nilpotence}
Any field object in $\Sp$ splits (additively) into a wedge of shifted copies of Morava $K$-theories. Moreover, if $R$ is a ring spectrum such that $K(n)_*(R) = 0$ for all $0 \le n \le \infty$, then $R \simeq 0$.
\end{theorem}
\noindent
For example, $K(0) =H\Q$, $K(\infty)=H\F_p$, and $K(1)$ is an Adams summand of mod $p$ $K$-theory. Informally speaking, the spectra $K(n)$ may be thought of as the homotopical residue fields of the sphere spectrum. 

\begin{remark}
As remarked in \cite{nilpotence2}, this theorem can be interpreted as providing a classification of prime fields of $\Sp$. However, there is the subtlety that the ring structure on $K(n)$ is not unique at $p=2$, even in the homotopy category, see \cite[Theorem B.7.4]{RavNil} for a summary and further references. The existence and uniqueness of $\mathbb{A}_{\infty}$-structures on $K(n)$ is studied in Angeltveit's paper~\cite{angeltveit_moravak}. Hopkins and Mahowald have proved that none of these multiplicative structures on $K(n)$ refine to an $\mathbb{E}_2$-ring structure (e.g., \cite[Corollary 5.4]{acb_thom}).
\end{remark}

In light of this theorem, there is a natural notion of support for a spectrum $X \in \Sp$, namely 
\[
\supp(X) = \{n\mid K(n)_*(X) \neq 0\} \subseteq \Z_{\geq 0} \cup \{\infty\}.
\]
This notion of support turns out to be particularly well-behaved for the category of finite spectra $\Sp^{\omega}$. Since $K(\infty)_*F \cong H_*(F,\F_p)= 0$ implies $F \simeq 0$ for finite $F$, for any non-trivial $F$ there exists an $n \in \N$ such that $n \in \supp(F)$. Ravenel \cite{Rav84} further proved that $n \in \supp(F)$ implies $(n+1) \in \supp(F)$, so the only subsets of $\Z_{\geq 0} \cup \{\infty\}$ that can be realized as the support of a finite spectra are the sets $\{n,n+1,n+2,\ldots,\infty\}$ with $n \in \N$. A result of Mitchell's \cite{mitchell_an} implies that all of these subsets can be realized by a finite spectrum.

Write $\cC_0 = \Sp^{\omega}$ and, for $n\ge 1$, let $\cC_n \subseteq \Sp^{\omega}$ be the thick subcategory of finite spectra $F$ with $\supp(F) \subseteq \{n,n+1,n+2,\ldots,\infty\}$ for $n \in \N$. The following consequence of Theorem~\ref{thm:nilpotence} is often called the \emph{thick subcategory theorem}, proven in \cite{nilpotence2}. It says in particular that the support function defined above detects the thick subcategories of $\Sp^{\omega}$:

\begin{theorem}[Hopkins--Smith]\label{thm:hopsmith}\index{thick subcategory theorem}
If $\cC \subseteq \Sp^{\omega}$ is a nonzero thick subcategory, then there exists an $n \ge 0$ such that $\cC  = \cC_n$. Moreover, there is a sequence of proper inclusions
\[
\Sp^{\omega} = \cC_0 \supset \cC_1 \supset \cC_2 \supset \cdots \supset (0),
\]
which completely describes $\Spc(\Sp)$. 
\end{theorem}

This categorical filtration gives rise to a sequence of functorial approximations of any finite spectrum $F$ by spectra that are supported on $\{0,1,\ldots,n\}$ for varying $n$, where the zeroth approximation is given by the rationalization $F\smsh H\Q$. This filtration should be understood as a homotopical incarnation of the geometric filtration of $\cM_{fg}$, so that the approximations of $F$ correspond to the restriction of the associated sheaf $BP_*(F) \in \QCoh(\cM_{fg})$ to $\cM_{fg}^{\le n}$.

The tool required to formulate this notion of approximation rigorously is provided by Bousfield localization\index{Bousfield localization}, which we briefly review here for the convenience of the reader. Let $E$ be a spectrum and consider the full subcategory $\langle E \rangle \subseteq \Sp$ of $E$-acyclic spectra, i.e., those spectra $A$ with $E \smsh A \simeq 0$. Bousfield \cite{bousfield_locspectra} proved that there exists a fiber sequence 
\[
\xymatrix{C_E \ar[r] & \mathrm{id} \ar[r] & L_E}
\]
of functors on $\Sp$ satisfying the following properties:
\begin{enumerate}[(1)]
	\item For any $X \in \Sp$, $C_EX$ is in $\langle E \rangle$.
	\item For any $X \in \Sp$, $L_EX$ is $E$-local, i.e., it does not admit any nonzero maps from an $E$-acyclic spectrum.
\end{enumerate}
It follows formally that $L_EX$ is the initial $E$-local spectrum equipped with a map from $X$, and it is called the $E$-localization of $X$. The full subcategory of $\Sp$ on the $E$-local spectra will be denoted by $\Sp_E$; by construction, it is the quotient of $\Sp$ by $\langle E \rangle$.

In order to extract the part of a spectrum $X$ that is supported on $\{0,1,\ldots,n\}$, i.e., the information of $X$ that is seen by the residue fields $K(0),K(1), \ldots, K(n)$, it is natural to consider the following Bousfield localization
\[
\xymatrix{X \ar[r] & L_nX := L_{K(0) \vee K(1) \vee \ldots \vee K(n)}X.}
\]
In fact, for every $n$ there exists a spectrum $E(n)$ with coefficients $E(n)_* = \Z_{(p)}[v_1,\ldots,v_n][v_n^{-1}]$ called \emph{Johnson--Wilson spectrum}\index{Johnson--Wilson spectrum} (of height $n$) which has the property that $\langle E (n) \rangle = \langle K(0) \vee K(1) \vee \ldots \vee K(n) \rangle$, hence $L_n = L_{E(n)}$. We let $\Sp_n = \Sp_{E(n)}$ denote the category of $E(n)$-local spectra. 

By construction, these localization functors fit into a \emph{chromtic tower}\index{chromatic tower} under $X$ as follows
\begin{equation}\label{eq:chromatictower}
\xymatrixcolsep{1.7pc}
\xymatrix{ & M_nX \ar[d] &  & M_2X \ar[d] & M_1X \ar[d] & M_0X  \simeq H\Q \smsh X \ar@{=}@<-3.5ex>[d] \\ 
\cdots \ar[r] & L_nX \ar[r] & \cdots \ar[r] & L_2X \ar[r] & L_1X \ar[r] & L_0S^0 \simeq  H\Q \smsh X,}
\end{equation}
where the \emph{monochromatic layers}\index{monochromatic layer} $M_nX$ are defined by the fiber sequence 
\[
M_nX  \longrightarrow L_nX  \longrightarrow L_{n-1}X.
\] 
Specializing to the sphere spectrum and applying homotopy groups, we arrive at the definition of the chromatic filtration. 

\begin{definition}\index{chromatic filtration}
The chromatic filtration on $\pi_*S^0$ is given by the descending filtration
\begin{equation}\label{eq:chromaticfiltration}
\pi_*S^0 \supseteq \C_0\pi_*S^0 \supseteq \C_1\pi_*S^0 \supseteq \cdots 
\end{equation}
defined as $\C_n\pi_*S^0 = \ker(\pi_*S^0 \to \pi_*L_nS^0)$.
\end{definition}

There is an important subtlety in the definition of the chromatic filtration, as there is an a priori different way of constructing a filtration of $\Sp$ from the thick subcategory theorem (Theorem~\ref{thm:hopsmith}). Indeed, without relying on the Morava $K$-theories $K(n)$, one may instead take the quotient of $\Sp$ by the localizing subcategories $\Loc(\cC_n) \subseteq \Sp$ for each $n$. The resulting localization functors $L_n^f$ can then be used as above to construct a descending filtration
\[
\pi_*S^0 \supseteq \C_0^f\pi_*S^0 \supseteq \C_1^f\pi_*S^0 \supseteq \cdots 
\]
with $\C_n^f\pi_*S^0 = \ker(\pi_*S^0 \to \pi_*L_n^fS^0)$, known as the \emph{geometric filtration}\index{geometric filtration}, see Ravenel \cite[Section 7.5]{RavNil}. If $X$ is spectrum such that $L_n^fX \simeq 0$, then also $L_nX \simeq 0$, so there are natural comparison transformations $L_n^f \to L_n$, leading to the following optimistic conjecture about the comparison between the two filtrations:

\begin{conjecture}[Telescope conjecture]\index{telescope conjecture}
The natural transformation $L_n^f \to L_n$ is an equivalence.
\end{conjecture}

A number of equivalent formulations of this conjecture and the current state of knowledge about it can be found in Mahowald--Ravenel--Schick \cite{mrstriple} and \cite{chromconj_survey}. The smash product theorem\index{smash product theorem} of Hopkins and Ravenel \cite[Section 8]{RavNil} states that $L_n$ is smashing, i.e., $L_n$ as an endofunctor on $\Sp$ commutes with colimits, while the analogous fact for $L_n^f$ was proven by Miller \cite{miller_finite}. It therefore suffices to show the telescope conjecture for $S^0$. This has been verified by explicit computations for $n=0$ and $n=1$ by work of Mahowald \cite{mahowald_bo} for $p=2$ and Miller \cite{miller_ass} for odd $p$, but the telescope conjecture is open in all other cases. It is known however that $L_n^fM \to L_nM$ is an equivalence for many spectra $M$, including $BP$-modules~\cite[Corollary 1.10]{cschov} and $E(m)$-local spectra~\cite[Corollary 6.10]{HovStrick} for any $m \ge 0$.

\subsection{The chromatic filtration: disassembly and reassembly}\label{ssec:chromfiltration2}

The goal of this subsection is to first demonstrate how the chromatic filtration decomposes the stable homotopy groups of spheres into periodic families and then to explain how these irreducible pieces reassemble into $\pi_*S^0$. The starting point is the chromatic convergence theorem due to Hopkins and Ravenel, proven in  \cite{RavNil}, whose content  is that the chromatic tower \eqref{eq:chromatictower} does not lose any information about $S^0$. In particular, the chromatic filtration \eqref{eq:chromaticfiltration} on $\pi_*S^0$ is exhaustive. We continue to follow Convention~\ref{conv:plocal}.
\begin{theorem}[Hopkins--Ravenel]\label{thm:chromaticconvergence}\index{chromatic convergence theorem}
The canonical map $X \to \lim_nL_nX$ is an equivalence for all finite spectra $X$.
\end{theorem}
\begin{remark}
For general $X$, this map can be far from being an equivalence. For example, the chromatic tower of $H\F_p$ or the Brown--Comenetz dual $IS^0$ of the sphere is identically zero. However, chromatic convergence is known to hold for a class of spectra larger than just finite ones, including $\mathbb{CP}^{\infty}$. See \cite{chromaticcompletion}.
\end{remark}

We now turn to the filtration quotients of the chromatic filtration, which correspond homotopically to the monochromatic layers $M_nX$. Much of the material in this section can be found in \cite{HovStrick}.

The layers $M_nX$ decompose into spectra which are periodic of periods a multiple of $2(p^n-1)$, thereby resembling the decomposition of light into waves of different frequencies. (This is the origin of the term \emph{chromatic} homotopy theory, coined by Ravenel.) More precisely, if $X$ is any spectrum, then its $n$th monochromatic layer is equivalent to a filtered colimit of spectra $F_{\alpha}$,
\[
\xymatrix{\colim_{\alpha}F_{\alpha} \ar[r]^-{\sim} & M_nX,}
\]
such that each $F_{\alpha}$ is periodic. That is, for each $\alpha$ there exists a natural number $\lambda(\alpha)$ and a homotopy equivalence $F_{\alpha} \simeq \Sigma^{2(p^n-1)p^{\lambda(\alpha)}}F_{\alpha}$. This follows from the fact that $M_n$ is equivalent to the colocalization of the $E(n)$-local category with respect to the $E(n)$-localization of a finite type $n$ spectrum, see for example \cite[Proposition~7.10]{HovStrick}, together with the periodicity theorem of Hopkins and Smith~\cite[Theorem 9]{nilpotence2}

Having resolved $S^0$ into its irreducible chromatic pieces $M_nS^0$, it is now time to consider the question of how to reassemble the pieces. For this, it is more convenient to consider the $K(n)$-localizations instead of the monochromatic layers, as we shall explain next. 

Write $\fM_n \subset \Sp$ for the essential image of the functor $M_n$ and let $\Sp_{K(n)}$ be the category of $K(n)$-local spectra. For any $n$, the functors $L_{K(n)}$ and $M_n$ restrict to an adjunction on the category $\Sp_n$ (with $M_n$ as the left adjoint) and then further to a symmetric monoidal equivalence~\cite[Theorem~6.19]{HovStrick}
\[
\xymatrix{M_n   \colon  \Sp_{K(n)} \ar@<0.5ex>[r] &  \ar@<0.5ex>[l] \fM_n  \noloc L_{K(n)}.   }
\]
So we may equivalently work with $L_{K(n)}S^0$ in place of $M_nS^0$. 

\begin{remark}
The more categorically minded reader may think of the situation as follows: The descending filtration of $\Sp^{\omega}$ of Theorem \ref{thm:hopsmith} extends to two descending filtrations of $\Sp$:
\[
\Sp = \ker(0) \supset \ker(L_0) \supset \ker(L_1) \supset \cdots \supset \ker(\id) = (0)
\]
and
\[
\Sp = \Loc(\cC_0) \supset \Loc(\cC_1) \supset \Loc(\cC_2) \supset \cdots \supset (0),
\]
which are equivalent if the telescope conjecture holds for all $n$. Focusing on the first filtration for concreteness and writing $\Sp_n$ for the essential image of $L_n$ as before, we could equivalently pass to the associated ascending filtration
\[
(0) = \im(0)  \subset \Sp_0 \subset \Sp_1 \subset \Sp_2 \subset \cdots \subset \im(\id) = \Sp.
\]
The consecutive subquotients $\Sp_{n}/\Sp_{n-1}$ can then be realized in two different ways as subcategories of $\Sp$, namely either as a localizing subcategory $\fM_n$ or as a colocalizing subcategory $\Sp_{K(n)}$. The resulting equivalence between $\fM_n$ and $\Sp_{K(n)}$ is an instance of a phenomenon called \emph{local duality}\index{local duality}, see \cite{bhv1}. 
\end{remark}

Suppose $X$ is a spectrum for which we have determined $L_{n-1}X$ and $L_{K(n)}X$, and we are interested in reassembling them to obtain $L_nX$. Motivated by the geometric model of Section~\ref{sec:geomodel}, we expect this process to be analogous to the way a sheaf on the open subset $\cM_{fg}^{\le n-1}$ and another sheaf on the stratum $\cH_n$ are glued together along the formal neighborhood $\widehat{\cH}_n$ of $\cH_n$ inside $\cM_{fg}^{\le n-1}$ to produce a sheaf on $\cM_{fg}^{\le n}$. This picture turns out to be faithfully reflected in stable homotopy theory: The chromatic reassembly process for $X \in \Sp$ is governed by the homotopy pullback square displayed on the left, usually called the \emph{chromatic fracture square}\index{chromatic fracture square} (see for example \cite{greenlees_axiomatictate}):
\begin{equation}\label{eq:chromaticfracture}
\xymatrix{L_{n}X \ar[r] \ar[d] & L_{K(n)} X \ar[d] & \Sp_n \ar[r]^-{L_{K(n)}} \ar[d]_{X \mapsto \iota_n(X)} & \Sp_{K(n)} \ar[d]^{L_{n-1}} \\
L_{n-1}X \ar[r]_-{\iota_n(X)} & L_{n-1}L_{K(n)}X & \Fun(\Delta^1,\Sp_{n-1}) \ar[r]_-{\mathrm{target}} & \Sp_{n-1}}
\end{equation}
In fact, by~\cite{acb_cubes} the category $\Sp_n$ itself admits a decomposition into chromatically simpler pieces, see the pullback square on the right of \eqref{eq:chromaticfracture}. Here, $\Fun(\Delta^1,\Sp_{n-1})$ is the arrow category of $\Sp_{n-1}$ and the pullback is taken in a suitably derived sense. The labels of the arrows in this diagram indicate how to translate from the chromatic fracture square of a spectrum $X$ to the categorical decomposition on the right of \eqref{eq:chromaticfracture}.

Based on computations of Shimomura--Yabe \cite{shimyabe}, Hopkins~\cite{cschov} conjectured that the chromatic reassembly process which recovers $L_{n}X$ from $L_{K(n)}X$ and $L_{n-1}X$ takes a particularly simple form:\index{chromatic splitting conjecture}

\begin{conjecture}[Weak Chromatic Splitting]\label{conj:wcsc}
The map 
\[
\xymatrix{\iota_n(S_p^0)\colon L_{n-1}S_p^0 \ar[r] & L_{n-1}L_{K(n)}S_p^0}
\]
in \eqref{eq:chromaticfracture} is split, i.e., it admits a section. Here, $S_p^0$ is the $p$-complete sphere spectrum.
\end{conjecture}

This conjecture, its variations, and its consequences are discussed in more detail in Section~\ref{sec:cscmore}. For now we note that Conjecture~\ref{conj:wcsc} is known to hold for $n\le 2$ and all primes $p$, and is wide open otherwise.\\

We can now summarize the chromatic approach\index{chromatic approach} as follows:
\begin{chromaticapproach*}
The chromatic approach to $\pi_*S^0_{(p)}$ consists of three steps:
\begin{enumerate}[(1)]
	\item Compute $\pi_*L_{K(n)}S^0$ for each $n$.
	\item Understand the gluing in the chromatic fracture square \eqref{eq:chromaticfracture}.
	\item Use chromatic convergence (Theorem~\ref{thm:chromaticconvergence}) to recover $S^0_{(p)}$.
\end{enumerate}
\end{chromaticapproach*}

Finally, the $p$-local sphere spectrum $S^0_{(p)}$ determines $S_p^0$ by $p$-completion. Together with $H\Q \smsh S^0 \simeq H\Q$, we can thus reassemble the sphere spectrum $S^0$ itself via the following homotopical analogue of the Hasse square \eqref{eq:hassesquare}:
\[
\xymatrix{S^0 \ar[r] \ar[d] & \prod_pS_p^0 \ar[d] \\
H\Q \smsh S^0 \ar[r] & H\Q \smsh \prod_pS_p^0.}
\]
In the next section, we discuss the first two steps of the chromatic approach.

\begin{remark}
As mentioned earlier, the deconstructive analysis of the stable homotopy category based on its spectrum $\Spc(\Sp)$ can be carried out in any tensor triangulated category; many examples can be found in~\cite{balmer_chapter}. This is the subject of \emph{prismatic algebra}. An especially interesting example is the stable module category $\StMod_{kG}$ of a finite $p$-group $G$ and field $k$ of characteristic $p$, whose spectrum $\Spc(\StMod_{kG})$ is homeomorphic to $\Proj(H^*(G;k))$, the Proj construction of the graded ring $H^*(G;k)$. This category is a good test case for chromatic questions: for instance, the analogues of both the telescope conjecture and the weak chromatic splitting conjecture are known to hold in $\StMod_{kG}$, see \cite{bik_fingps} and \cite{bhv3}. 
\end{remark}

\section{Local chromatic homotopy theory}\label{ssec:lcht}

We begin this section by introducing  the main players of local chromatic homotopy theory: Morava $E$-theory $E_n$, the Morava stabilizer group $\GG_n$ and its action on $E_n$, and the resulting descent spectral sequence computing $\pi_*L_{K(n)}S^0$. We then summarize the key algebraic features of the Morava stabilizer group, its continuous cohomology, and its action on the coefficients of Morava $E$-theory. In order to have a toy case in mind for the general constructions to follow, we study in detail the case of height $1$. 

\subsection{Morava $E$-theory and the descent spectral sequence}\label{sec:moravaEthy}

The chromatic program has led us naturally and inevitably to the study of the $K(n)$-local categories, which should be thought of as an analog of $(\cD_{\Z})_p^{\wedge}$ for abelian groups. Formally, we note that $\Sp_{K(n)}$ is a closed symmetric monoidal stable category. Moreover, in close analogy with Section~\ref{ssec:toy}, the $K(n)$-local categories have the following properties:
\begin{enumerate}[(1)]
	\item The category $\Sp_{K(n)}$ is compactly generated by $L_nF(n)$ for any $F(n) \in \cC_n \setminus \cC_{n+1}$ for $\cC_n$ as in Theorem~\ref{thm:hopsmith}, and an object $X \in \Sp_{K(n)}$ is trivial if and only if $X \smsh K(n)$ is trivial.
	\item The only proper localizing subcategory of $\Sp_{K(n)}$ is $(0)$.
	\item A spectrum $X \in \Sp_n$ can be reassembled from $L_{K(n)}X$, $L_{n-1}X$, together with the gluing information specified in the pullback square displayed on the right of \eqref{eq:chromaticfracture}.
\end{enumerate}
This confirms the idea that the $K(n)$-local categories play the role of the irreducible pieces of $\Sp$. With the techniques developed so far, both the finer structural properties of $\Sp_{K(n)}$ as well as any concrete calculations would be essentially inapproachable: Incipit \emph{Morava $E$-theory}.

We let $\Gamma_n$ denote the Honda formal group law\index{Honda formal group law} of height $n$. It is the formal group law classified by the map 
\[
BP_*  \cong \Z_{(p)}[v_1, v_2, v_3, \ldots ] \longrightarrow \F_p
\] 
which sends $v_n$ to $1$ and $(p, v_1, \ldots, v_{n-1}, v_{n+1}, v_{n+2}, \ldots )$ to zero. In fact, it is the unique $p$-typical formal group law over $\F_{p^n}$ whose $p$-series satisfies $[p]_{\Gamma_n}(x) =x^{p^n}$. A good reference on formal group laws for homotopy theorists is \cite[Appendix A2]{ravgreen}.

Let $\WW_n= W(\F_{p^n}) $ be the ring of Witt vectors of $\F_{p^n}$\index{Witt vectors}, which is isomorphic to the ring of integers in an unramified extension of $\mathbb{Q}_p$ of degree $n$. Lubin and Tate~\cite{lubintate} showed that there exists a $p$-typical universal deformation $F_{n}$ of $\Gamma_n$ to complete local rings with residue field $\FF_{p^n}$, whose formal group law $F_{n}(x,y)$ is defined over the ring 
\begin{equation}\label{eq:Estariso}
(E_n)_0 = \WW_n[\![u_1, \ldots, u_{n-1}]\!]   \ \ \ u_i \in (E_n)_0.   
\end{equation}
Introducing a formal variable $u$ in degree $-2$ then allows to extend $F_{n}$ to a graded formal group law $F_{n, \mathrm{gr}}(x,y) = uF_{n}(u^{-1}x,u^{-1}y)$ defined over $(E_n)_* = (E_n)_0[u^{\pm1}]$, classified by the ring homomorphism  
\[ 
BP_* \cong \Z_{(p)}[v_1, v_2, v_3, \ldots ]  \longrightarrow (E_n)_*
\] 
which sends $(v_{n+1}, v_{n+2}, \ldots)$ to zero, $v_n$ to $u^{1-p^n}$, and $v_k$ to $u_k u^{1-p^k}$ for $k<n$. Here, we are using the Araki generators for $BP_*$. See \cite[A2.2]{ravgreen} for more details.

In order to lift this construction to stable homotopy theory, one first shows that the functor
\[
X \mapsto (E_n)_* \otimes_{BP_*} BP_*(X)
\] 
is a homology theory, represented by a complex orientable ring spectrum $E_n = E(\F_{p^n}, \Gamma_n)$, called Morava $E$-theory\index{Morava $E$-theory} or the \emph{Lubin--Tate spectrum}\index{Lubin--Tate spectrum} because of its connection to Lubin--Tate deformation theory, see Rezk \cite{rezk_hm}. This is an instance of the Landweber exact functor theorem\index{Landweber exact functor theorem} mentioned in Remark~\ref{rem:landweber}. The spectrum $E_n$ is a completed and  2-periodized version of the Johnson--Wilson spectrum $E(n)$ from Section \ref{ssec:chromfiltration1} and it turns out that $L_{E(n)} = L_{E_n}$ for all $n$; in particular, the terms $E(n)$-local and $E_n$-local are synonymous.

Since $(E_n)_*$ is a regular graded commutative ring which is concentrated in even degrees, reduction modulo the maximal ideal $\fm = (p,u_1,\ldots,u_{n-1})$ can be realized by a (homotopy) ring map
\[
E_n \longrightarrow E_n/\fm =: K_n
\] 
with $\pi_*K_n \cong \F_{p^n}[u^{\pm 1}]$, see for example Chapter V of \cite{ekmm}. The spectrum $K_n$ splits as a wedge of equivalent spectra, which are shifts of the Morava $K$-theory $K(n)$ of Theorem~\ref{thm:nilpotence}, with homotopy groups $K(n)_* \cong \F_p [v_n^{\pm 1}]$ for $v_n = u^{1-p^n}$.

\begin{definition}\label{defn:stabilizergrop}
The \emph{small Morava stabilizer group} $\mathbb{S}_n:=\Aut_{\mathbb{F}_{p^n}} (\Gamma_n)$ is the group of automorphisms of $\Gamma_n$ with coefficients in $\F_{p^n}$ 
\[
\bbS_n = \{f(x) \in \F_{p^n}[\![x]\!] : f(\Gamma_n(x,y)) = \Gamma_n(f(x),f(y)), \ f'(0) \neq 0\}.
\]
Since $\Gamma_n$ is defined over $\F_p$, the Galois group $ \Gal=\Gal(\F_{p^n}/\F_p)$ acts on $\mathbb{S}_n$ by acting on the coefficients of an automorphism. The \emph{big Morava Stabilizer group} $\GG_n $ is the extension $\mathbb{S}_n \rtimes \Gal $. Equivalently, $\GG_n$ is the group of automorphisms of the pair $(\F_{p^n}, \Gamma_n)$.\index{Morava stabilizer group}
\end{definition}

The construction $E(\F_{p^n}, \Gamma_n)$ is natural in the formal group law $\Gamma_n$, so there is an up to homotopy action of $\Aut_{\mathbb{F}_{p^n}} (\Gamma_n)$ on $E(\F_{p^n}, \Gamma_n)$. This action can be promoted to an action through $\mathbb{E}_{\infty}$-ring maps in a unique way: By Goerss--Hopkins--Miller obstruction theory~\cite{hm_elliptic, gh_modulien}, $E_n$ admits an essentially unique structure of an $\mathbb{E}_{\infty}$-ring spectrum and $\GG_n$ acts on it through $\mathbb{E}_{\infty}$-ring maps. In fact, $\GG_n$ gives essentially all such automorphisms of $E_n$. A new proof of these results from the perspective of derived algebraic geometry has recently appeared in Lurie~\cite{ell2}. The connection between $K(n)$-local homotopy theory and Morava $E$-theory is then illustrated in the diagram
\begin{align}\label{eq:galoisspectra}
\xymatrix{L_{K(n)}S^0 \ar[r] & E_n \ar[r] & \wK.}
\end{align}
The first map is a pro-Galois extension of ring spectra\index{pro-Galois extension of ring spectra} with Galois group $\GG_n$ in the sense of Rognes. In particular, $L_{K(n)}S^0  \simeq E_n^{h\GG_n}$ and the extension $L_{K(n)}S^0 \to E_n$  behaves like an unramified field extension. The second map in \eqref{eq:galoisspectra} corresponds to the passage to the residue field. 
See \cite{rognes,br_absgal} for precise definitions on pro-Galois extensions and \cite{DH} for a definition of homotopy fixed points for profinite groups. Further results and alternative approaches to the construction of (continuous) homotopy fixed points in the generality needed for chromatic homotopy theory can be found in \cite{davis_hfp, quick_chfp, davisquick} and the references therein.

\begin{remark}
Note also that the extension $L_{K(n)}S^0 \to E_n $ can be broken into two pro-Galois extensions
\[
\xymatrix{L_{K(n)}S^0 \ar[rr]^-{\abGal} & &  E_n^{h\wS} \simeq L_{K(n)}S^0(\omega)  \ar[rr]^-{\wS} & & E_n,}
\]
where the arrows are labelled by the structure group of the extension. Here, $L_{K(n)}S^0(\omega) $ is an $\mathbb{E}_{\infty}$-ring obtained by adjoining a primitive $(p^n-1)$th root of unity $\omega$ to the $K(n)$-local sphere. See \cite[5.4.6]{rognes} and \cite[Section 1.6]{BobkovaGoerss} for details on this. 
\end{remark}

From the fact that the first map of \eqref{eq:galoisspectra} is a Galois extension, it follows that
\begin{align}\label{eq:EnEn}
\pi_*L_{K(n)}(E_n \smsh E_n) &\cong \Map^c(\GG_n, (E_n)_*) ,
\end{align}
where $\Map^c$ denotes the continuous functions as profinite sets. See for example \cite[Theorem 4.11]{hovey_ops}. In fact, the functor $(E_n)_*^{\vee}(-):=\pi_*L_{K(n)}(E_n \smsh - )$ 
takes values in a category of \emph{Morava modules}, which are $(E_n)_*$-modules equipped with a continuous action by $\GG_n$ (see Definition~\ref{defn:moravamodules}). Furthermore, a map $f$ is a $K(n)$-local equivalence (i.e., $K(n)_*(f)$ is an isomorphism) if and only if $(E_n)_*^{\vee}(f)$ is an isomorphism. The resulting relationship between the topological category $\Sp_{K(n)}$ and the algebraic category of Morava modules provides very powerful tools for computations in the $K(n)$-local category. In particular, it gives rise to a \emph{homotopy fixed point spectral sequence}\index{homotopy fixed point spectral sequence}, also called the \emph{descent spectral sequence}\index{descent spectral sequence}. 

\begin{theorem}[Hopkins--Ravenel \cite{RavNil}, Devinatz--Hopkins \cite{DH}, Rognes \cite{rognes}]\label{thm:descentss}
The unit map $L_{K(n)}S^0 \to E_n$ is a pro-Galois extension with Galois group $\GG_n$. There is a convergent descent spectral sequence 
\begin{equation}\label{eq:hfpss}
E_2^{s,t} \cong H^s_c(\GG_n, (E_n)_t)  \Longrightarrow \pi_{t-s} L_{K(n)}S^0,
\end{equation}
which collapses with a horizontal vanishing on a finite page. 
\end{theorem}

The spectral sequence \eqref{eq:hfpss} is the \emph{$K(n)$-local $E_n$-based Adams--Novikov spectral sequence}\index{$K(n)$-local $E_n$-based Adams--Novikov spectral sequence}, which for a general $X$ has the form
\begin{align}\label{rem:KnlocalANSS}
E_1^{s,t} = \pi_tL_{K(n)}(E_n \smsh E_n^{s} \smsh X)  \Longrightarrow \pi_{t-s} L_{K(n)}X. 
\end{align}
It is constructed in \cite[Appendix A]{DH}. The description of the $E_2$-page in terms of continuous group cohomology $H^*_c$ for $X=S^0$ uses \eqref{eq:EnEn} to identify the $E_1$-term with the cobar complex. More generally, if the $(E_n)_*$-module $(E_n)_*^{\vee}(X)$ is flat, or finitely generated, or if there exists $k\ge 1$ such that $\fm^k[(E_n)_*^{\vee}(X)] = 0$, then there is an isomorphism~\cite{bh_e2adams}
\[
E_2^{s,t} \cong H^s_c(\GG_n, (E_n)_t^{\vee}(X)).
\]  
Section~\ref{sssec:moravahomalg} below further discusses homological algebra over profinite groups and properties of this spectral sequence.

In fact, as discussed in \cite{mathew_galois}, the $\GG_n$-action on $E_n$ lifts to an action on the $\infty$-category $\Mod_{E_n}$, which yields a categorical reformulation of the theorem as a canonical equivalence
\[
\xymatrix{\Sp_{K(n)} \ar[r]^-{\sim} & \Mod_{E_n}^{h\GG_n}, }
\]
where the right hand side denotes the homotopy fixed points taken in the $\infty$-category of $\infty$-categories. These observations demonstrate the fundamental role of $E_n$-theory and the Morava stabilizer group $\GG_n$ together with its cohomology in chromatic homotopy theory.

\begin{remark}\label{rem:Ethygeneral}\label{rem:connextensions}
Other choices for Morava $E$-theory are possible. For any perfect field $k$ of characteristic $p$ and formal group law $\Gamma$ of height $n$ defined over $k$, there is an associated spectrum $E(k, \Gamma)$ whose formal group law is a universal deformation of $\Gamma$ to complete local rings with residue field isomorphic to $k$. There is an associated Morava $K$-theory $K(k, \Gamma)$, stabilizer group $\GG(k,\Gamma)$, and $\GG(k,\Gamma)$-Galois extension $L_{K(k, \Gamma)}S^0 \to E(k, \Gamma)$. The localization functor $L_{(k,\Gamma)}$ is independent of the choice of formal group law $\Gamma$ and extension $k$ of $\F_p$, so one can make any convenient choice to study $L_{K(k, \Gamma)}S^0$.

Recall the Galois extension $\Def(\overline{\F}_p,\Gamma) \longrightarrow  \widehat{\cH}(n)$ from \eqref{eq:galgeometric}. By definition, $\GG(\bF_p, \Gamma)$ is the group $ \Aut_{\overline{\F}_p}(\Gamma) \rtimes \Gal(\overline{\F}_p/\F_p) $, and the Galois extension $ E(\bF_p, \Gamma)  \leftarrow L_{K(\bF_p, \Gamma)}S^0 $ is a homotopical lift of the pro-Galois extension \eqref{eq:galgeometric}. The coefficients of the Lubin--Tate spectrum $E(\bF_p, \Gamma)$ correspond to the global sections of $\Def(\overline{\F}_p,\Gamma)$ (hence the reversal of the arrow direction). A more thorough discussion of Morava $E$-theory\index{Morava $E$-theory} is given in \cite{stapleton_chapter}.  See also Remark~\ref{rem:choiceGamma} for more on this point. 
\end{remark}

The first step in the Chromatic Approach described in Section~\ref{ssec:chromfiltration2} is to compute the homotopy groups of $L_{K(n)}S^0 \simeq E_n^{h\GG_n}$. As for any Galois extension, it makes sense to first study intermediate extensions. In general, if $H$ and $K$ are closed subgroups of $\GG_n$ and $H$ is normal in $K$, then $E_n^{hK} \to E_n^{hH}$ is a $K/H$-Galois extension 
and there is a descent spectral sequence
\begin{equation}\label{eq:hfpssgen}
E_2^{s,t} \cong H^s_c(K/H, (E_n^{hH})_t)  \Longrightarrow \pi_{t-s} E_n^{hK}.
\end{equation}
See for example Devinatz \cite{devinatz_lhs}. It seems natural to consider the following kinds of intermediate extensions:
\begin{enumerate}[(a)]
\item The $\GG_n/K$ Galois extensions $L_{K(n)}S^0 \to E_n^{hK}$ for $K \subseteq \GG_n$ normal closed subgroups.
\item The $F$-Galois extensions $E_n^{hF} \to E_n$ for finite subgroups $F$ of $\GG_n$.
\end{enumerate}

An important example of an intermediate extension of the form (a) is given in Remark~\ref{rem:3termfiber} below. These kinds of extensions are conceptually important, but the homotopy groups of spectra of the form $E_n^{hK}$ are generally out of reach at heights $n\geq 3$. An exception is when $K=F \subseteq \GG_n$ is finite, which brings us to extensions of the form (b), in which case the intermediate extensions $E_n^{hF}\to E_n$ and computations of the homotopy groups of $E_n^{hF}$ are more accessible. 

In fact, there are many computations of the homotopy groups of $E_n^{hF}$ at various heights and the recent developments in equivariant homotopy theory by Hill, Hopkins and Ravenel \cite{HHR}, followed by the work on real orientations for $E$-theory of Hahn and Shi \cite{hahnshi} make these computations even more accessible. A non-exhaustive list of reference related to these types of computations is given by \cite{bauer_tmf, beaudrybobkovahillstojanoska, BehrensOrmsby, BobkovaGoerss, hahnshi, heard_tate, HHRC4, c4e4, hill_chapter, MR}.

In view of this, an approach to studying the $K(n)$-local sphere is to approximate it by spectra of the form $E_n^{hF}$ for finite subgroups $F \subseteq \GG_n$. These approximations fit together to form so-called \emph{finite resolutions}. This is the philosophy established by Goerss, Henn, Mahowald, and Rezk (GHMR) in \cite{ghmr}. It has proven to be very effective for organizing computations and clarifying the structure of the $K(2)$-local category. In the next sections, we will describe the study of chromatic homotopy theory using the finite resolution perspective, starting with explicit examples at height $n=1$. In particular, finite resolutions will be discussed at length in Section~\ref{sec:fin}.

\subsection{The Morava stabilizer group}\label{ssec:gn}\label{ssec:finsubgps}

In this section, we give more details on the structure of the Morava stabilizer groups\index{Morava stabilizer group} $\GG_n$ and $\mathbb{S}_n$, which were introduced in Definition~\ref{defn:stabilizergrop}. We also discuss homological algebra in this context. More detail on this material can be found, for example, in \cite{henn_minicourse}.

\subsubsection{The structure of $\GG_n$}\label{sec:GGn}

Recall that $\Gamma_n$ denotes the Honda formal group law\index{Honda formal group law} of height $n$. We write 
\[x+_{\Gamma_n}y := \Gamma_n(x,y).\]
By definition, $\mathbb{S}_n$ is the group of the units in $\End_{\F_{p^n}}(\Gamma_n)$. In fact, 
\[\End_{\bF_p}(\Gamma_n)  \cong \End_{\F_{p^n}}(\Gamma_n).\] 
That is, all endomorphisms of $\Gamma_n$ have coefficients in $\F_{p^n}$. We give a brief description of the endomorphism ring here, originally due to Dieudonn{\'e}~\cite{dieudonne_lie7} and Lubin~\cite{lubin_oneparam}. A good reference for this material from the perspective of homotopy theory is \cite[Appendix A2.2]{ravgreen}.

Recall that $\WW_n$ denotes the ring of Witt vectors on $\F_{p^n}$.
It is isomorphic to the ring of integers $ \Z_p(\omega)$ of the unramified extension of $\mathbb{Q}_p(\omega)$ obtained from $\Q_p$ by adjoining a primitive $(p^n-1)$th root of unity $\omega$. The residue field is $\F_{p^n}$ and we also let $\omega$ denote its reduction in $\F_{p^n}$. 

The series $\omega(x) = \omega x$ and $\xi(x) = x^p$ are elements of $\End_{\F_{p^n}}(\Gamma_n)$ and, in fact, the endomorphism ring is a $\Z_p$-module generated by these elements:
\begin{equation}\label{eq:stabpres}
\End_{\F_{p^n}}(\Gamma_n) \cong \WW_n\langle \xi \rangle/ ( \xi \omega  - \omega^{p}\xi,  \xi^n- p)
\end{equation}
The identification of \eqref{eq:stabpres} is given explicitly as follows. An element of the right hand side can be written uniquely as
\begin{equation}\label{eq:alphaexpressions}
f= \sum_{j=0}^{n-1}f_j \xi^j 
\end{equation}
for $f_j \in \WW_n$. Further, $f_j = \sum_{i=0}^{\infty}a_{j+in}p^i$ for unique elements $a_i \in \WW_n$ such that $a_i^{p^n}-a_i=0$. Using the fact that $\xi^n=p$, the element $f$ can also be written uniquely as 
\[
f = \sum_{i =0}^{\infty} a_i \xi^i.
\]  
The series
\[
[f](x) = \sum\nolimits_{i\geq 0}^{\Gamma_n} a_ix^{p^i} = a_0 x +_{\Gamma_n} a_1 x^p +_{\Gamma_n} \cdots
\]
is the endomorphism of $\Gamma_n$ corresponding to $f$. 

Let $\mathbb{D}_n = \Q \otimes \End_{\F_{p^n}}(\Gamma_n) $. There is a valuation $v \colon \D_n^{\times} \to \frac{1}{n}\Z$ normalized so that $v(\xi)=1/n$. The ring $\mathbb{D}_n$ is a central division algebra algebra over $\mathbb{Q}_p$ of Hasse invariant $1/n$. The ring of integers of $\mathcal{O}_{\D_n}$ is defined to be those $x\in \D_n$ such that $v(x)\geq 0$, so that $ \mathcal{O}_{\mathbb{D}_n} \cong \End_{\F_{p^n}}(\Gamma_n)$.

The element $\xi$ is invertible in $\mathbb{D}_n$ and conjugation by $\xi$ preserves $ \mathcal{O}_{\mathbb{D}_n}$. In fact, conjugation by $\xi$ corresponds to the action of a generator of $\Gal = \Gal(\F_{p^n}/\F_p)$ on $\mathbb{S}_n \cong \mathcal{O}_{\mathbb{D}_n}^{\times}$. From this, we get a presentation
\[  
\mathbb{D}_n^{\times}/(\xi^n) \cong \GG_n. 
\]

The problem of determining the isomorphism types and conjugacy classes of maximal finite subgroups of $\mathbb{S}_n$ was studied by Hewett \cite{hewett_groups, hewett_normalizers} and was revisited by Bujard \cite{bujard}. We have listed the conjugacy classes of maximal finite subgroups of $\mathbb{S}_n$ in Table~\ref{tab:finitesub}. Note that the list is rather restricted and that the groups which appear all have periodic cohomology in characteristic $p$.

The kind of finite subgroups of $\GG_n$ that have appeared in the construction of finite resolutions so far are extensions of finite subgroups of $\mathbb{S}_n$ in the following sense.
\begin{definition}\label{defn:extendsubgroups}
For $F_0$ a finite subgroup of $\mathbb{S}_n$, an \emph{extension of $F_0$ to $\mathbb{G}_n$} is a subgroup $F$ of $\mathbb{G}_n$ which contains $F_0$ as a normal subgroup and such that the following diagram commutes:
\[
\xymatrix{0 \ar[r] & F_0 \ar[d] \ar[r] & F \ar[r]  \ar[d] & \Gal \ar[r] \ar[d]^-{\cong} & 0 \\
0 \ar[r]  & \mathbb{S}_n \ar[r] & \GG_n \ar[r] & \Gal \ar[r] & 0}
\]
Here, the rows are exact, the left and middle vertical arrows are the inclusions, and the induced right vertical map is an isomorphism.
\end{definition}
The question of when a finite subgroup $F_0$ of $\mathbb{S}_n$ extends to a finite subgroup of $\GG_n$ is subtle and largely addressed by Bujard in \cite{bujard}. We do not give it much attention here.
\begin{table}
\centering
\caption{The table below lists isomorphism types of maximal finite subgroups of $\mathbb{S}_n$ at various heights and primes. Each isomorphism type listed below belongs to a unique conjugacy class. Here, $C_q$ denotes a cyclic group of order $q$ and $T_{24} \cong Q_8 \rtimes C_3$ is the binary tetrahedral group (the action of $C_3$ on $Q_8$ permutes a choice of generators $i$, $j$ and $k$).
See Hewett \cite{hewett_groups, hewett_normalizers} and Bujard \cite{bujard} for more details. In particular, see \cite{hewett_groups} for the isomorphism type of the semi-direct product on the list below. }
\ \\
\begin{tabular}{|l|l|l|l|}
\hline
Prime $p$ & Height $n$: $k\geq 1$, $p\ndiv m$ &\makecell[l]{Isomorphism Types of Maximal \\
Finite Subgroups in $\mathbb{S}_n$} \\
\hline
\hline
$p\neq 2$ & $n$ not divisible by $p-1$ &$C_{p^n-1}$   \\
\hline
$p\neq 2$ & $n=(p-1)p^{k-1}m$ &\makecell[l]{ $C_{p^n-1}$, and \\ $C_{p^\ell} \rtimes C_{(p^{p^{k-\ell}m} -1 )(p-1)}$, $1\leq \ell \leq k$}  \\
\hline
$p=2$ & $n$ odd & $C_{2(2^n-1)}$ \\
\hline
$p=2$ & $n=2^{k-1}m$ and $k\neq 2$ &$C_{2^\ell(2^{2^{k-\ell}m} -1 )}$, $1\leq \ell \leq k $  \\
\hline
$p=2$ & $n=2 m$ and $m\neq 1$ &\makecell[l]{$C_{2(2^m-1)}$, and \\ $T_{24} \times C_{2^m-1}$} \\
\hline
$p=2$ & $n=2$ &$T_{24}$ \\
\hline
\end{tabular}
\label{tab:finitesub}
\end{table}

\begin{remark}\label{rem:choiceGamma}
For any formal group law $\Gamma$ of height $n$ defined over a perfect field extension $k$ of $\F_p$, one can define the group
\[
\GG(k,\Gamma) = \{(f,i) : \sigma \in \Gal(k/\F_p), \ f \in k[\![x]\!]  \colon \sigma^*\Gamma \xrightarrow{\cong} \Gamma\}.
\]
With this definition, $\GG_n = \GG(\F_{p^n}, \Gamma_n)$. This group was mentioned in Remark~\ref{rem:Ethygeneral}. The group $\mathbb{S}(k, \Gamma) = \End_k(\Gamma)^{\times}$ is the subgroup of $\GG(k,\Gamma)$ consisting of pairs for which $\sigma=\id$. 

In general, both $\mathbb{S}(k, \Gamma)$ and $\GG(k, \Gamma)$ depend on the pair $(k, \Gamma)$. However, since any two formal group laws of height $n$ are isomorphic over $\bF_p$, $\End_{\bF_p}(\Gamma)$ is independent of $\Gamma$, and hence so are $\GG(\bF_p, \Gamma)$ and $\mathbb{S}(\bF_p, \Gamma)$. Since 
\[
\mathbb{S}_n = \mathbb{S}(\F_{p^n}, \Gamma_n) \cong \mathbb{S}(\bF_p, \Gamma_n), 
\]
it follows that for any formal group law $\Gamma$ as above, there is an isomorphism $\mathbb{S}_n \cong \mathbb{S}(\bF_p, \Gamma)$. So, Table~\ref{tab:finitesub} is canonical in the sense that it classifies conjugacy classes of finite subgroups of $ \mathbb{S}(\bF_p, \Gamma)$ for any formal group law $\Gamma$ of height $n$ defined over $\bF_p$. 

However, even if all of the automorphisms of $\Gamma$ are defined over $\F_{p^n}$, so that
\[
\mathbb{S}(\F_{p^n}, \Gamma) \cong \mathbb{S}(\bF_{p}, \Gamma)   \cong \mathbb{S}_n,
\] 
it can still be the case that $\GG(\F_{p^n}, \Gamma)$ and $\GG_n$ are not isomorphic. If this is the case, extensions of a finite subgroup of $\mathbb{S}_n \cong \mathbb{S}(\F_{p^n}, \Gamma)$ to $\GG(\F_{p^n}, \Gamma)$ and $\GG_n$ can have different isomorphism types.
\end{remark}

We now turn to the definition of a few group homomorphisms that play a role in the rest of this paper.
 \begin{definition}\label{defn:determinant}
 The \emph{determinants}\index{determinant} 
 \begin{align*} \det &\colon  \GG_n    \longrightarrow \Z_p^{\times} &  \det &\colon  \mathbb{S}_n   \longrightarrow \Z_p^{\times} 
 \end{align*}
are the homomorphisms defined as follows. The group $ \mathbb{S}_n$ acts on $\mathcal{O}_{\D_n}$ by right multiplication. This action gives a representation $\rho\colon \mathbb{S}_n \to GL_n(\W_n)$. The composite $\det \circ \rho$ has image in the Galois invariants of $\W_n^{\times}$ (see \cite[Section 5.4]{henn_minicourse}), so it induces a homomorphism $\mathbb{S}_n \to \Z_p^{\times}$, which we also denote by $\det$. We extend this homomorphism to $\GG_n$ via the composite
\[ 
\det \colon \GG_n \cong    \mathbb{S}_n\rtimes \Gal \xrightarrow{ \det \times \id }\Z_p^{\times}    \times \Gal \to \Z_p^{\times},    
\]
where the second map is the projection.  \end{definition}

Composing $\det \colon \GG_n \to \Z_p^{\times}$ with the quotient map to $\Z_p^{\times}/\mu \cong \Z_p$ gives a homomorphism 
\begin{equation}\label{eq:zetan}
\zeta_n \colon \GG_n \longrightarrow \Z_p
\end{equation}
where $\mu=C_2$ if $p=2$ and $\mu = C_{p-1}$ if $p$ is odd. \index{$\zeta_n$} 
This corresponds to a class 
\[\zeta_n \in H_c^1(\GG_n, \Z_p) \cong \Hom^c(\GG_n, \Z_p),\] 
where $\Hom^c$ denotes continuous group homomorphisms and $H^1_c$ the continuous cohomology (see Section~\ref{sssec:moravahomalg}). If $p=2$, the determinant also induces a map
\begin{equation}\label{eq:chin}
\chi_n \colon \GG_n \to (\Z_2/4)^{\times}  \cong \Z/2.
\end{equation}
which then represents a class $\chi_n \in H_c^1(\GG_n, \Z/2)$. Let $\widetilde{\chi}_n \in H_c^2(\GG_n, \Z_2)$ be the Bockstein of $\chi_n$, and note that $2\widetilde{\chi}_n =0$.

Denote by $\GG_n^1$ the kernel of $\zeta_n$ and let $\mathbb{S}_n^1 = \mathbb{S}_n \cap \GG_n^1$. The homomorphism $\zeta_n$ is surjective, and necessarily split since $\Z_p$ is topologically free. Therefore,
\begin{align}\label{eq:G1ext}
\GG_n &\cong \GG_n^1 \rtimes \Z_p, &  \mathbb{S}_n &\cong \mathbb{S}_n^1 \rtimes \Z_p.
\end{align}
If $n$ is coprime to $p$, then the splitting is trivial and this is a product.

\begin{remark}\label{rem:3termfiber}
As a consequence of the fact that $\GG_n/\GG_n^1\cong \Z_p$, there is an equivalence $L_{K(n)}S^0 \simeq (E_n^{h\GG_n^1})^{h\Z_p}$. If $\psi \in \GG_n$ is such that $\zeta_n(\psi)$ is a topological generator of $\Z_p$, then we get an exact triangle
\[ 
\xymatrix{ L_{K(n)}S^0 \ar[r] & E_n^{h\GG_n^1} \ar[r]^-{\psi-1} &  E_n^{h\GG_n^1} \ar[r]^-{\delta} & \Sigma L_{K(n)}S^0.}
\]

We also denote by $\zeta_n$ its image $H_c^*(\GG_n, (E_n)_0)$. It is known that $\zeta_n$ is a permanent cycle in the homotopy fixed point spectral sequence, see \cite[Section 8]{DH}.
It detects the composite $S^0 \to E_n^{h\GG_n^1} \xrightarrow{\delta}  \Sigma L_{K(n)}S^0 $ (where the first map is the unit), which is also denoted by $\zeta_n \in \pi_{-1}L_{K(n)}S^0$.
\end{remark}

\subsubsection{The action of the Morava stabilizer group}
\index{Action of the Morava stabilizer group}
We now discuss the action of $\GG_n$ on $(E_n)_*$. Most notably, this problem was first attacked in depth by Devinatz and Hopkins in \cite{DH_action} using the Gross--Hopkins period map (Remark~\ref{rem:grosshop}). A very nice summary of this approach is given by Kohlhaase \cite{Kohlhaase} and we discuss some of the consequences here.

Let $F_n$ be the formal group law over $(E_n)_{0}$ which is a universal deformation of $\Gamma_n$ and was defined in Section~\ref{sec:moravaEthy}. For $\alpha \in \GG_n$ given by a pair $(f,\sigma)$ where $\sigma \in \Gal(\F_{p^n}/\F_p)$ and $f \in \mathbb{S}_n$, the universal property of the deformation $F_n$ implies that there exists a unique pair $(f_{\alpha}, \alpha_*)$ consisting of a continuous ring isomorphism $\alpha_* \colon (E_n)_0 \to (E_n)_0$ and an isomorphism of formal group laws $f_\alpha \colon \alpha_*F_n  \to F_n$
such that 
\begin{align}\label{eq:cong}
({f}_{\alpha},\alpha_*) \equiv (f, \sigma) \mod (p, u_1, \ldots, u_{n-1}).\end{align}
The isomorphism $\alpha_*$ is extended to $(E_n)_*$ by defining $\alpha_*(u) = f'(0)u$. The assignment $\alpha \mapsto \alpha_*$ gives a left action of $\mathbb{G}_n$ on $(E_n)_*$. The action of an element $(\id, \sigma)$ corresponds to the natural action of the Galois group on the coefficients $\WW_n$ in $(E_n)_* \cong \WW_n[\![u_1, \ldots, u_{n-1}]\!][u^{\pm 1}] $, and we denote it by $\sigma_*$. Similarly, if $\alpha = (f, \id)$, we let $f_*$ denote the isomorphism $\alpha_*$.

Computing the action explicitly is difficult and there exists no general formula. However, three cases are fairly simple to deduce from the general description above:
\begin{enumerate}[(a)]
\item If $\alpha$ for $\sigma \in \Gal(\F_{p^n}/\F_p)$, then $\sigma_*$ is the action of the Galois group on the coefficients $\WW_n$. For $x \in \WW_n$, we write $x^{\sigma} = \sigma_*(x)$.
\item If $ \omega \in \mathbb{S}_n$ is a primitive $(p^n-1)$th root of unity, then $\omega_*(u_i)=\omega^{p^i-1}u_i$ and $\omega_*(u) = \omega u$. 
\item If $\psi\in \Z_p^{\times} \subseteq \mathbb{S}_n$ is in the center, then $\psi_*(u_i)=u_i$ and $\psi_*(u)= \psi u$.
\end{enumerate}
Understanding the action more generally is difficult, but we say a few words on this here. 

For $f \in \mathbb{S}_n$, write
$f= \sum\nolimits_{j= 0}^{n-1} f_j \xi^{j}$
for $f_j \in \WW_n$ with $f_0 \in \WW_n^{\times}$ as in \eqref{eq:alphaexpressions}.
The following results due to Devinatz and Hopkins \cite{DH_action} are also given in Theorem 1.3 and Theorem 1.19 of \cite{Kohlhaase}.

\begin{theorem}[Devinatz--Hopkins]\label{thm:DHact}
Let $1\leq i \leq n-1$ and $f_j$ be as above. Then, modulo $(p,u_1, \ldots, u_{n-1})^2$, 
\begin{align*}
f_*(u) &\equiv f_0 u+\sum_{j=1}^{n-1}f_{n-j}^{\sigma^{j}}uu_j , & f_*(uu_i) &\equiv \sum_{j=1}^{i} f_{i-j}^{\sigma^{j}}uu_j + \sum_{j=i+1}^{n} p f_{n+i-j}^{\sigma^{j}}uu_j \ .
\end{align*}
Further, if $f \in \WW_n^{\times} \subseteq \mathbb{S}_n$, so that $f=f_0$ then $f_*(u_i) \equiv f_0^{\sigma^i}f_0^{-1}u_i$ modulo $(u_1, \ldots, u_{n-1})^2$.
\end{theorem}
An example of an immediate consequence of Theorem~\ref{thm:DHact} is the following result. See \cite[Lemma 1.33]{BobkovaGoerss} for a surprisingly simple proof.
\begin{corollary}\label{cor:easyfixpoints}
For all primes $p$ and all heights $n$, the unit $\Z_p \to (E_n)_*$ induces an isomorphism on $\GG_n$ fixed points $\Z_p \cong (E_n)_*^{\GG_n}$. 
\end{corollary}

\begin{remark}[Gross--Hopkins period map]\label{rem:grosshop}
The proof of Theorem~\ref{thm:DHact} relies on one of the deepest results in chromatic homotopy theory, due to Gross and Hopkins \cite{HopkinsGross}\index{Gross--Hopkins period map}, which points towards the mysterious interplay between this subject and arithmetic geometry. Let $K$ be the quotient field of $\WW_n$ and $\Spf((E_n)_0)_K$ be the generic fiber of the formal scheme associated to $(E_n)_0$. Since the division algebra $\mathbb{D}_n$ splits over $K$, i.e., $\mathbb{D}_n \otimes_{\Q_p}K$ is isomorphic to a matrix algebra $M_n(K)$, there is a natural $n$-dimensional $\GG_n$-representation $V_K$. It follows that $\GG_n$ acts on the corresponding projective space $\bP(V_K)$ through projective linear transformations. In \cite{HopkinsGross,GrossHopkins}, Gross and Hopkins construct a \emph{period mapping} that linearizes the action of $\GG_n$ on $\Spf((E_n)_0)_K$: They prove that 
there is an {\'e}tale and $\GG_n$-equivariant map of rigid analytic varieties
\begin{equation}\label{eq:permap}
\xymatrix{\Phi\colon \Spf((E_n)_0)_K \ar[r] & \bP(V_K).}
\end{equation}
Devinatz and Hopkins use this map to prove Theorem~\ref{thm:DHact} and it also features in the computations of Kohlhaase \cite{Kohlhaase}.
\end{remark}

One often needs more precision than that provided by Theorem~\ref{thm:DHact}. Since $f_{\alpha}$ is a morphism of formal group laws, it follows that
\[
f_\alpha( [p]_{\alpha_*F_n}(x)) = [p]_{F_n}(f_\alpha(x)).
\]
This relation contains a lot of information. In practice, it gives a recursive formula to compute the morphism $\alpha_*$ as a function of the $\alpha_j$s. This method is applied explicitly in Section 4 of the paper~\cite{HKM} by Henn--Karamanov--Mahowald. 

However, even with these methods, it is difficult to get good approximations for the action of $\GG_n$. 
If one restricts attention to finite subgroups $F \subseteq \GG_n$, it is sometimes possible to do much better than these kinds of approximations. Recent developments suggest that working with a formal group law other than the Honda formal group law $\Gamma_n$ 
may be better suited to this task. For example, when $n=2$, one can choose to work with the formal group law of a super-singular elliptic curve. The automorphisms of the curve embed in the associated Morava Stabilizer group and one can use geometric information to write explicit formulas for their action on the associated $E$-theory. See Strickland \cite{strickland_level} and \cite[Section 2]{beaudry_towards}. In fact, the spectra $E_2^{hF}$ at height $2$ are the $K(2)$-localizations of various spectra of topological modular forms with level structures. See, for example, \cite{behrens_modular} and Remark~\ref{rem:tmf}. Elliptic curves are not available at higher heights, but there is a hope that the theory of automorphic forms will provide a replacement. This is the subject of \cite{bl_taf}, see also \cite{behrens_chapter}. Finally, equivariant homotopy theory also seems to provide better choices of formal group laws for studying the action of the finite subgroups. See, for example, \cite{HHR, HHRC4} together with \cite{hahnshi}, \cite{c4e4}, and \cite{beaudrybobkovahillstojanoska}.

\subsubsection{Morava stabilizer group: homological algebra}\label{sssec:moravahomalg}

Recall that the $E_2$-term of the descent spectral sequence in Theorem~\ref{thm:descentss} is given by the continuous cohomology of the Morava stabilizer group with coefficients in $(E_n)_*$. The goal of this section is to summarize the homological algebra required to construct these cohomology groups and to then discuss some features specific to $\GG_n$. An important subtlety arising from the homotopical applications we have in mind is that we have to study the continuous cohomology of $\GG_n$ with profinite coefficients, and not merely discrete ones. This theory has been systematically developed by Lazard~\cite{lazard_analyticgroups}; our exposition follows the more modern treatment of Symonds and Weigel~\cite{sw_analyticgroups}. 

Let $G = \lim_i G_i$ be a profinite group, given as an inverse limit of a system of finite groups $(G_i)$ and write $\CC_p(G)$ for the category of profinite modules over  
\[\Z_p\llbracket G \rrbracket = \lim_{i,j}\Z/p^j[G_i]\]
and continuous homomorphisms. The category $\CC_p(G)$ is abelian and has enough projective objects. Moreover, the completed tensor product equips $\CC_p(G)$ with the structure of a symmetric monoidal category with unit $\Z_p$. In order to define a well-behaved notion of continuous cohomology\index{continuous cohomology} for $G$, assume that $G$ is a compact $p$-analytic Lie group in the sense of~\cite{lazard_analyticgroups}. A good reference for properties of $p$-adic analytic groups is \cite{dixon}. Lazard then shows that:
\begin{itemize}
	\item $G$ is of type $p\text{-}\mathrm{FP}_{\infty}$, i.e., $\Z_p$ admits a resolution by finitely generated projective $\Z_p\llbracket G \rrbracket$-modules. It follows that the continuous cohomology of $G$ with coefficients in $M \in \CC_p(G)$, defined as 
	\[
	H_c^*(G,M) = \Ext_{\Z_p\llbracket G \rrbracket}^*(\Z_p,M),
	\]
	is a well-behaved cohomological functor, where the (continuous) $\Ext$-group is computed in $\CC_p(G)$. In particular, there is a Lyndon--Hochschild--Serre spectral sequence and an Eckmann--Shapiro type lemma for open normal subgroups \cite[Theorem 4.2.6 and Lemma 4.2.8]{sw_analyticgroups}. Similarly, continuous homology\index{continuous homology} is defined as 
	\[H^c_*(G,M) = \Tor^{\Z_p\llbracket G \rrbracket}_*(\Z_p,M)\]
	where the (continuous) $\Tor$-group is computed in $\CC_p(G)$.
	\item $G$ is a virtual Poincar\'e duality group\index{virtual Poincar\'e duality group} in dimension $d = \dim(G)$ \cite[Theorem 5.1.9]{sw_analyticgroups}, i.e., there exists an open subgroup $H$ in $G$ such that 
	\[
	H_c^*(H,\Z_p\llbracket H \rrbracket) \cong 
		\begin{cases}
			\Z_p & \text{if } * = d, \\
			0 & \text{otherwise},
		\end{cases}
	\]
	and the length of a projective resolution of $\Z_p \in \CC_p(H)$ can be taken to be $d$. The second property is referred to by saying that the cohomological dimension of $H$ is $d$ and that the virtual cohomological dimension\index{virtual cohomological dimension} of $G$ is $d$; in symbols, $\cdim_p(H) = d$ and $\vcdim_p(G) = d$. The Poincar\'e duality property gives rise to a non-degenerate pairing 
\[
H_c^*(H,\F_p) \otimes H_c^{d-*}(H,\F_p) \longrightarrow H_c^{d}(H,\F_p)\cong \F_p,
\]
thereby justifying the terminology. 
\end{itemize}
The key theorem, proved by Morava~\cite[\S 2.2]{morava} and relying on work by Lazard~\cite{lazard_analyticgroups}, allows us to apply this theory to the Morava stabilizer group:

\begin{theorem}[Lazard, Morava]\label{thm:vcdim}
The Morava stabilizer group $\mathbb{S}_n$ is a compact $p$-analytic virtual Poincar\'e duality group of dimension $n^2$. Further, the group $\mathbb{S}_n$ is $p$-torsion-free if and only if $p-1$ does not divide $n$, and in this case $\vcdim_p(G) = \cdim_p(G) =n^2$. 
\end{theorem}

We note an important immediate consequence of this theorem, which is the underlying reason for the small prime vs.~large prime dichotomy in chromatic homotopy theory. See also Figure~\ref{fig:sse2}:
\begin{corollary}\label{cor:hfpsslarge}
If $p>2$ is such that $2(p-1) > n^2$, then the descent spectral sequence \eqref{eq:hfpss}\index{homotopy fixed point spectral sequence}\index{descent spectral sequence} for $S^0$ collapses at the $E_2$-page with a horizontal vanishing line\index{horizontal vanishing line} of intercept $s=n^2$ (meaning that $E_{2}^{s,t}=0$ for $s>n^2$) and there are no non-trivial extensions.
\end{corollary}
\begin{remark}
The condition $2(p-1) > n^2$ can be improved to $2(p-1) \geq n^2$ using Corollary~\ref{cor:easyfixpoints}. 
\end{remark}

\begin{remark}\label{rem:strongvanishing}
An extremely powerful result of Devinatz--Hopkins is that, for any prime $p$ and any height $n$, there exists an integer $N$ such that, for all spectra $X$, the $K(n)$-local $E_n$-based Adams--Novikov spectral sequence\index{$K(n)$-local $E_n$-based Adams--Novikov spectral sequence} for $X$ (see \eqref{rem:KnlocalANSS})
has a horizontal vanishing line on the $E_{\infty}$-term at $s=N$, although the minimal such $N$ may be greater than $n^2$. For example, when $n=1$ and $p=2$, the homotopy fixed point spectral sequence \eqref{eq:hfpss} has non-trivial elements on the $s=2>1^2$ line at $E_{\infty}$. See \cite[Section 2.3]{BGH} for a proof of the existence of the vanishing line.

Note further that it follows from Corollary~\ref{cor:easyfixpoints} and the existence of the vanishing line that the natural map $\Z_p \to \pi_0L_{K(n)}S^0$ is a nilpotent extension of rings.
\end{remark}

\begin{center}
\begin{figure}[h]
\begin{tikzpicture}[scale=0.35]

\draw[thick] (-8,0) -- (20,0) node[anchor=north west] {$t-s$};

\draw[thick] (0,0) -- (0,6) node[anchor=north west] {$s$};

\draw (0,0) -- (-4,4);

\draw (8,0) -- (4,4);

\draw (16,0) -- (12,4);

\draw[dashed,thick] (-8,4) -- (18,4);

\node[anchor=north] at (0,0) {$0$};

\node[anchor=north] at (8,0) {$2(p-1)$};

\node[anchor=north] at (16,0) {$4(p-1)$};

\node[anchor=north] at (-8,0) {$-2(p-1)$};

\node[anchor=north west] at (0,4) {$n^2$};
\node[anchor=north] at (-4,0) {$-n^2$};
\end{tikzpicture}
\caption{{The $E_2$-term of $E_2^{s,t} \cong H^s_c(\GG_n, (E_n)_t)  \Rightarrow \pi_{t-s} L_{K(n)}S^0$ for $p>2$ and $ 2(p-1) > n^2$. The dashed line indicates the horizontal vanishing line at $E_2$, that is, $E_2^{s,t}=0$ when $s>n^2$. The non-zero contributions are concentrated on the lines of slope $-1$ that intercept the $(t-s)$-axis at multiples of $2(p-1)$.
}}
\label{fig:sse2}
\end{figure}
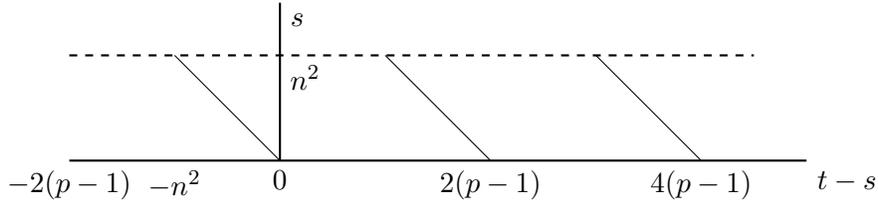
\end{center}

In order to run the descent spectral sequence computing $\pi_*L_{K(n)}S^0$, we have to come to grips with $H_c^*(\GG_n, (E_n)_*)$, an extremely difficult problem. However, if one restricts attention to $H_c^*(\GG_n, (E_n)_0)$, the computation appears to radically simplify in a completely unexpected way. Let $\iota \colon \WW_n \to (E_n)_0$ be the natural inclusion. The following has been shown to be true at all primes when $n\leq 2$, see \cite{shimyabe, behse2, Kohlhaase, GoerssSplit, BGH, Beaudry2018}:
\begin{conjecture}[Vanishing conjecture]\label{conj:vanishing}\index{Vanishing conjecture}
Let $p$ be any prime and $n$ be any height. The map $\iota$ induces an isomorphism
\[ 
\xymatrix{\iota_* \colon H_c^*(\GG_n, \WW_n) \ar[r]^-{\cong} & H_c^*(\GG_n, (E_n)_0).}
\]
\end{conjecture}

\begin{remark}
The conjecture is so named because it implies that the cohomology of the $\GG_n$-module $(E_n)_0/\WW_n$ vanishes in all degrees. Note further that if one proves that $\WW_n/p \to (E_n)_0/p$ induces an isomorphism on cohomology, then Conjecture~\ref{conj:vanishing} follows formally. 
\end{remark}

As we will see in Section~\ref{sec:cscmore}, this conjecture and the accompanying computations informs our understanding of $L_{K(n)}S^0$, the essence of which is distilled in the formulation of the chromatic splitting conjecture. In fact, what makes Conjecture~\ref{conj:vanishing} particularly appealing is the fact that $H_c^*(\GG_n,\WW_n)$ appears to be rather simple when $p$ is large with respect to $n$. 

Rationally, we have some partial understanding due to work of Lazard \cite[Remark 2.2.5]{morava} and \cite[Rem.~2.2.5]{morava}, who established an isomorphism for all heights and primes
\begin{equation}\label{eq:rationaliso}
H_c^*(\GG_n,\WW_n) \otimes \Q \cong \Lambda_{\Z_p}(x_1,\ldots,x_n)\otimes \Q,
\end{equation}
where $\Lambda_{\Z_p}(x_1,\ldots,x_n)$ is the exterior algebra over $\Z_p$ on $n$ generators in degrees $\deg(x_i) = 2i-1$. Here, the class $x_1$ is $\zeta_n$ as defined in \eqref{eq:zetan}. Furthermore, when $p$ is large with respect to $n$, it is believed that there is an isomorphism \eqref{eq:rationaliso} before rationalization. 
\begin{conjecture}
If $p\gg 0$, then $H_c^*(\GG_n,\WW_n) \cong \Lambda_{\Z_p}(x_1,\ldots,x_n)$.
\end{conjecture}

\begin{remark}\label{rem:WGmodules}
For our chromatic applications, we need a mild extension of the setup presented above. Here and below, $\WW_n = W(\F_{p^n})$ denotes the $\GG_n$-module whose action is the restriction along
$\GG_n \to  \Gal$ of the natural action of $\Gal$ on $W(\F_{p^n})$. We write $w^g=g(w)$ for this action. For $ G\subseteq \GG_n$, define the \emph{twisted group ring}\index{twisted group ring} to be
\[
\twistgr{\WW_n}{G} := \lim_{i,j} \WW_n/\mathfrak{m}^j\langle G_i\rangle
\]
with $G$-twisted multiplication determined by the relations $g\cdot r = g(r)g$ for $r\in \WW_n$ and $g\in G$. We let $\GMod{\WW_n}{G}$ be the category of profinite left $\twistgr{\WW_n}{G}$-modules. These are profinite  abelian groups $M = \lim_kM_k$ with a continuous action $\twistgr{\WW_n}{G} \times M \to M$.  
If $H \subseteq G$ is a closed subgroup and $M$ is a left $\twistgr{\WW_n}{H}$-module, then
\begin{align}\label{eq:uparrow}
M\upar{H}{G} := \twistgr{\WW_n}{G} \otimes_{\twistgr{\WW_n}{H}}M =\lim_{i,j,k} \left(\WW_n/\mathfrak{m}^j[G_i] \otimes_{\twistgr{\WW_n}{H}} M_k\right) 
\end{align}
is a left $\twistgr{\WW_n}{G}$-module.

One can show that the homological algebra summarized above also works in the context of profinite modules over twisted group rings. Note that there is a functor from
$\Mod_{\Z_p\llbracket G \rrbracket}$ to $\GMod{\WW_n}{G}$ which sends a $\Z_p\llbracket G \rrbracket$-module $M$ to the $\WG{G}$-module $\WW_n\otimes_{\Z_p}M$ with action given by $g( w \otimes m) = w^{g}\otimes g(m)$. This allows us to transport constructions in $\GMod{\Z_p}{G}$ to constructions in $\GMod{\WW_n}{G}$. 
\end{remark}

We now come to another important construction in chromatic homotopy theory, namely the $(E_n)_*$-module 
\[
(E_n)_*^{\vee}X= \pi_*L_{K(n)}(E_n\smsh X)
\]
associated to a spectrum $X$. 
The action of $\GG_n$ on $E_n$ induces an action on $(E_n)_*^{\vee}X$ compatible with the $(E_n)_*$-action. Moreover, let $\fm = (p,u_1,\ldots,u_{n-1})$ be the maximal ideal of $(E_n)_0$ and, for $s\ge 0$, let $\mathbb{L}_s$ be the $s$th left derived functor of $\fm$-adic completion on $\Mod_{(E_n)_*}$. There is a strongly convergent spectral sequence
\[
\mathbb{L}_s(\pi_*(E_n \smsh X))_t \Longrightarrow (E_n)_{s+t}^{\vee}X
\]
which in particular implies that the canonical map $(E_n)_*^{\vee}X \to \mathbb{L}_0((E_n)_*^{\vee}X)$ is an isomorphism. Such $(E_n)_*$-modules are called $\mathbb{L}$-complete and we refer the interested reader to \cite[Appendix A]{HovStrick} for a more thorough treatment. Taken together, this structure is called the \emph{Morava module of $X$}: 

\begin{definition}[Morava modules]\label{defn:moravamodules}\index{Morava module}
A \emph{Morava module} $M$ is an $\mathbb{L}$-complete $(E_n)_*$-module equipped with an action by $\GG_n$ in $\mathbb{L}$-complete modules that is compatible with the action on $(E_n)_*$. That is, for every $g\in \GG_n$, $e\in (E_n)_*$ and $m \in M$, $g(em) = g(e)g(m)$. A morphism of Morava modules is a continuous map of $(E_n)_*$-modules that preserves the action. We denote the category of Morava modules by $\moravamod{n}$.
\end{definition}
By the discussion above, $(E_n)_*^{\vee}X$ is a Morava module for any spectrum $X$ and we obtain a functor
\begin{align}\label{eq:moravamodulefunctor}
(E_n)_*^{\vee} (-) := \pi_*L_{K(n)}(E_n\smsh -) \colon \Sp \longrightarrow \moravamod{n}. \end{align}
This functor detects and reflects $K(n)$-local equivalences, but has the advantage that $(E_n)_*^{\vee} (-) $ comes equipped with an action of $\GG_n$. This extra structure proves to be extremely powerful for computations, and is one of the reasons why Morava modules play a central role in the field.

For more information on Morava modules, we refer the reader to \cite[Section 1.3]{BobkovaGoerss} and \cite[Section 2]{ghmr}, noting that authors often simply write $(E_n)_*X= \pi_*L_{K(n)}(E_n\smsh X)$ as opposed to the non-completed homology $(E_n)_*X= \pi_*(E_n \smsh X)$, but we will not do so here. Note also that, if $X$ is finite, then $(E_n)_*X \cong (E_n)_*^{\vee}X$.

\begin{remark}\label{rem:periodicity}\index{periodicity}
For $F$ a finite subgroup of $\GG_n$, the action of $\GG_n$ on $(E_n)_*$ restricts to an action of $F$. 
We can also consider the category $\Mod_{(E_n)_*}^F$ of $\mathbb{L}$-complete $(E_n)_*$-modules equipped with an action of $F$. Then $(E_n)_*$ is periodic as an object in $\Mod_{(E_n)_*}^F$ since the element 
$N=\prod_{g\in F} g(u)$
for $u \in (E_{n})_{-2}$ as in \eqref{eq:Estariso} is an invariant unit. Let $d^{\alg}_F$ be the smallest integer such that $(E_n)_* \cong (E_n)_{*+d^{\alg}_F}$ in $\Mod_{(E_n)_*}^F$. 
This leads to an isomorphism of Morava modules
\begin{align*}
 (E_n)_*^{\vee}E_n^{hF} & \cong \Map^c(\GG_n/F, (E_n)_*)   \cong \Hom^c_{\WW_n}(\WW_n\upar{F}{\GG_n}, (E_n)_*)  
\end{align*}
closely related to \eqref{eq:EnEn} and it implies that
\[ (E_n)_*^{\vee}E_n^{hF}  \cong  (E_n)_*^{\vee}\Sigma^{d^{\alg}_F}E_n^{hF}. \]
However, $E_n^{hF} $ need not be equivalent to $\Sigma^{d^{\alg}_F}E_n^{hF}$. Nonetheless, because of the strong vanishing line discussed in Remark~\ref{rem:strongvanishing}, some power of $N$ is a permanent cycle and gives rise to a periodicity generator for $E_n^{hF}$, so for some multiple $d_{F}^{\mathrm{top}}$ of $d^{\alg}_F$, there is an equivalence
$E_n^{hF}\simeq \Sigma^{d^{\mathrm{top}}_F}E_n^{hF}$.

For example, at $p=2$, $E_1$ is $2$-complete complex $K$-theory\index{complex $K$-theory} and $E_1^{hC_2}$ is the $2$-complete real $K$-theory spectrum $KO$.\index{real $K$-theory} We have:
\begin{align*}
K_*^{\vee}KO &\cong  K_*^{\vee}\Sigma^4KO, &  KO &\not\simeq \Sigma^4KO,  & KO &\simeq \Sigma^8KO.
 \end{align*}
 \end{remark}

\section{$K(1)$-local homotopy theory}\label{sec:K1}\label{sec:cscheight1}

In this section, we tell a part of the chromatic story at height $n=1$ as a warm up for the more complicated ideas needed to study higher heights. The contents of this section are classical and can be found in various forms throughout the literature, for example, Adams and Baird \cite{adamsbaird}, Bousfield \cite{bousfield_locspectra, bousfield_odd}, Ravenel \cite[Theorem 8.10, 8.15]{Rav84}. See \cite[Section 6]{henn_minicourse} for a more recent treatment, and \cite[Section 4]{BGH} for more details on the case $p=2$.

\subsection{Morava $E$-theory and the stabilizer group at $n=1$}
At height $n=1$, Morava $E$-theory is the $p$-completed complex $K$-theory\index{complex $K$-theory} spectrum, which we simply denote by $K$. There is an isomorphism $K_* \cong \Z_p[u^{\pm 1}]$ for a unit $u \in K_{-2}$ which can be chosen so that $u^{-1}\in K_2$ is the Bott element.\index{Bott element}

In this case, $\GG_1  = \mathbb{S}_1 \cong \Z_p^{\times}$  corresponds to the $p$-completed Adams operations.\index{Adams operations} The action of $\mathbb{S}_1$ on $K_*$ is the $\Z_p$-algebra isomorphism determined by 
\begin{align}\label{eq:actadams}
\alpha_*(u) = \alpha u 
\end{align}
for $\alpha \in \Z_p^{\times}$. 
The keen reader will notice that this is the inverse of the action of the Adams operations, which is given by $\alpha_*(u) = \alpha^{-1} u$. This comes from switching a right action to a left action.

The map $L_{K(1)}S^0 \to K^{h\Z_p^{\times}}$ of \eqref{eq:galoisspectra}
is a $\Z_p^{\times}$ pro-Galois extension. We use this extension to compute the homotopy groups of $\pi_*L_{K(1)}S^0$.\index{$K(1)$-local sphere} One can take the direct approach of computing the spectral sequence of \eqref{eq:hfpss}
\begin{equation}\label{eq:hfpssk1}
E_2^{s,t} = H_c^s(\Z_p^{\times}, K_t ) \Longrightarrow \pi_{t-s}L_{K(1)}S^0.
\end{equation}
In fact, this spectral sequence collapses at the $E_2$-page at odd primes and at the $E_4$-page at the prime $2$. This is not a hard computation, but we take a different path in order to illustrate the finite resolution philosophy.

\subsection{Finite resolution at height $n=1$}\label{sssec:resn1}

Here, we describe our first example of a finite resolution.\index{finite resolution} Let $C_m$ denotes a cyclic group of order $m$, $\mu =C_2$ if $p=2$, and $\mu =C_{p-1}$ if $p$ is odd. Then, $\Z_p^{\times} \cong \mu \times \Z_p$, where the $\Z_p$ corresponds to the subgroup of units congruent to $1$ modulo $p$ if $p$ is odd, and to those congruent to $1$ modulo $4$ is $p=2$. We let $\psi$ be a topological generator for this factor of $\Z_p$. The notation is meant to be reminiscent of the Adams operations. We will make a choice for $\psi$ below in \eqref{eq:psichoice}. 
\begin{remark} 
The spectrum $K^{hC_{p-1}} $ is the unit component in the splitting of the $p$-completed complex $K$-theory spectrum $K$ into Adams summands if $p$ is odd, and $K^{hC_2}$ is the $2$-completed real $K$-theory spectrum if $p=2$.
\end{remark}
The $K(1)$-local sphere can be obtained by an iterated fixed points construction:
\[
L_{K(1)} S^0 \simeq K^{h\Z_p^{\times}} \simeq  (K^{h\mu})^{h\Z_p}.
\] 
Since $\psi \in \Z_p$ is a topological generator, taking homotopy fixed point with respect to $\Z_p$ is equivalent to taking the homotopy fiber of the map $\psi-1$. Therefore, there is a fiber sequence
\begin{align}\label{eq:K1localfiber}
\xymatrix{L_{K(1)}S^0 \ar[r] & K^{h\mu} \ar[rr]^-{ \psi-1} & & K^{h\mu} \ar[r] & \Sigma L_{K(1)}S^0 }. \end{align}
This is a \emph{finite resolution} of $L_{K(1)}S^0$ as will be defined in Definition~\ref{defn:finresolution} below. 

To construct finite resolutions at higher heights where the structure of the Morava stabilizer group is more intricate, we start by attacking the problem in algebra and then we transfer algebraic constructions to topology. We give a quick overview of how this would happen at height $1$ to give the reader something to think of while reading Section~\ref{sec:fin}.

\ \\
\noindent
\emph{Step 1: Algebraic resolution.}\label{algebraic resolution}
The group $\Z_p$ is topologically free of rank one and there is an exact sequence of left $\Z_p^{\times}$-modules
\begin{align}\label{eq:algres1}
\xymatrix{ 0\ar[r] & \Z_p\upar{\mu}{\Z_p^{\times}} \ar[r]^-{\psi-1} &\Z_p\upar{\mu}{\Z_p^{\times}} \ar[r] & \Z_p  \cong \Z_p\upar{\Z_p^{\times}}{\Z_p^{\times}} \ar[r] & 0.}
\end{align}
Here, $\Z_p[\![\Z_p^{\times}]\!] = \lim_{i,j} \Z/p^i[(\Z/p^j)^{\times}]$ is the completed group ring, which was discussed in Section~\ref{sssec:moravahomalg}, and $\Z_p\upar{\mu}{\Z_p^{\times}} \cong \Z_p[\![\Z_p^{\times}]\!] \otimes_{\Z_p[\mu]} \Z_p$.
This is a projective resolution of $\Z_p$ as a $\Z_p^{\times}$-module if and only if $p>2$. See Remark~\ref{rem:notproj} below on this point.
Applying $\Hom^c_{\Z_p}( -,  K_* )$ to \eqref{eq:algres1} 
gives a short exact sequence of Morava modules
\begin{align}\label{eq:morava}
\xymatrix@C=1.5pc{ K_* \ar[r] &\Hom^c_{\Z_p} (\Z_p\upar{\mu}{\Z_p^{\times}}  ,K_* )  \ar[rrr]^-{\Hom_{\Z_p}^c( \psi-1,  K_* ) }  & &  & \Hom^c_{\Z_p} (\Z_p\upar{\mu}{\Z_p^{\times}} , K_* ) }
\end{align}

\ \\
\noindent
\emph{Step 2: Topological Resolution.}
The second step is to prove that the algebraic resolution has a topological realization. More precisely, \eqref{eq:morava} is an exact sequence in the category of Morava modules $\moravamod{1}$. 
As was described in \eqref{eq:moravamodulefunctor}, there is a functor 
\[
K_*^{\vee}(-)=(E_1)_*^{\vee}(-) \colon \Sp \longrightarrow  \moravamod{1}.
\]
When we have an algebraic resolution of length $1$, a \emph{topological realization of \eqref{eq:algres1}} is a choice of fiber sequence in the category of $K(1)$-local spectra
\begin{align}\label{eq:diag}
\xymatrix{ \EE_{-1} \ar[r]^{\delta_{-1}} & \EE_{0} \ar[r]^{\delta_0} &  \EE_{1}& }  \end{align}
where $\EE_{0}$ and $\EE_{1}$ are finite wedges of suspensions of spectra of the form $K^{hF}$ for $F\subseteq \GG_1$ a finite subgroup, such that, up to isomorphism of complexes, \eqref{eq:morava} is the complex of Morava modules obtained from \eqref{eq:diag} by applying $K_*^{\vee}(-)$. If $\EE_{-1} = L_{K(1)}S^0$, then this is an example of a finite resolution of the $K(1)$-local sphere. 

\begin{remark}The case when the algebraic resolution has length $d\geq 1$ is discussed in the next section. We will see in Definition~\ref{defn:finresolution} that the definition of a finite topological resolution of length greater than $1$ is more subtle. See also Example~\ref{ex:finres1}.
\end{remark}

There is no algorithm for finding a topological realization. A priori, one may not exist, and if it does, it may not be unique.  
Without a priori knowledge of the existence of \eqref{eq:K1localfiber}, the key observations for finding a topological realization of \eqref{eq:algres1} are  
\begin{itemize}
\item the isomorphism of Morava modules
\[  K_*^{\vee}K^{h\mu}  \cong \Map^c (\Z_{p}^{\times}/\mu,K_* )
\cong \Hom^c_{\Z_p} (\Z_p\upar{\mu}{\Z_p^{\times}} ,K_* ) ,\] 
and 
\item the fact that $\Hom^c_{\Z_p}( \psi-1,  K_* )  = K_*^{\vee}(\psi-1)$.
\end{itemize}
Knowing these facts, \eqref{eq:morava} can be identified with the short exact sequence of Morava modules
\begin{align}\label{eq:resheight1algK}
\xymatrix{ K_* \ar[r] & K_*^{\vee}K^{h\mu} \ar[rrr]^-{K_*^{\vee}(\psi-1)}  & &  & K_*^{\vee}K^{h\mu}}.
\end{align}
Given this, we let $\EE_{-1} = L_{K(1)}S^0 $, $\EE_0 = \EE_1 = K^{h\mu}$. We let $\delta_{-1}$ be the unit and $\delta_{0}$ be $\psi-1$.
It follows that the fiber sequence 
\begin{align*}
\xymatrix{L_{K(1)}S^0 \ar[r] & K^{h\mu} \ar[rrr]^-{\psi-1}  & &  &  K^{h\mu}}.
\end{align*}
is a topological realization as it gives rise to \eqref{eq:morava} under the functor $K_*^{\vee}(-)$. This is our first example of a finite resolution of $L_{K(1)}S^0$.

\begin{remark} We \emph{did} make choices here and different choices could have given a different topological realization. For example, for $p=2$, $K^{hC_2}  \simeq KO$ and $K_*^{\vee} KO \cong  K_*^{\vee}\Sigma^4 KO$, yet $KO \not\simeq \Sigma^4 KO $. In fact, we could have constructed a topological realization using $ \Sigma^4 KO$ instead of $KO$. Such a resolution is described below in \eqref{eq:realizationP1}. The resolution described there is a topological realization of the algebraic resolution \eqref{eq:algres1}, but it is not a finite resolution of the sphere as $\EE_{-1} = P_1 \not\simeq L_{K(1)}S^0$.
\end{remark}

\subsection{Homotopy groups and chromatic reassembly}\label{sec:reassemlyn1}

The long exact sequence on homotopy groups associated to \eqref{eq:K1localfiber} allows one to compute $\pi_*L_{K(1)}S^0$ from $\pi_*K^{h\mu}$ and knowledge of the action of $\psi$. The homotopy groups of $K^{h{\mu}}$ are computed using the homotopy fixed point spectral sequence
\[
E_2^{s,t} \cong H^s({\mu}, \pi_tK) \Longrightarrow \pi_{t-s}K^{h{\mu}}.  
\]
Recall that $\mu = C_{p-1}$ if $p$ is odd and $C_2$ if $p=2$. So computing group cohomology with coefficients in $K_* = \Z_p[u^{\pm1}]$ is not so bad given the explicit formula \eqref{eq:actadams}. We get
\[
H^*({\mu}, \pi_*K) \cong \begin{cases} \Z_p[u^{\pm (p-1)}]  & p \neq 2 \\
 \Z_2[\eta, u^{\pm 2}]/(2\eta)  & p=2, 
\end{cases}  
\]
where the $(s,t)$ bidegree of $\eta$ is $(1,2)$. The element $\eta$ detects the Hopf map in $\pi_1S^0$.
For $p$ odd, the spectral sequence collapses for degree reasons. For $p=2$, the fact that $\eta^4=0$ in $\pi_*S^0$ implies a differential $d_3( u^{-2}) =\eta^3$, and the spectral sequence collapses at $E_4$ for degree reasons. So, we have
\[
\pi_{*}K^{h{\mu}} \cong \begin{cases} \Z_p[\beta^{\pm 1}]  & p \neq 2, \ |\beta|=2(p-1)  \\
 \Z_p[\eta, \alpha, \beta^{\pm 1}]/(2\eta, \eta^3, \alpha^2-4\beta)  & p=2, \ |\eta|=1, \ |\alpha|=4, \  |\beta| = 8.
\end{cases}  
\]
If $p$ is odd, $\beta \in \pi_{2(p-1)}$ is detected by $u^{1-p}$. If $p=2$, $\eta \in \pi_1$ is detected by the same-named class on the $E_2$-page, $\alpha\in \pi_4$ is detected by $2u^{-2}$ and $\beta \in \pi_8$ is detected by $u^{-4}$.

\begin{remark}
The differential $d_3( u^{-2}) =\eta^3$ can be obtained as a consequence of the \emph{slice differentials theorem} \cite[Theorem 9.9]{HHR}. This is an overkill for this particular example which follows from classical considerations. However, we mention this here since the slice differentials theorem also implies differentials at higher heights in spectral sequences computing $\pi_*E_n^{hF}$ for finite subgroups $F \subseteq \GG_n$.
\end{remark}

Now, we turn to computing the long exact sequence on homotopy groups associated to \eqref{eq:K1localfiber}. Choose an element $\psi$ of $\Z_p^{\times}$ which satisfies
\begin{align}\label{eq:psichoice}
\psi^{-1}= \begin{cases}(1+p)& p\neq 2 \\
5  & p= 2.
\end{cases}
\end{align}
There are other possible choices: One could choose any element in $\Z_p^{\times}$ such that the image of $\psi^{-1}$ in $\Z_p^{\times}/\mu$ is a topological generator. The outcome of these calculations are independent of the choice.

\begin{figure} 
\center
\includegraphics[width=\textwidth]{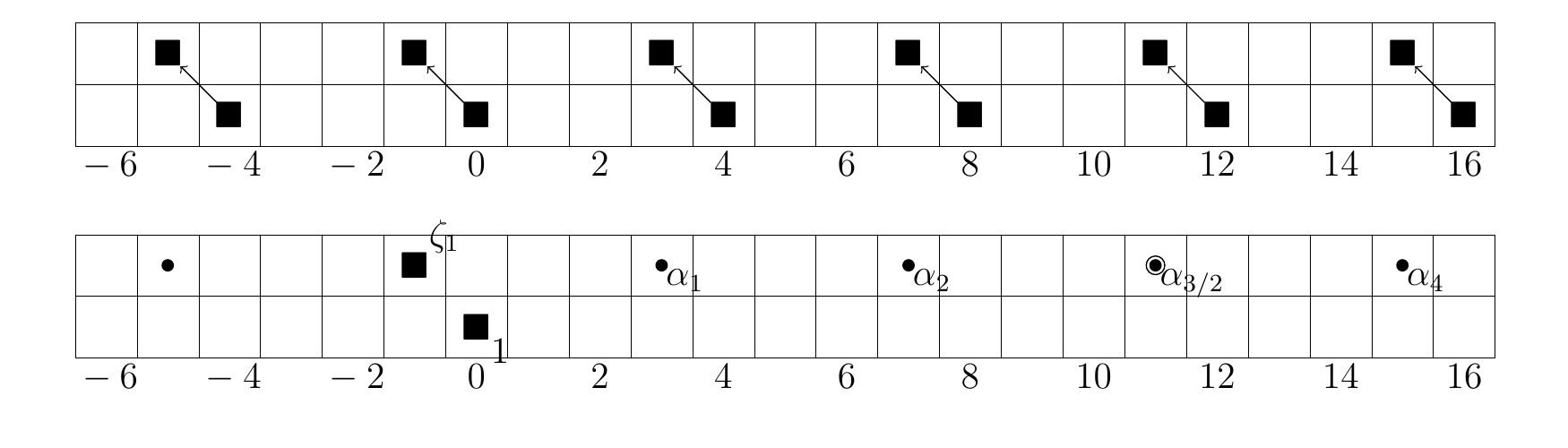}
\captionsetup{width=\textwidth}
\caption{The long exact sequence on homotopy groups associated to \eqref{eq:K1localfiber} and computing $\pi_{*}L_{K(1)}S^0$ at $p=3$. It is drawn as a spectral sequence in Adams grading $(t-s,s)$ with with $E_1^{s, t} \cong \pi_tK^{h\mu}$  for $s=0, 1$. The arrows denote the $d_1$-differentials, $d_1 \colon E_1^{0, t} \to E_1^{1, t}$, which is just the connecting homomorphism.
The top chart is the $E_1$-term and the bottom chart the $E_{\infty}$-term. A $\blacksquare$ is a $\Z_3$, a $\bullet$ a $\Z/3$ and a \circled{$\bullet$} a $\Z/9$.}
\label{fig:K1p3}
\end{figure}

\begin{figure} 
\center
\includegraphics[width=\textwidth]{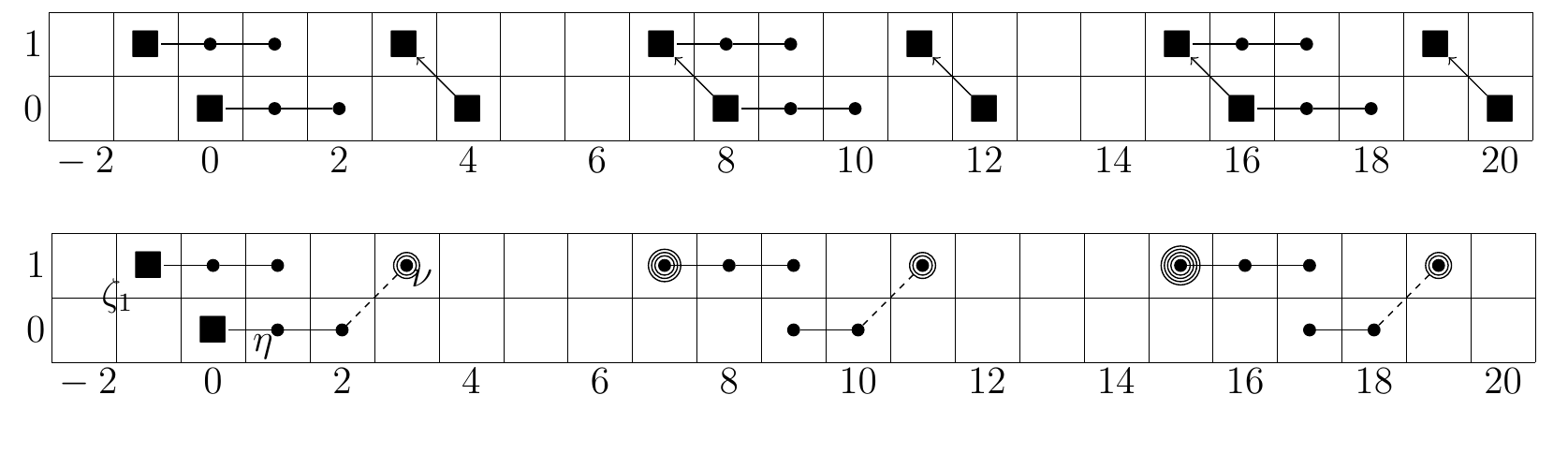}
\captionsetup{width=\textwidth}
\caption{The long exact sequence on homotopy groups associated to \eqref{eq:K1localfiber} and computing $\pi_{*}L_{K(1)}S^0$ at $p=2$. It is drawn using the same convention as in Figure~\ref{fig:K1p3}, except that a $\blacksquare$ is a $\Z_2$, a $\bullet$ a $\Z/2$, a \circled{$\bullet$} a $\Z/4$, etc. Dashed arrows denote exotic multiplications by $\eta$.
}
\label{fig:K1p2}
\end{figure}

From \eqref{eq:actadams}, we deduce that the action of $\psi$ is then given by
\begin{align*}\psi_*(\beta) &= \begin{cases} (1+p)^{p-1}\beta & p\neq  2\\
5^{4} \beta & p=2
\end{cases}, & \psi_*(\alpha) &= 5^{2} \alpha, & \psi_*(\eta ) &=\eta.
\end{align*}
Let $v_p(k)$ denote the $p$-adic valuation of $k \in \Z$. For $p$ odd, the long exact sequence on homotopy groups gives
\[ 
\pi_*L_{K(1)}S^0 = \begin{cases} \Z_p  & *=0,-1 \\
\Z/p^{v_p(k)+1} & *=2k(p-1)-1 \\
0 & \text{otherwise.}
\end{cases}
\]
This is depicted in Figure~\ref{fig:K1p3} for $p=3$.
 For $p=2$, we have
\[ 
\pi_*L_{K(1)}S^0 = \begin{cases} \Z_2  & *=-1 \\
\Z_2 \oplus \Z/2 & *=0 \\
\Z/2 & *=0, 2\mod 8, *\neq 0 \\
\Z/2 \oplus \Z/2 & *=1\mod 8 \\
  \Z/8 & *=3\mod 8 \\
\Z/2^{v_2(k)+4} & *=-1+8k, k\neq 0 \\
0 & *=4,5,6 \mod 8 .
\end{cases}
\]
This is depicted in Figure~\ref{fig:K1p2}. One has to argue that there is no additive extension in degrees $1 \mod 8$ but we do not do this here.

\begin{remark}\label{rem:notproj}
The dichotomy between $p=2$ and odd primes in the computations is an instance of the general phenomena which was discussed in Section~\ref{sssec:moravahomalg} and is revisited in Section~\ref{sec:assymptotics} below.
That is, when $p$ is large with respect to $n$, chromatic homotopy theory becomes algebraic (see for example Corollary~\ref{cor:hfpsslarge}). On the other hand, when $p$ is small the stabilizer group $\mathbb{S}_n$ might contain $p$-torsion and this appears to reflect interesting topological phenomena. Here, $\mathbb{S}_1 \cong\Z_2^{\times}$ contains $2$-torsion at $p=2$ and there are differentials in the spectral sequence computing the homotopy groups of $KO \simeq K^{hC_2}$, a much more intricate spectrum than the Adams summand $K^{hC_{p-1}}$ at odd primes whose homotopy fixed point spectral sequence collapses at the $E_2$-page.
\end{remark}

\begin{remark}
The $\Z_p$ summand in $\pi_{-1}L_{K(1)}S^0$ is generated by the image of the composite $S^0 \to K^{h\mu} \to \Sigma^{-1}L_{K(1)}S^0$ where the first map is the unit and the second is the connecting homomorphism of \eqref{eq:K1localfiber}. We call this map and the homotopy class it represents $\zeta_1 \in \pi_{-1}L_{K(1)}S^0$. It is detected in \eqref{eq:hfpssk1} by the same-named class 
\[\zeta_1 \in H_c^1(\Z_p^{\times}, K_0) \cong \Hom^c(\Z_p^{\times}, \Z_p) \]
corresponding to the projection $ \Z_p^{\times} \to \Z_p^{\times}/\mu \cong \Z_p$. See \eqref{eq:zetan} and Remark~\ref{rem:3termfiber} for analogues at higher heights.
\end{remark}
\begin{remark}\label{rem:L1V0}
An easy computation that will be relevant later is that of $\pi_*(L_{K(1)}S^0/p)$ for $p$ odd. The descent spectral sequence
\[
H_c^s(\Z_p^{\times}, K_t/p) \Longrightarrow \pi_{t-s}(L_{K(1)}S^0/p).
\]
collapses with no extensions and
\[ 
\pi_*(L_{K(1)}S^0/p) \cong  \F_p[v_1^{\pm 1}]\otimes \Lambda_{\F_p}(\zeta_1),
\] 
where $v_1 = u^{1-p}$. We abuse notation by denoting the composite $\Z_p^{\times} \xrightarrow{\zeta_1} \Z_p^{\times} /\mu \cong \Z_p \to \Z/p$ also by $\zeta_1 \in H^1_c(\Z_p^{\times}, K_0/p)$.
\end{remark}

Finally, we turn to the problem of chromatic reassembly at height $n=1$. The chromatic fracture square\index{chromatic fracture square} \eqref{eq:chromaticfracture} in this case gives
\[\xymatrix@=1.5pc{F_1 \ar[r] \ar[d]_-{\simeq} & L_1S^0 \ar[r] \ar[d]  & L_{K(1)}S^0 \ar[d] \\
F_1 \ar[r] & L_0S^0_p \ar[r]  & L_0L_{K(1)}S^0
} \]
where $F_1$ is the fiber of the horizontal maps. In particular, it is the fiber of the map $L_0S^0_p \to L_0L_{K(1)}S^0$ induced by the unit. Since $L_0$ is rationalization, there is an isomorphism $\pi_*L_0L_{K(1)}S^0  \cong p^{-1}\pi_*L_{K(1)}S^0$. From the above calculations, we see that the map $1 \vee \zeta_1 \colon S^0 \vee S^{-1} \to L_{K(1)}S^0$ induces an equivalence
\begin{equation}\label{eq:cscn1first}
L_0L_{K(1)}S^0  \simeq L_0S^0_p \vee L_0S^{-1}_p.
\end{equation}
In particular, $L_0S^0_p \to L_0L_{K(1)}S^0$ is split and $\Sigma F_1 \simeq  L_0 S_p^{-1}$.
This proves the \emph{strong chromatic splitting conjecture}\index{chromatic splitting conjecture} for $n=1$, which will be stated in general in Section~\ref{sec:cscmore}.

We get 
the following diagram from the long exact sequence on homotopy groups associated to the fiber sequence $L_1S^0_p \to L_{K(1)}S^0 \to \Sigma  F_1\simeq  L_0 S_p^{-1}$:
\begin{align*}
\xymatrix{\pi_{-1}L_1S^0_p \ar[r]  \ar@{=}[d] & \pi_{-1} L_{K(1)}S^0 \ar@{=}[d] \ar[r] & \pi_{-2} F_1 \ar@{=}[d] \ar[r] &  \pi_{-2} L_1S^0_p \ar[r] \ar@{==}[d] & \pi_{-2} L_{K(1)}S^0  \ar@{=}[d] \\
0\ar[r] & \Z_p \ar[r] &\Q_p \ar[r] & \Q_p/\Z_p \ar[r] & 0.  }
\end{align*}
Piecing the rest of the long exact sequence on homotopy groups together gives
\[ 
\pi_*L_1S^0_p \cong \Z_p  \oplus \Sigma^{-2}\Q_p/\Z_p \oplus \mathrm{Tor}(\pi_*L_{K(1)}S^0 ), 
\]
where the $\Z_p$ is in degree $0$ and comes from the summand $\Z_p \subseteq \pi_0L_{K(1)}S^0$, this inclusion being an isomorphism when $p$ is odd. The summand $\Q_p/\Z_p$ is in degree $-2$ and $ \mathrm{Tor}$ denotes the torsion subgroup.

\bigskip

In the next sections, we review these topics at higher heights. While we are not able to do such an explicit analysis for $n\geq 2$, the tools and ideas described above do generalize and we give an overview of some of the techniques available to study the $K(n)$-local category and the $K(n)$-local sphere.

\section{Finite resolutions and their spectral sequences}\label{sec:fin}
\index{finite resolution}
We now describe a \emph{recipe} for the construction of finite resolutions of the $K(n)$-local sphere. We note that every step of this procedure requires hard work specific to the height and the prime. We then illustrate the general formalism with many examples at height $n=2$ in Section~\ref{sec:resheight2}. 
Some applications of these finite resolutions will then be discussed in the next section on the chromatic splitting conjecture and local Picard groups. References for this material are \cite{ghmr, henn_res, henn_minicourse}.

\subsection{What is a finite resolution}

\begin{definition}\label{defn:finresolution}
A \emph{finite resolution of $L_{K(n)}S^0$} of length $d$ is a diagram 
\begin{align}\label{eq:finres}
\xymatrix{ 
L_{K(n)}S^0 =X_d \ar[r]^-{i_d} & X_{d-1} \ar[r]^-{i_{d-1}} & \cdots \ar[r]^-{i_2} & X_1 \ar[r]^-{i_1} & X_0 \\
 \Sigma^{-d} \EE_d \ar[u]^-{j_d} & \Sigma^{-(d-1)} \EE_{d-1} \ar[u]^-{j_{d-1}} &  & \Sigma^{-1} \EE_1 \ar[u]^-{j_1} &  \EE_0 \ar[u]^{j_0}_-{\simeq} }\end{align}
 in the $K(n)$-local category such that
 \begin{enumerate}[(a)]
 \item the sequences
 \begin{equation}\label{eq:exactfirst}
\xymatrix@C=1pc{   \Sigma^{-k} \EE_k \ar[rr]^-{ j_k} & &X_k  \ar[rr]^-{ i_k} & & X_{k-1} \ar[rr]^-{\ell_k} & & \Sigma^{-k+1} \EE_k   }\end{equation}
 are exact triangles, and
 \item the $\EE_k$s are finite wedges of suspensions of spectra of the form $E_n^{hF}$ for finite subgroups $F$ of $\GG_n$. 
  \end{enumerate}
 \end{definition}
 In other words, a finite resolution is a tower of fibrations resolving $L_{K(n)}S^0$ by spectra of the form $E_n^{hF}$ in a finite number of steps using a finite number of pieces. Typically, $d=n^2$. Note that the tower \eqref{eq:finres} gives a diagram

\begin{align}\label{eq:resunrav}
\xymatrix{  L_{K(n)}S^0 \ar[r]^-{\delta_0} & \EE_{0} \ar[d]_-{\simeq}^-{j_0}  \ar[r]^-{\delta_{1}}  & \EE_{1} \ar[r]^-{\delta_{2}}  \ar[d]^-{\Sigma j_1}  &  \EE_{2}   \ar[r] & \cdots   \ar[r]^-{\delta_{d}}    & \EE_d  \ar[d]^-{\Sigma^d j_d}  & \\
   & X_0  \ar[ur]_{\ell_1} & \Sigma X_1 \ar[ur]_-{\Sigma \ell_2}  & & & \Sigma^d L_{K(n)}S^0 & }  
\end{align}
where $\delta_0$ is defined so that $j_0  \delta_0 = i_1  \ldots  i_d$.
We often denote the finite resolution by the top line of this diagram.
 
We can also smash \eqref{eq:finres} (in the $K(n)$-local category) with a spectrum $Y$ to obtain a tower of fibrations resolving $L_{K(n)} Y $. 
 
For any $X \in \Sp$, a resolution of the form \eqref{eq:finres} gives rise to a strongly convergent spectral sequence
 \begin{align*}
 E_1^{s,t} &= [X,  \EE_s \smsh Y ]_t \Longrightarrow [X, L_{K(n)} Y ]_{t-s},
 \end{align*}
 with differentials $d_r\colon  E_r^{s,t} \to E_r^{s+r, t+r-1}$ that collapses at the $E_{d+1}$-page. There is also a similar spectral sequence computing $ [L_{K(n)}Y,X  ]$.

\begin{example}\label{ex:finres1}
The proto-example of such a resolution is the resolution \eqref{eq:K1localfiber}. Recall that $E_1$ is $p$-completed $K$-theory and let $\mu$ be as in Section~\ref{sec:K1}. The fiber sequence \eqref{eq:K1localfiber} can be rearranged into a (very short) tower of fibrations
\[
\xymatrix{ L_{K(1)}S^0 \ar[r]^-{i_1} & E_1^{h\mu}=X_0 \\
\Sigma^{-1}E_1^{h\mu} = \Sigma^{-1}\EE_1. \ar[u]^-{j_1}
 &   E_1^{h\mu} = \EE_0\ar[u]^-{j_0}_{\simeq} }
 \]
In this case, the associated Bousfield--Kan spectral sequence degenerates to the long exact sequence on homotopy groups.
\end{example} 

For the rest of this section, we give an overview of how such resolutions are constructed. Note that the art of building finite resolutions has evolved in the last fifteen years. For a long time, the role of the Galois group was not as clear as it has become recently in the work of Henn in \cite{henn_centr}, so we give a revised recipe here.

\subsection{Algebraic resolutions}\index{algebraic resolution}

In practice, the first step to constructing a finite topological resolution is to construct its algebraic ``reflection''. These are the finite algebraic resolution. 
In practice, experts do not work from a definition, but rather know a finite algebraic resolution when they see one. Because of this, we give the following loose \emph{description} as opposed to \emph{definition}.

\begin{descript}\label{def:finalgres}
 A \emph{finite algebraic resolution} of length $d$ is an exact sequence
\begin{align}\label{eq:algres}
\xymatrix{0 \ar[r] & \cC_d  \ar[r]^-{\partial_{d}} & \cC_{d-1} \ar[r]^{\partial_{d-1}} & \cdots \ar[r]^-{\partial_1} & \cC_0 \ar[r]^-{\partial_{0}} & \mathcal{C}_{-1}= \WW_n \ar[r] & 0, }
\end{align}
where the $\cC_k$s are $\WG{\GG_n}$-modules that have the property that, for some $\EE_k$ as in Definition~\ref{eq:finres} (b), there is an isomorphism
\begin{equation}\label{eq:fortheCk} (E_n)_*^{\vee}\EE_k \cong \Hom_{\WW_n}^c(\cC_k, (E_n)_*).\end{equation}
\end{descript}

Roughly, a topological resolution \emph{realizes} an algebraic topological resolution if there is an isomorphism of exact sequences
\[\Hom_{\WW_n}^c( \cC_\bullet , (E_n)_*)  \cong  (E_n)_*^{\vee}(\EE_\bullet).\]
Here $\cC_{\bullet}$ is as in \eqref{eq:algres} and $\EE_\bullet$ is the top row of \eqref{eq:resunrav}. In this sense, the algebraic resolution is a ``reflection'' of the topological resolution.

\begin{remark}
Recall from \eqref{eq:uparrow} that $M \upar{F}{\GG_n}  = \WG{\GG_n} \otimes_{\WW_n\langle F \rangle} M $. 
Typical examples for the modules $\cC_k$ are among the following:
\begin{itemize}
\item
If $\cC_k$ is a direct sum of modules of the form $\WW_n\upar{F}{\GG_n}$ for $F$ a finite subgroup of $\GG_n$, then $\cC_k$ satisfies \eqref{eq:fortheCk}. 
Indeed, it was mentioned in \eqref{rem:periodicity} that for any $m\in \Z$ and $F$ a finite subgroup of $\GG_n$, there are isomorphisms
\begin{align*}
\Hom_{\WW_n}^c(\WW_n\upar{F}{\GG_n}, (E_n)_*)  & \cong \Map^c(\GG_n/F, (E_n)_*) \cong (E_n)_*^{\vee}\Sigma^{md_F^{\alg}}E_n^{hF}.
\end{align*}

\item By a \emph{character} $\chi$ of $\W_n \langle F \rangle$, we will mean a $\W_n \langle F \rangle$-module which has rank one over $\W_n$.
Suppose that $\chi$ is a summand (as a $\WW_n\langle F \rangle$-module) in $\WW_n\langle F \rangle$ and that $e_{\chi}$ is an idempotent of $\WW_n\langle F \rangle$ that picks up $\chi$. Let $E_n^{\chi}$ be wedge summand of $E_n$ associated to this idempotent, obtained as the telescope on $e_\chi \colon E_n \to E_n$. Then,
\[(E_n)_*^{\vee}E_n^{\chi} \cong \Hom_{\WW_n}^c(\chi\upar{F}{\GG_n}, (E_n)_*).\]
In existing examples, some of the summands of the terms $\cC_k$s are built out of projective characters $\chi$ of $\W \langle F \rangle$, such that $E_n^\chi$ is a suspension of $E_n^{hF}$. See, for example, \cite[Section 5]{ghmr} and Section~\ref{sec:resheight2} below. 
\end{itemize}
\end{remark}

\begin{remark}
One reason for using $\WW_n$-coefficients (which don't seem to play a role in the topological story) rather than $\Z_p$-coefficients in these constructions is that, if $p$ divides $n$, $\GG_n$ is ``cohomologically larger'' than $\mathbb{S}_n$ over $\Z_p$, but not over $\WW_n$ since the later is free over $\Gal$. So, if one wants to construct a resolution of length $n^2$ for $L_{K(n)}S^0 \simeq E_n^{h\GG_n}$ in cases when $p$ divides $n$, the right approach appears to be to work over $\WW_n$, and not over $\Z_p$. See also Remark~\ref{rem:HWinsight} below.
\end{remark}

We now give an outline of the steps one follows to construct a finite algebraic resolution. In practice, to construct such a resolution, it is essential to have some control over the homology 
$H^c_*(U, \WW_n)$ for an open subgroup $U$ of $\GG_n$ of finite cohomological dimension. In fact, all the examples of finite algebraic resolutions which we describe below restrict to a projective resolution of $\WW_n$ as a $\WG{U}$-module for some choice of $U$. This motivates the name of \emph{resolutions} for these exact sequences. In practice, if $p$ is large with respect to $n$ so that $\mathbb{S}_n$ has finite cohomological dimension, the finite algebraic resolutions are projective resolutions of $\WW_n$ as a $\WG{\GG_n}$-modules.

The process is inductive and goes as follows. Suppose that the $\WG{\GG_n}$-modules $\cC_{i}$ for $i\leq k-1$ together with maps $\partial_{k-1} \colon \cC_{k-1} \to \cC_{k-2}$ of  $\WG{\GG_n}$-modules have been defined so that 
\begin{align*}
\xymatrix{ \cC_{k-1} \ar[r]^{\partial_{k-1}} & \cC_{k-2} \ar[r]^-{\partial_{k-2}} & \cdots \ar[r]^-{\partial_0} & \mathcal{C}_{-1}= \WW_n \ar[r] & 0 }
\end{align*}
is an exact sequence. Suppose further that each term restricts to a projective $\WG{U}$-module. Let $N_{k-1} = \ker(\partial_{k-1})$. The projectivity assumption implies that 
\[
\Tor^{\WG{U}}_0(\WW_n, N_{k-1}) = H_{0}^c(U, N_{k-1}) \cong H_{k}^c(U, \WW_n).
\]
This isomorphism, the knowledge of $H_{k}^c(U, \WW_n)$ and a generalized form of Nakayama's Lemma \cite[Lemma 4.3]{ghmr} allows us to identify a set of $\WG{\GG_n}$-generators for $N_{k-1} \subseteq \cC_{k-1}$. The trick then is to choose a $\WG{\GG_n}$-module $\cC_{k}$ of the desired form (preferably as ``small'' as possible) and to construct a map $f \colon \cC_k \to N_{k-1}$ which surjects onto this set of generators. The map $ f$ is surjective by construction since it is chosen to make $\Tor_{\WG{U}}^{0}(f, \F_{p^n})$ surjective. The map $\partial_{k} \colon \cC_{k} \to \cC_{k-1}$ is then defined to be the composite $ \cC_{k} \to N_{k-1} \to \cC_{k-1}$, completing the inductive step. 

The process stops once $\partial_{d-1}$ has been defined. At this point, we define $\cC_d= N_{d-1}=\ker(\partial_{d-1})$ and prove that $\cC_d$ is a $\WG{\GG_n}$-module of the required type. Of course, this need not be the case and proving that this happens for some series of choices of modules $\cC_k$ and maps $\partial_k$ is usually difficult.

\begin{remark}[Algebraic resolution spectral sequence]\label{rem:ARSS}
If one resolves \eqref{eq:algres} into a double complex $P_{\bullet,\bullet}$ where $P_{\bullet,k} \to  \mathcal{C}_k$ for $0\leq k \leq d$ is a projective resolution as $\WG{\GG_n}$-modules, then the totalization of the double complex $P_{\bullet,\bullet}$ is a projective resolution of $\WW_n$. 
For a (graded) profinite $\WG{\GG_n}$-module $M = \{M_t\}_{t\in \Z}$, let $E_0^{s,k,t} = \Hom_{\WG{\GG_n}}^c( P_{k,s}, M_t)$ and take the vertical cohomology (i.e., with $k$ fixed). The result is the $E_1$-term of a spectral sequence
\[E_1^{s,k,t} =\Ext_{\WG{\GG_n}}^s(\cC_k, M_t) \Longrightarrow H^{s+k}_c(\GG_n, M_t)\]
with differentials $d_r \colon E_1^{s,k,t} \to E_1^{s+r, k+r-1,t}$. If the $\cC_k$s are direct sums of modules of the form $\chi \upar{F}{\GG_n}$ for characters $\chi$, then the $E_1$-term is easy to compute since by a version of Shapiro's lemma, we have
\begin{align*}\Ext_{\WG{\GG_n}}^s(\chi \upar{F}{\GG_n} , M_t) \cong \Ext_{\WW_n\langle F\rangle}^s(\chi , M_t) .
\end{align*}
We call this an \emph{algebraic resolution spectral sequence}.\index{algebraic resolution spectral sequence}
\end{remark}
Finally, applying the functors $\Hom_{\WW_n}^c(-, (E_n)_*)$ to \eqref{eq:algres} gives an exact sequence in the category of Morava modules $\moravamod{n}$:
\begin{align}\label{eq:moravan} 
\xymatrix@C10pt{0 \ar[r] & (E_n)_* \ar[r]^-{\partial^{0}}  &  \Hom_{\WW_n}^c(\cC_0, (E_n)_*) \ar[r]^-{\partial^{1}} & \cdots \ar[r]^-{\partial^{d}} & \Hom_{\WW_n}^c(\cC_d, (E_n)_*) \ar[r] & 0,}
\end{align}
where the maps $\partial^{k}$ are induced by $\partial_k$.

\subsection{Topological resolutions}\index{topological resolution}
With an algebraic resolution  \eqref{eq:algres} in hand, the next step is to prove that it has a topological realization which is a finite resolution of $L_{K(n)}S^0$. That is, one wants to construct a finite resolution 
\begin{equation}\label{eq:seqEE}
\xymatrix{\EE_{-1}=L_{K(n)}S^0 \ar[r]^-{\delta_{0}} & \EE_0 \ar[r]^-{\delta_1} & \cdots \ar[r]^-{\delta_{d-1}} & \EE_d}
\end{equation}
in the sense of Definition~\ref{defn:finresolution}
with the property that applying the functor 
\[
(E_n)_*^{\vee}(-) \colon \Sp \longrightarrow  \moravamod{n} 
\]
to this sequence gives rise to a complex of Morava modules isomorphic to \eqref{eq:moravan}.

By our choice of $\cC_k$s (see Description~\ref{def:finalgres}), there are isomorphisms of Morava modules
$ (E_n)_*^{\vee}\EE_k \cong \Hom_{\WW_n}^c(\cC_k, (E_n)_*)$ for non-uniquely determined
spectra $\EE_k$ of the form specified in part (b) of Definition~\ref{defn:finresolution}. The non-uniqueness of the $\EE_k$s comes from the freedom in choosing the values of $m$ above. (Note that the spectrum $E_n^{hF}$ itself is periodic with periodicity some multiple $d_{F}^{\mathrm{top}}$ of $d_{F}^{\alg}$ so there is a limited number of choices.)
Fixing some choice of $\EE_k$s, we can identify \eqref{eq:moravan} with
\[
\xymatrix@C15pt{ 0\ar[r] & (E_n)_*^{\vee}  \ar[r]^-{\partial^{0}} &  (E_n)_*^{\vee}\EE_0 \ar[r]^-{\partial^{1}}  & (E_n)_*^{\vee}\EE_1 \ar[r]  & \cdots \ar[r]^-{\partial^{d}} & (E_n)_*^{\vee}\EE_d \ar[r] & 0 .}
\]

To obtain a topological realization, one must also show that the maps $\partial^k$ are of the form $(E_n)_*^{\vee}(\delta_k)$ for maps of spectra $\delta_k \colon \EE_{k-1} \to \EE_k$. Note that this being the case can depend on the choices of $\EE_k$s. The existence of $\delta_k$ is established using a Hurewicz homomorphism
\[
\xymatrix{[\EE_{k-1}, \EE_{k}] \ar[r] & \Hom_{(E_n)_*^{\vee}E_n}( (E_n)_*^{\vee}\EE_{k-1},(E_n)_*^{\vee}\EE_{k}).} 
\]
See Proposition 2.7 \cite{ghmr} for more details. 
%
%

Even if the $\delta_k$s exist, it still does not imply that any choice of $\EE_k$s and $\delta_k$s give a finite resolution in the sense of Definition~\ref{defn:finresolution}. For this to be the case, one must have first that the compositions $\delta_k \circ \delta_{k-1}$ are null-homotopic. If such choices exists, then we inductively define
 spectra $X_k$ and maps $\ell_k$ so that
\[ 
\xymatrix{\Sigma^{k-1} X_{k-1} \ar@{.>}[rr]^-{\Sigma^{k-1} \ell_{k}} & & \EE_k \ar[rr]^-{\Sigma^k j_k} & & \Sigma^k X_k  \ar[rr]^-{\Sigma^k i_k}  & & \Sigma^{k} X_{k-1}}
\]
are exact triangles (see \eqref{eq:exactfirst}). That is, if the map $\ell_k$ can be chosen so that $\delta_{k+1} \circ \Sigma^{k-1} \ell_{k}$ is null-homotopic, then $X_{k+1}$ is defined as the cofiber of $ \Sigma^{-1} \ell_k \colon \Sigma^{-1}  X_{k-1} \to \Sigma^{-k}\EE_k$ and there exists a map $\Sigma^{k}X_{k}  \xrightarrow{\Sigma^{k}\ell_{k+1}} \EE_{k+1}$ which factorizes $\delta_{k+1}$.

To prove that $X_d \simeq L_{K(n)}S^0$, one needs to check that the map $\delta_0$ lifts along the tower
\[ \xymatrix{L_{K(n)}S^0 \ar[d]_-{\delta_0} \ar@{.>}[dr]  \ar@{.>}[drrr] & & &   \\
X_0  & X_1 \ar[l]^-{i_1} & \ldots  \ar[l] & X_d \ar[l]^-{i_d} }\]
 to a map $L_{K(n)}S^0 \to X_d$. If the lift exists, it will induce an isomorphism $(E_n)_*^{\vee}(L_{K(n)}S^0) \xrightarrow{\cong} (E_n)_*^{\vee}(X_d)$ so will be a $K(n)$-local equivalence.

%
%

\begin{remark}[Doubling up]\label{rem:doublingup}
In \eqref{eq:G1ext} above, we defined a normal subgroup $\GG_n^1 \subseteq \GG_n$ with the property that $\GG_n \cong \GG_n^1 \rtimes \Z_p$ where the extension splits whenever $n$ is coprime to $p$. 
In practice, one first constructs a finite resolution of $\WW_n$ as a $\WG{\GG_n^1}$-module and then upgrades it to a resolution of $\WW_n$ as a $\WG{\GG_n}$-module. See Corollary 4.2 of \cite{ghmr} for an example.
\end{remark}

\subsection{Diagram of resolution spectral sequences}

The resolutions whose construction is described above give rise to spectral sequences which fit in a diagram:
\begin{equation*}
\xymatrix@C=1.5pc{E_1^{s,k,t}=\Ext_{\WG{\GG_n}}^s(\cC_k, (E_n)_t) \ar@{=>}[rr]^-{\mathrm{ARSS}} \ar@{=>}[d]_-{\mathrm{LHFPSS}} &  & H_{c}^{s+k}(\GG_n,(E_n)_t) \cong E_2^{s+k,t} \ar@{=>}[d]^-{\mathrm{HFPSS}} \\
E_1^{k,t-s}=\pi_{t-s}\EE_k \ar@{=>}[rr]_-{\mathrm{TRSS}}  & & \pi_{t-(s+k)}L_{K(n)}S^0.}
\end{equation*}
Here ARSS stands for \emph{algebraic resolution spectral sequence}\index{algebraic resolution spectral sequence}, TRSS for \emph{topological resolution spectral sequence}\index{topological resolution spectral sequence}, HFPSS for \emph{homotopy fixed point spectral sequence}\index{homotopy fixed point spectral sequence} and in LHFPSS, the L is for \emph{level-wise}. The horizontal spectral sequence have the advantage of being first quadrant spectral sequences which are zero in degrees $k>d$ and so collapse at the $E_{d+1}$-page. By Remark~\ref{rem:strongvanishing}, the vertical spectral sequences also collapse at some finite stage with a horizontal vanishing line.\index{horizontal vanishing line}

\subsection{Finite resolutions at height $n=2$}\label{sec:resheight2}

Now we give examples of some of the finite resolutions at height $n=2$ that exist in the literature. In the references cited, the algebraic resolutions are usually constructed in the category $\Z_p[\![G]\!]$ for $G=\GG_2$ or $G=\GG_2^1$. As is explained Remark~\ref{rem:WGmodules}, we can transport the constructions to the category of $\GMod{\WW_2}{G}$ and this is what we do here. The reason for this change is explained in Remark~\ref{rem:HWinsight}.

\begin{notation}[Finite subgroups and their modules]
The maximal finite subgroups of  $\mathbb{S}_n$ were given in Table~\ref{tab:finitesub}. Here, we discuss them more specifically in the case $n=2$. Note that in the cases $p=2,3$, $\bbS_2$ contains $p$-torsion and so has more interesting finite subgroups (see (2) and (3) below).
\begin{enumerate}[(1)]
\item Let $p$ be odd.
Let $\sigma \in \abGal = \Gal(\F_{p^2}/\F_p)$ be the Frobenius and $\omega \in \F_{p^2}^{\times}$ be a primitive $(q=p^2-1)$th root of unity. The group $\abGal$ acts on $\F_{p^2}^{\times}$ by $\sigma(\omega) = \omega^\sigma =\omega^p$. We define
\[
F=F_{2q}:=  \F_{p^2}^{\times} \rtimes \abGal.
\]
For example, if $p=3$, then $F \cong SD_{16}$, the semi-dihedral group of order 16.
We let $\omega \in \WW_2\cong \Z_p(\omega)$ denote the Teichm\"uller lift of the same named class in $\FF_{p^2}$. The Teichm\"uller lifts then specify an embedding of $ \F_{p^2}^{\times}$ in $\mathbb{S}_n$, and so of $F$ in $\GG_n$.

We define $\WFT{F}$-modules $\chi_i^{+}$ and $\chi_i^{-}$ for $0\leq i\leq q-1$ as follows. The underlying $\WW_2$-module of $\chi_i^{\pm}$ is $\WW_2$. For $x\in \chi_i^{\pm}$, define $\omega_*(x)=\omega^i x$ and
\begin{align*}
\sigma_*(x)&= \begin{cases} x^{\sigma} & x\in \chi_i^+ \\
-x^{\sigma} & x\in \chi_i^-. 
\end{cases} 
\end{align*}
Let $\chi_i = \chi_i^{+} \oplus  \chi_i^{-} $. The twisted group ring completely decomposes as a $\WFT{F}$-modules as 
\[
\WFT{F} \cong \bigoplus_{i=0}^{q-1} \chi_i = \bigoplus_{i=0}^{q-1} \chi_i^{+} \oplus  \chi_i^{-}.  
\]
To see this isomorphism, let $x_i \in \WFT{F}$ be given by
\[
x_i = [e] +\omega^{-i}[\omega]+\omega^{-2i}[\omega^2] + \cdots + \omega^{-(q-2)i} [\omega^{q-2}]
\]
for $0\leq i \leq q-2$. The elements $x_i$ together with the elements $x_i[\sigma]$ generated $\WFT{F}$ as a $\WW_2$-module.
Furthermore, $\omega_*(x_i) = \omega^i x_i$ and $\sigma_*(x_i) =x_i[\sigma]$. So, there are isomorphisms $\chi_i^{+} \cong \WW_2\{x_i + x_i [\sigma]\}$ and $\chi_i^{-} \cong \WW_2\{x_i - x_i[\sigma]\}$.

(Note that the $\Z_p$-module $\lambda_i$ of \cite{henn_res} has the property that 
\[\WW_2 \otimes_{\Z_p} \lambda_i \cong \chi_{-i}^+ \oplus \chi_{-pi}^{+}\] 
when viewed as a $\WFT{F}$-module as described in Remark~\ref{rem:WGmodules}.)

\item
For $p=3$, let $G_{24}$ be an extension of $C_3 \rtimes C_4 \subseteq \bbS_2$ to $\mathbb{G}_2$ in the sense of Definition~\ref{defn:extendsubgroups}. We can give an explicit choice as follows. 
The subgroup of $\bbS_2$ generated by $ s= \frac{1}{2}(1+\omega \xi)$ and  $t=\omega^2$ is isomorphic to $C_3 \rtimes C_4$. Here, $s$ is of order $3$, $t$ is of order $4$, and $t s t^{-1}=s^2$. We let
$G_{24}$ be the group generated by $s$,  $t$,  $\psi=\omega \xi$ in  $\D_2^{\times}/\xi^2 \cong \GG_2$. Note that $\psi s =s\psi$ and $t \psi = \psi t^3$. The group $G_{24}$ is an extension of $C_3 \rtimes C_4$ in $\GG_2$. See Section 1.1 of \cite{ghmr}. Note that $C_3$ is normal in $G_{24}$. Therefore, the $\chi_{i}^{\pm 1}$ inherit a $G_{24}$-module structure via the map 
\[
G_{24} \to G_{24}/C_3 \cong (\F_9^{\times})^2 \times \Gal \xrightarrow{\subseteq} SD_{16},
\]
where $(\F_9^{\times})^2 \cong C_4$ denotes the subgroup of squares in $\F_9^{\times}$.

\item
For $p=2$, the group $\mathbb{S}_2$ contains a unique conjugacy class of maximal finite subgroups isomorphic to the binary tetrahedral group $T_{24} \cong Q_8 \rtimes \F_4^{\times}$. There is a choice of $T_{24}$ generated by $\omega \in \F_4^{\times}$ and an element of order four which we denote by $i \in Q_8$ with the property that $i^2=-1 \in \bbS_2$. For $j = \omega i \omega^{-1}$, the elements $i$ and $j$ satisfy the usual quaternion relations.
We let $G_{48}$ be the extension of $T_{24} \subseteq \bbS_2$ to $\GG_2$ given by $G_{48} = \langle \omega, 1+i \rangle \subseteq \D_2^{\times}/\xi^2 \cong \GG_2$. The group  $G_{48}$ is isomorphic to the binary octahedral group.
See \cite[Lemma 2.1, 2.2]{henn_centr}. 

We also let $C_2 = (\pm 1)\subseteq \bbS_2$. Define $V_4=C_2 \times \Gal$ and
$G_{12}= (C_2 \times \F_4^{\times}) \rtimes \Gal \cong C_2 \times \Sigma_3$. The group $C_4 = \langle i\rangle \subseteq \mathbb{S}_2$ also extends to a finite subgroup of $\GG_2$ which is a cyclic group of order eight given by $C_8  =\langle 1+i\rangle   \subseteq \mathbb{D}_2^{\times}/\xi^2$.

Finally, we let $\pi=1+2\omega$, which has the property that $\det(\pi)=\pi \pi^{\sigma}=3 \in \Z_2^{\times}/(\pm 1)$ is a topological generator. We let $G_{48}' = \pi G_{48} \pi^{-1}$. This group is conjugate to $G_{48}$ in $\GG_2$, but not in $\GG_2^1$.

\end{enumerate}
\end{notation}

\begin{remark}\label{rem:tmf}
The spectra $E_2^{hF}$ for $F$ finite are often equivalent to the $K(2)$-localization of topological modular forms with level structures (see \cite{behrens_chapter}). For example, 
\begin{align*}
L_{K(2)} {TMF} &\simeq E_2^{hG_{24}}  ,  &  L_{K(2)} {TMF}_0(2) &\simeq E_2^{hQ_8}   \simeq E_2^{hSD_{16}} \vee \Sigma^8 E_2^{hSD_{16}} 
\end{align*}
at $p=3$, and at $p=2$, 
\begin{align*}
L_{K(2)} {TMF}  &\simeq E_2^{hG_{48}},  &  L_{K(2)} {TMF}_0(3) &\simeq  E_2^{hG_{12}},    & L_{K(2)} {TMF}_0(5)&\simeq  E_2^{hC_8}. 
\end{align*}
Note that, in \cite{behrens_modular}, Behrens writes $ L_{K(2)} {TMF}_0(2) \simeq E_2^{hD_8} $ (for $p=3$). The difference comes from the fact that he is using the formal group law $\Gamma_C$ of a super-singular curve while we are using the Honda formal group law $\Gamma_2$. The groups $Q_8$ and $D_8$ are extensions to $\GG(\F_8, \Gamma_2)$ and $\GG(\F_8, \Gamma_C)$ respectively of the subgroup $C_4 \subseteq \Aut_{\bF_3}(\Gamma_2) \cong \Aut_{\bF_3}(\Gamma_C) $. See Remark~\ref{rem:choiceGamma}.
\end{remark}

Recall from Remark~\ref{rem:doublingup} that in practice, one begins by constructing a finite resolution of the group $\GG_2^1$. (The groups $\GG_2^1$ and $\bbS_2^1$ were defined right before \eqref{eq:G1ext}.)
\begin{example}[The Duality Resolutions]\label{ex:duality}\index{Duality resolution}
The following examples of resolutions of $\WW_2$ as a $\WGT{\GG_2^1}$-module have been coined the \emph{duality resolutions}. They are self-dual in a suitable sense (see \cite[Section 3.4]{HKM} or \cite[Section 3.3]{beaudry_res}). In fact, this duality is related to the virtual Poincar\'e duality of the group $\GG_2^1$.
 They are given by exact sequences
\begin{equation}\label{eq:dualityalg}
\xymatrix{0 \ar[r] & \mathcal{D}_3 \ar[r] &  \mathcal{D}_2 \ar[r] &  \mathcal{D}_1 \ar[r] &  \mathcal{D}_0 \ar[r] & \WW_2 \ar[r] & 0}
\end{equation}
such that each $\mathcal{D}_i$ is isomorphic to a direct sum of modules of the form 
\begin{align}\label{eq:DCtype}
\chi\upar{H}{\GG_2^1} := \WGT{\GG_2^1}\otimes_{\WFT{H}} \chi
\end{align}
for $H$ an extension to $\GG_2^1$ of a finite subgroup of $\bbS_2^1$ and $\chi$ is $\WFT{H}$-module which restricts to a free module of rank one over $\WW_2$.

They are minimal in the sense that their associated algebraic resolution spectral sequence
\begin{equation}\label{eq:dualADSS}
E_1^{r,q} = \Ext_{\WGT{\GG_2^1}}^q(\mathcal{D}_r, \F_{p^2}) \Longrightarrow H_c^{p+q}(\GG_2^1, \F_{p^2})  \end{equation}
collapses at the $E_1$-term. They can be realized as finite resolutions of $E_2^{h\GG_2^1}$
\begin{equation}\label{eq:dualitytop}
\xymatrix{E_2^{h\GG_2^1} \ar[r] & \EED_0 \ar[r] &  \EED_1 \ar[r] & \EED_2 \ar[r] & \EED_3.}
\end{equation}

\begin{enumerate}[(a)]
\item Let $p\geq 5$. There is an exact sequence \eqref{eq:dualityalg} with 
\begin{align*}
\mathcal{D}_0 &\cong \mathcal{D}_3 \cong \WW_2\upar{F}{\GG_2^1}, & \mathcal{D}_1 &\cong \mathcal{D}_2 \cong (\chi_{p-1}^+\oplus  \chi_{1-p}^{+}) \upar{F}{\GG_2^1} \ .
\end{align*}
This can be realized as a finite resolution \eqref{eq:dualitytop} with 
\begin{align*}
\EED_0& \simeq \EED_3 \simeq E_2^{hF}, & \EED_1 &\simeq \EED_2 \simeq \Sigma^{2(p-1)} E_2^{hF} \vee \Sigma^{2(1-p)}  E_2^{hF}.\end{align*}
These were constructed by Henn \cite[Theorem 5]{henn_res}. See also Lader \cite{lader}.

 \item Let $p=3$. There is an exact sequence \eqref{eq:dualityalg} with 
 \begin{align*}
 \mathcal{D}_0 &\cong \mathcal{D}_3 \cong \WW_2\upar{G_{24}}{\GG_2^1},   & \mathcal{D}_1 &\cong \mathcal{D}_2 \cong \chi_4^{-} \upar{SD_{16}}{\GG_2^1} \ . \end{align*}
This can be realized as a finite resolution \eqref{eq:dualitytop} with 
\begin{align*}
\EED_0  &\simeq E_2^{hG_{24}}, & \EED_1 &\simeq \Sigma^8 E_2^{hSD_{16}}, & \EED_2 &\simeq \Sigma^{40}E_2^{hSD_{16}}, & \EED_3 &\simeq \Sigma^{48}E_2^{hG_{24}}.
\end{align*}
These were constructed by Goerss, Henn, Mahowald and Rezk in \cite{ghmr}.

\item Let $p=2$. There is an exact sequence \eqref{eq:dualityalg} with 
\begin{align*}
\mathcal{D}_0 &\cong  \WW_2\upar{G_{48}}{\GG_2^1}, & \mathcal{D}_1 &\cong \mathcal{D}_2 \cong \WW_2\upar{G_{12}}{\GG_2^1}, &  \mathcal{D}_3 \cong  \WW_2\upar{G_{48}'}{\GG_2^1} \ .
\end{align*}
This can be realized as a finite resolution \eqref{eq:dualitytop} with 
\begin{align*}
\EED_0  &\simeq E_2^{hG_{48}}, & \EED_1 &\simeq  E_2^{hG_{12}}, & \EED_2 &\simeq \Sigma^{48}E_2^{hG_{12}}, & \EED_3 &\simeq \Sigma^{48}E_2^{hG_{48}}.
\end{align*}
These were constructed by Beaudry, Bobkova, Goerss, Henn, Mahowald, and Rezk in \cite{henn_res, beaudry_res, BobkovaGoerss}. In these references, the resolution is constructed for $E_2^{h\mathbb{S}_2}$. However, using the ideas of \cite{henn_centr} it is now straightforward to construct it for $E_2^{h\GG_2}$.\end{enumerate}
\end{example}

\begin{remark}\label{rem:dualityfullresp3}
If $p\neq 2$, the algebraic resolution can be doubled up in the sense of Remark~\ref{rem:doublingup} and the result can be realized topologically. For $p\geq 5$, this gives a resolution
\begin{equation*}
\xymatrix@=1.5pc{
L_{K(2)}S^0 \ar[r] & E_2^{hF} \ar[r]^-{\delta_0} & X  \vee   E_2^{hF} \ar[r]^-{\delta_1} &   X \vee   X   
    \ar[r]^-{\delta_2} &E_2^{hF} \vee X   \ar[r]^-{\delta_3}  &  E_2^{hF}
}
\end{equation*}
where $X =  \Sigma^{2(p-1)} E_2^{hF} \vee \Sigma^{2(1-p)}  E_2^{hF}  $,
and for $p=3$, we get
\begin{equation*}
\xymatrix@=0.6pc{
L_{K(2)}S^0 \ar[r] & E_2^{hG_{24}} \ar[r]^-{\delta_0} & E_2^{hSD_{16}}  \vee   E_2^{hG_{24}}   \ar[r]^-{\delta_1} &    \Sigma^{48}E_2^{hSD_{16}}\vee    E_2^{hSD_{16}}   \\
&    \ar[r]^-{\delta_2} & \Sigma^{48}(E_2^{hG_{24}} \vee   E_2^{hSD_{16}})  \ar[r]^-{\delta_3}  &   \Sigma^{48}E_2^{hG_{24}}  
}
\end{equation*}
However, the duality resolution at $p=2$ cannot be doubled up. 
\end{remark}

\begin{example}[The Centralizer Resolutions]\index{centralizer resolution}
The following two resolutions of the trivial $\WGT{\GG_2^1}$-modules are called \emph{centralizer resolutions} because their construction has as a key input Henn's Centralizer Approximation Theorem \cite[Theorem 1.4]{henn_duke}. 
 They are given by exact sequences
\begin{equation}\label{eq:centalg}
\xymatrix{0 \ar[r] & \mathcal{C}_3 \ar[r] &  \mathcal{C}_2 \ar[r] &  \mathcal{C}_1 \ar[r] &  \mathcal{C}_0 \ar[r] & \WW_2 \ar[r] & 0,}
\end{equation}
where the $\cC_i$s are of the form described in \eqref{eq:DCtype}. 
They can be realized as 
finite resolutions
\begin{equation}\label{eq:centtop}
\xymatrix{E_2^{h\GG_2^1} \ar[r] & \EEC_0 \ar[r] &  \EEC_1 \ar[r] & \EEC_2 \ar[r] & \EEC_3.}
\end{equation}
The algebraic centralizer exact sequences \eqref{eq:centalg} described below are \emph{$\mathcal{F}$-resolutions} in the sense of \cite[\S 3.5.1]{henn_res} and \cite[\S 1.2]{henn_centr}. We will not explain what this means, but a consequence of this fact is that the algebraic centralizer resolutions can be doubled up in the sense of Remark~\ref{rem:doublingup}. The downside of the centralizer resolutions is that they are ``larger'' than the duality resolutions. For example, the analogues of \eqref{eq:dualADSS} for the centralizer resolutions do not collapse. As a consequence, the associated algebraic and topological spectral sequences are less efficient for computations. Nonetheless, having different resolutions offers different perspectives and the centralizer resolutions have been crucial in recent computations. See for example \cite{BobkovaGoerss, gh_bcdual}.

\begin{enumerate}[(a)]
\item Let $p=3$. There is an exact sequence \eqref{eq:centalg} with 
\begin{align*}
\mathcal{C}_0 &\cong  \WW_2\upar{G_{24}}{\GG_2^1}, 
& \mathcal{C}_1 &\cong  \chi_4^{-} \upar{SD_{16}}{\GG_2^1} \oplus \chi_2^{+} \upar{G_{24}}{\GG_2^1},\\
 \mathcal{C}_2 &\cong  (\chi_{2}^+ \oplus \chi_{-2}^+) \upar{SD_{16}}{\GG_2^1}, 
&  \mathcal{C}_3 &\cong  \WW_2\upar{SD_{16}}{\GG_2^1} \ .
\end{align*}  
This can be realized as a finite resolution \eqref{eq:centtop} with 
\begin{align*}
\EEC_0  &\simeq E_2^{hG_{24}}, & \EEC_1 
&\simeq  \Sigma^8 E_2^{hSD_{16}} \vee \Sigma^{36} E_2^{hG_{24}}, \\
 \EEC_2 &\simeq \Sigma^{36}E_2^{hSD_{16}} \vee \Sigma^{44} E_2^{hSD_{16}}, 
 & \EEC_3 &\simeq \Sigma^{48}E_2^{hSD_{16}}.
\end{align*}
This is constructed by Henn in \cite{henn_res}. See also \cite[Section 4]{gh_bcdual}. Since $\GG_2 \cong \GG_2^1 \times \Z_3$, Remark~\ref{rem:doublingup} applies and we get a resolution of $L_{K(2)}S^0$.
\item  Let $p=2$. There is an exact sequence \eqref{eq:centalg} with 
\begin{align*}
\mathcal{C}_0 &\cong  \WW_2\upar{G_{48}}{\GG_2^1} \oplus  \WW_2\upar{G_{48}'}{\GG_2^1} , 
& \mathcal{C}_1 &\cong  \WW_2\upar{C_8}{\GG_2^1} \oplus \WW_2\upar{G_{12}}{\GG_2^1},\\
 \mathcal{C}_2 &\cong  \WW_2\upar{V_4}{\GG_2^1}, 
&  \mathcal{C}_3 &\cong  \WW_2\upar{G_{12}}{\GG_2^1} \ .
\end{align*}  
This can be realized as a finite resolution \eqref{eq:centtop} with 
\begin{align*}
\EEC_0  &\simeq E_2^{hG_{48}} \vee E_2^{hG_{48}'}, 
& \EEC_1 &\simeq   E_2^{hC_8} \vee E_2^{hG_{12}}, &
 \EEC_2 &\simeq E_2^{hV_4}, 
 & \EEC_3 &\simeq E_2^{hG_{12}}.
\end{align*}
\end{enumerate}
This is constructed by Henn \cite{henn_centr}. Note again that, as opposed to the duality resolution at $p=2$, the algebraic centralizer resolution \emph{can} be doubled up and the resulting sequence can be realized topologically to give a finite resolution of $L_{K(2)}S^0$. See \cite[Theorem 1.1, 1.5]{henn_centr}. 
\end{example}

\begin{remark}
The doubled up centralizer resolution at $n=p=2$ is very large compared to the duality resolution available at odd primes. However, there is a handicraft way to glue a duality resolution with a centralizer resolution to obtain a much smaller resolution called the \emph{hybrid resolution}.\index{hybrid resolution} It can be realized as a resolution of $L_{K(2)}S^0$ 
\begin{equation*}
\xymatrix@=0.6pc{
L_{K(2)}S^0 \ar[r] &  E_2^{hG_{48}}  \ar[r]^-{\delta_0} &  E_2^{hC_{8}} \vee  E_2^{hG_{12}}    \ar[r]^-{\delta_1} &  E_2^{hV_4}\vee    E_2^{hG_{12}}   \\
&    \ar[r]^-{\delta_2} & E_2^{hG_{12}} \vee    \Sigma^{48} E_2^{hG_{12}}   \ar[r]^-{\delta_3}  &   \Sigma^{48}E_2^{hG_{48}}.  
}
\end{equation*}
The construction of this resolution is not published but is due to the second author and Henn.
\end{remark}

\begin{remark}\label{rem:HWinsight}
Until recently, at $n=p=2$, topological and algebraic resolutions existed only for $\mathbb{S}_2$ and not $\GG_2$. In some loose sense, the reason for this was our inability to average over the action of $\Gal$. By an insight of Hans-Werner Henn \cite{henn_res} if one switches to $\WW_n$-coefficients, one can form ``weighted averages'' in $\WW_n$ and this has allowed us to upgrade our resolutions for $\mathbb{S}_2$ to resolutions for $\GG_2$. Note further that the description of the duality resolution at $n=2$ and $p\geq 5$ is much cleaner if one works over the Witt vectors.
\end{remark}

\section{Chromatic splitting, duality, and algebraicity}\label{sec:thms}

This section discusses two of the major areas of applications of the techniques developed above in local chromatic homotopy theory: the chromatic splitting conjecture and the study of duality phenomena in $K(n)$-local homotopy theory. In both cases, we start with an outline of the general picture, before specializing to a summary of the known results at height 2. We then conclude with a brief outlook to the asymptotic behavior of chromatic homotopy theory for primes large with respect to the height.

\subsection{Chromatic splitting and reassembly}\label{sec:cscmore}

In this section, we discuss the chromatic splitting conjecture\index{chromatic splitting conjecture} in more detail, putting an emphasis on new developments and points that were not discussed in \cite{cschov}.

The chromatic splitting conjecture (CSC) gives a fairly simple prediction of $L_{n-1}L_{K(n)}S^0$ in the chromatic fracture square. 
Although the original conjecture does not hold for the prime $p=2$ and $n=2$ as it was stated in \cite{cschov}, the fundamental idea behind the conjecture remains intact. The philosophy behind the CSC is that chromatic reassembly is governed by the structure of $H^*_c(\GG_n, \W)$, via the map $H^*_c(\GG_n, \W) \to H^*_c(\GG_n, (E_n)_*)$ to the $E_2$-term of the homotopy fixed point spectral sequence \eqref{eq:hfpss}. As discussed in Conjecture~\ref{conj:vanishing}, this map is expected to be an isomorphism onto $H^*_c(\GG_n, (E_n)_0)$.

The isomorphism \eqref{eq:rationaliso} implies at the very least that there is an inclusion
\[
\Lambda_{\Z_p}(x_1, x_2,  \ldots, x_{n}) \subseteq H_c^*(\GG_n, \WW_n) 
\]
for generators $x_i$ of cohomological degree $2i-1$ at all primes and heights.
Here, we always choose $x_1=\zeta_n$ for $\zeta_n$ as in \eqref{eq:zetan} and choose the other $x_i$s so that they do not map to zero in $H_c^*(\GG_n, \WW_n)/p$. Conjecture~\ref{conj:vanishing} then implies the existence of non-zero classes $x_{i}\in E_2^{2i-1,0}$ in
\[ 
E_2^{s,t} \cong H_c^{s}(\GG_n, (E_n)_t) \Longrightarrow \pi_{t-s}L_{K(n)}S^0. 
\]
Further, at heights $n\leq 2$ the following phenomena have been observed.
\begin{conjecture} If $p$ is odd, or if $p=2$ and $n$ is odd, then 
\[
E_{\infty}^{*,0} \cong  \Lambda_{\Z_p}(e_1, e_2,  \ldots, e_{n})\subseteq \pi_*L_{K(n)}S^0
\] 
for some choice of classes $e_i$ detected by a multiple of $x_i$.
If $p=2$ and $n$ is even, then 
\[
E_{\infty}^{*,0} \cong  \Lambda_{\Z_2}(f, e_1, e_2,  \ldots, e_{n})/(2f, e_{n} f) \subseteq \pi_*L_{K(n)}S^0.
\]
for $e_i$ as above and some choice of class $f$ detected by $\widetilde{\chi}_n$ (see \eqref{eq:chin}).  
\end{conjecture}

\begin{remark}
If $p$ is large with respect to $n$, then there is no ambiguity about the choice of classes $e_i$ because of the sparsity of the spectral sequence.
At $n=2$ and $p=3$, once can choose $e_2$ to be detected by $3x_2$ and at $n=p=2$, by $4x_2$.
\end{remark}

The dichotomy between odd and even heights for the prime $2$ comes from the following observations. If $n$ is odd, then the inclusion of $C_2 \cong (\pm 1)$ in $\GG_n$ splits the map $\chi_n$ and
\[  
H^*(C_2, \Z_2) \cong \Z_2[\widetilde{\chi}_n]/(2\widetilde{\chi}_n) \to  H^*(C_2, \WW_n)  \xrightarrow{\chi_n^*} H_c^*(\GG_n, \WW_n)
\] 
is an inclusion. However, if the image of $\widetilde{\chi}_n$ in $H_c^2(\GG_n, (E_n)_0)$ is non-trivial, then it must support a non-trivial $d_3$ differential since its image in the HFPSS computing $E_n^{hC_2}$ has this property by \cite[Theorem 1.3]{hahnshi}. (The image of $\widetilde{\chi}_n$ would be the class $u_{2\sigma}^{-1}a_{\sigma}^2$ which supports a non-zero $d_3$ differential.)

At $n=2$, what we observe is that $(\widetilde{\chi}_2^2)$ is the kernel of $\chi_2^*$ so that the latter induces the inclusion of $\Z_2[\widetilde{\chi}_2]/(2\widetilde{\chi}_2,\widetilde{\chi}_2^2)$ into $ H_c^*(\GG_2, \WW_2)$.
We do not know how this generalizes at even heights $n>2$.

As is discussed in \cite{cschov}, the induced maps $e_i \colon S^{1-2i} \to L_{n-1}S^0$ (for some choice of $e_i$) are conjectured to factor through $L_{n-i}S^0$. 
Since $f$ has order $2$, it induces a map $f \colon S^{-2}/2 \to L_{K(n)}S^0$, so after localizing at $E_{n-1}$, we get a map $f\colon L_{n-1} S^{-2}/2 \to L_{n-1}L_{K(n)}S^0$. The CSC as stated in \cite{cschov} did not take into account this class $f$.  Based on what we see in the case $n=p=2$, Conjecture~\ref{conj:csc2} below is a suggestion for a revised version of the CSC in its strongest form which includes $f$. Below, for spectra $X_i$, we let $\Lambda_{L_{n-1}S^0_p}(X_1, \ldots, X_n)$ be the wedge of $L_{n-1}S^0_p$ and of $X_{i_1}\smsh \ldots \smsh X_{i_j}$ for $1\leq i_1<\ldots <i_j \leq n$. Let 
\[\iota \colon L_{n-1}S^0_p \to \Lambda_{L_{n-1}S^0_p}(X_1, \ldots, X_n)\] 
be the inclusion of the $L_{n-1}S^0_p$ summand.

\begin{conjecture}[Strong CSC]\label{conj:csc2}
There is an equivalence in the category of $E_{n-1}$-local spectra
\[
L_{n-1}L_{K(n)}S^0  \simeq \Lambda_{L_{n-1}S^0_p} \left( L_{n-i}S^{1-2i} : 1\leq i \leq n  \right)
\]
if $p\neq 2$, or $p=2$ and $n$ odd. The map $\iota$ corresponds to the unit $L_{n-1}S^0_p  \to L_{n-1}L_{K(n)}S^0_p$.
If $p=2$ and $n$ is even, there is an $E_{n-1}$-local equivalence
\[
L_{n-1}L_{K(n)}S^0  \simeq \Lambda_{L_{n-1}S^0_2} \left( L_{n-i}S^{1-2i} : 1\leq i \leq n  \right)\smsh \Lambda_{L_{n-1}S^0_2}  \left(L_{n-1}S^{-2}/2 \right) . 
\]
In this case, the map $\iota \smsh \iota$ corresponds to the unit $L_{n-1}S^0_p  \to L_{n-1}L_{K(n)}S^0_p$.
\end{conjecture}

\begin{remark}
A criterion for this revision is for the conjecture at odd primes to remain as stated in \cite{cschov}. We have made what we think is a minimal modification to reflect what we see at $n=p=2$. However, other reformulations are possible and we concede that this is a somewhat arbitrary choice.
\end{remark}

Note that Conjecture~\ref{conj:csc2} implies the weak CSC (Conjecture~\ref{conj:wcsc}), saying that $L_{n-1}S^0 \to L_{n-1}L_{K(n)}S^0$ splits. 
Further, both Conjecture~\ref{conj:csc2} and Conjecture~\ref{conj:wcsc} hold if one replaces the sphere by any finite complex $X$. However, 
even Conjecture~\ref{conj:wcsc} does not hold for arbitrary spectra. In \cite{devcounterBP}, Devinatz proves that it fails for the $p$-completion of $BP$.

Before giving examples we would like to point out that, among its many consequences, the strong form of the chromatic splitting conjecture would also imply
\begin{conjecture}\label{conj:fingen}
For any $n \ge 0$ and any prime $p$, the homotopy groups $\pi_*L_{K(n)}S^0$ are degreewise finitely generated $\Z_p$-modules.
\end{conjecture}

\begin{example}\label{ex:csc}
At height $n=1$, the equivalence $L_0L_{K(1)}S^0 \simeq L_0(S^0_p \vee S^{-1}_p)$ holds for all primes 
and was discussed in Section~\ref{sec:reassemlyn1}, see \eqref{eq:cscn1first}. 
\end{example}
\begin{theorem}
At height $n=2$, if $p$ is odd, then
\begin{equation*}
L_1L_{K(2)}S^0 \simeq L_1(S^0_p \vee S^{-1}_p) \vee L_0(S_p^{-3} \vee S_p^{-4}) .
\end{equation*}
See \cite{behse2, GoerssSplit}.
If $p=2$, there is an equivalence
\begin{equation*}
L_1L_{K(2)}S^0 \simeq L_1(S^0_2\vee S^{-1}_2  \vee S^{-2}/2 \vee S^{-3}/2) \vee L_0(S_2^{-3} \vee S_2^{-4} ).
\end{equation*}
See \cite{BGH}. 
\end{theorem}

\begin{remark}
The fact that the CSC in its original form would most likely fail at $n=p=2$ was first noticed by Mark Mahowald. His intuition was based on the computations of Shimomura and Wang \cite{shimwang}, who identify $v_1$-torsion-free summands in the $E_2$-term of the Adams-Novikov Spectral Sequence for $\pi_\ast L_{K(2)}V(0)$ that are not predicted by the CSC. Their work, however, did not preclude the possibility of differentials that could have eliminated the summands not accounted for in the original statement of the CSC.
\end{remark}

The CSC is one of the key inputs for chromatic reassembly, which recovers $L_nS^0$ from $L_{K(n)}S^0$ and $L_{n-1}S^0$ via the chromatic fracture square \eqref{eq:chromaticfracture}. We discuss this a little more here. 

Let $F_n$ be the fiber of the map $L_{n-1}S^0 \to L_{n-1}L_{K(n)}S^0$, whose homotopy type is predicted by the CSC. The chromatic fracture square implies that $F_n$ is also the fiber of $L_nS^0 \to L_{K(n)}S^0$. The homotopy groups of $L_nS^0$ can be reassembled from the long exact sequence on homotopy groups 
\[
\xymatrix{\pi_{k}\Sigma^{-1} L_{K(n)} S^0 \ar[r] & \pi_k F_n \ar[r] & \pi_{k}L_nS^0 \ar[r] & \pi_k L_{K(n)}S^0 \ar[r] & \pi_{k}\Sigma F_n.}
\]

We explained chromatic reassembly at height $n=1$ in Section~\ref{sec:reassemlyn1}. At this point, we would like to at least give the reader an idea of chromatic reassembly at chromatic level $n=2$. A description of reassembly for $L_2S^0_p$ itself would be very technical, so instead, we describe the reassembly process for $L_2S^0/p$ for primes $p\geq 5$, which is significantly simpler. To do this, we first give a qualitative description of the homotopy groups of $\pi_*(L_{K(2)}S^0/p)$.

For primes $p \ge 5$, the homotopy fixed point spectral sequence is too sparse for differentials or extensions, so collapses to give
\[
\pi_m(L_{K(2)}S^0/p) \cong \bigoplus_{m=t-s} H_c^s(\GG_2, (E_2)_t/p).
\]
The groups on the right side of this isomorphism can be deduced from the computation by Shimomura and Yabe in \cite{shimyabe} and were discussed in Sadofsky \cite{sadofsky_picture}. They are computed directly using the finite resolution (1) of Example~\ref{ex:duality} by Lader in \cite{lader}. See \cite[Corollaire 4.4]{lader} and the discussion before it for an explicit description. We make a few observations about the answer:
\begin{enumerate}[(a)]
\item The homotopy groups $\pi_{*}(L_{K(2)}S^0/p)$ form a module over $\F_p[v_1]\otimes \Lambda_{\F_p} ( \zeta_2)$ for a class $v_1 = u^{1-p} \in \pi_{2(p-1)}(L_{K(2)}S^0/p)$ and $\zeta_2 \in \pi_{-1}(L_{K(2)}S^0/p)$ as in Remark~\ref{rem:3termfiber}. The group $\pi_{m}(L_{K(2)}S^0/p)$ is zero if 
\[
2k(p-1) <m < 2(k+1)(p-1) -4.
\]
\item Since $L_0 S^0/p \simeq \ast$, the chromatic fracture square gives an equivalence $ L_{1}S^0/p \simeq L_{K(1)}S^0/p$ and, similarly, $L_1L_{K(2)}S^0/p\simeq L_{K(1)}L_{K(2)}S^0/p$. Furthermore, on homotopy groups, the effect of $E_1$-localization on $L_{K(2)}S^0/p$ is to invert $v_1$. 
\item There is unbounded $v_1$-torsion in $\pi_{*}(L_{K(2)}S^0/p)$. However, the homotopy groups are finite in each degree $m$. In fact, for any class $x$ detected in $H_c^s(\GG_2, (E_2)_t/p)$ for $t<0$, if $t+ 2k(p-1)\geq 0$, then $v_1^{k}x=0$. That is, multiplication by $v_1$ never ``crosses'' the $s=t$ line in $H_c^s(\GG_2, (E_2)_t/p)$.
\item The only homotopy classes in $\pi_*L_{K(2)}S^0$ that are not $v_1$-torsion are given by $\F_p[v_1]\otimes \Lambda_{\F_p}(\zeta_2, h_0)$ for a class $h_0 \in \pi_{ 2(p-1)-1}(L_{K(2)}S^0/p)$ that is the image of the homotopy class $\alpha_1 \in \pi_{2(p-1)-1}S^0_{(p)}$. Furthermore,
\[
\pi_*(L_1L_{K(2)}S^0/p) \cong v_1^{-1}\pi_{*}(L_{K(2)}S^0/p) \cong \F_p[v_1^{\pm 1}] \otimes \Lambda_{\F_p}(\zeta_2, h_0).
\]
Under the canonical map $\pi_{*}(L_{1}S^0/p) \to \pi_{*}(L_{1}L_{K(2)}S^0/p)$, the class $\zeta_1 \in \pi_{-1}(L_{1}S^0/p)$ maps to $v_1^{-1}h_0$.
\end{enumerate}
We use the long exact sequence on homotopy groups associated to the fiber sequence
\[ F_2/p \to L_{2}S^0/p \to L_{K(2)}S^0/p \]
to deduce that
\begin{align*}
\pi_*(L_2S^0/p) &\cong \mathrm{Tor}_{v_1}(\pi_{*}(L_{K(2)}S^0/p))\oplus  \F_p[v_1] \{1,h_0\}    \oplus  \Sigma^{-1}\F_p[v_1]/(v_1^{\infty}) \{\zeta_2,h_0\zeta_2\},
\end{align*}
where $\mathrm{Tor}_{v_1}(\pi_{*}(L_{K(2)}S^0/p))$ is the $v_1$-power torsion subgroup of $\pi_{*}(L_{K(2)}S^0/p)$ and $\F_p[v_1]/(v_1^{\infty})$ is the cokernel of the canonical map in the following short exact sequence
\[
\xymatrix{0 \ar[r] & \F_p[v_1] \ar[r] & \F_p[v_1^{\pm 1}] \ar[r] & \F_p[v_1]/(v_1^{\infty}) \ar[r] & 0.}
\]

All available proofs of the CSC, even in its weakest form, have been brutally computational. Short of simply computing $\pi_*L_{K(2)}S^0$ explicitly, the steps have been: 
\begin{enumerate}[(a)]
\item Prove that there are non-zero homotopy classes $e_1$ and $e_2$, and if $p=2$ an additional class $f$ detecting a non-trivial class of order $2$. This gives the map
\[
\xymatrix{S^0 \vee S^{-1} \vee S^{-3} \vee S^{-4} \ar[r]^-{\varphi} & L_{K(2)}S^0,}
\] 
where $\varphi = 1 \vee e_1 \vee e_2 \vee e_1e_2$. If $p=2$, there is an additional factor of $S^{-2}/2 \vee S^{-3}/2 \xrightarrow{f \vee e_1f} L_{K(2)}S^0$.
\item Compute $v_1^{-1}\pi_*L_{K(2)}(\varphi \smsh X)$ for a finite type $1$ complex, usually $X=S^0/p$.
\item Compute $p^{-1}\pi_*L_{K(2)}S^0$. 
\item Reassemble the fracture square.
\end{enumerate}

\begin{remark}
Let $p$ be an odd prime. In general, the CSC predicts that $1 \vee \zeta_n$ induces an equivalence $L_{K(n-1)}S^0 \vee L_{K(n-1)}S^{-1} \xrightarrow{\simeq} L_{K(n-1)}L_{K(n)}S^0$. In particular, it implies that $L_{K(n-1)}S^0 \simeq L_{K(n-1)}E_n^{h\GG_{n}^1}$.

There is another conjecture of Hopkins related to chromatic splitting called the \emph{algebraic chromatic splitting conjecture} \cite[Section 14]{petersoneric}.
It states:
\begin{conjecture}[Algebraic CSC]\label{conj:ACSC}
Let $p$ be an odd prime, possibly large with respect to $n$. Then 
\[
\underset{i}{\lim}{}_{(E_n)_*E_n}^s (E_n)_t/(p, v_1, \ldots, v_{n-2}, v_{n-1}^i ) \cong \begin{cases}  (E_n)_*/(p, v_1, \ldots, v_{n-2}) & s=0 \\
v_{n-1}^{-1}(E_n)_*/(p, v_1, \ldots, v_{n-2}) & s=1 \\
0 & s>1.
\end{cases}
\]
\end{conjecture}
\noindent
Here,  the limit (and its derived functors) is taken in the category of $(E_n)_*E_n$-comodules, where $(E_n)_*E_n  = \pi_*(E_n \smsh E_n)$ is the group of non-completed $E_n$-cooperations. Provided that $(E_n)_*\zeta_n$ generates the ${\lim}^1$-term, Conjecture~\ref{conj:ACSC} implies that $1 \vee \zeta_n$ is an $E_n$-local equivalence, thereby verifying the chromatic splitting conjecture at height $n$ for finite type $n-1$ complexes. 
\end{remark}

\begin{remark}
As explained above, the chromatic splitting conjecture is a fundamentally transchromatic statement. In a series of papers, Torii uses generalized character theory to study the relation between adjacent strata in the chromatic filtration. In particular, he shows in \cite{torii_zeta} that under the canonical map 
\[
\xymatrix{\pi_{-1}L_{K(n-1)}S^0 \ar[r] & \pi_{-1}L_{K(n-1)}L_{K(n)}S^0}
\]
the class $\zeta_{n-1}$ maps non-trivially. 
\end{remark}

\subsection{Invertibility and duality}\label{sec:invdual}

In analogy with the problem of computing the group of units of a classical ring, an important aspect of understanding a symmetric monoidal category $(\mathcal{C}, \otimes, I )$ with unit $I$ is to classify its invertible objects. An object $X \in \cC$ is invertible if there exists another object $Y \in \cC$ such that $X \otimes Y \cong I$ where $I$ is the unit of the symmetric monoidal structure. If the collection of invertible objects forms a set, then it is an abelian group under $\otimes$ and this group is called the \emph{Picard group} of 
$\mathcal{C}$, denoted $\Pic(\mathcal{C})$. The Picard group\index{Picard group} of a symmetric monoidal $\infty$-category $\cC$ is defined to be the Picard group of the homotopy category of $\cC$, i.e., we set $\Pic(\cC) = \Pic(\mathrm{Ho}(\cC))$.

If $\cC$ is a triangulated category, the Picard group always contains a cyclic subgroup generated by the shift of the unit $I[1]$. For example, the Picard group of the stable homotopy category contains a copy of $\Z$ generated by $S^1$. In fact, in this case, there is nothing else and $\Pic(\Sp)\cong \Z
\langle S^1\rangle$, see \cite{hms_pic}. We can view $\Sp_E$ as a symmetric monoidal category with product $L_E(-\smsh -)$ and unit $L_ES^0$. The objects $L_ES^n$ for  $n\in \Z$ are invertible in this category. One of the fascinating aspects of $E$-local homotopy theory is that, for some choices of $E$, there are invertible object in $\Sp_E$ 
which are not of the form $L_ES^n$ for some $n\in \Z$. 
When $E=K(n)$ for $0< n <\infty$, the Picard group is in fact much larger; we remark that for $E=K(n)$ and $E=E(n)$ and arbitrary $n$, the collection of isomorphism classes of invertible objects in $\Sp_E$ indeed forms a set, see \cite[Proposition 7.6]{hms_pic} and \cite[Proposition 1.4]{hs_pic}.

The Picard group of the $K(n)$-local category (with $p$ fixed and suppressed from the notation) is usually denoted by $\Pic_n$. Note that, if $X \in \Pic_n$, its inverse is the Spanier--Whitehead dual of $X$ in the $K(n)$-local category: $D_nX = F(X, L_{K(n)}S^0)$. By Galois descent~\cite[Proposition 10.10]{mathew_galois}, there is an isomorphism
\[
\Pic_n \cong \Pic_{K(n)}(\Mod_{E_n}^{\GG_n}),
\]
where $\Mod_{E_n}^{\GG_n}$ is the $K(n)$-local category of $\wG$-twisted $E$-module spectra. The right hand side has a natural algebraic analogue given by
\[
\Pic_n^{\alg} := \Pic(\moravamod{n})
\]
where $\moravamod{n}$ is the category of Morava modules (see Definition~\ref{defn:moravamodules}).

A Morava module $M$ is in $\Pic_n^{\alg}$ if and only if it is free of rank one over $(E_n)_*$. Since $(E_n)_*$ is two periodic, $\Pic_n^{\alg}$ is naturally $\Z/2$-graded. Let $\Pic_n^{\alg,0}$ be the subgroup of elements such that $M \cong (E_n)_*$ as $(E_n)_*$-modules.  
The latter can then be described (but not easily computed) as
\[
\Pic_n^{\alg,0} \cong H_c^1(\wG, (E_n)_0^{\times}). 
\]

The functor which sends $X  \in  \Sp_{K(n)}$ to $(E_n)_*^{\vee}X$ induces a map,  $\Pic_n \to \Pic_n^{\alg}$, and we define $\kappa_n$ to be the kernel:
\begin{equation}\label{eq:picseq}
\xymatrix{0 \ar[r] & \kappa_n \ar[r] & \Pic_n \ar[r] & \Pic_n^{\alg}}
\end{equation}
It is called the \emph{exotic} Picard group\index{exotic Picard group $\kappa_n$} fo $\Sp_{K(n)}$. Elements of $\Pic_n$ which are in $\kappa_n$ are called exotic.

For $2(p-1) \geq n^2$, an argument that uses the sparseness of \eqref{eq:hfpss} shows that $\kappa_n = 0$ so that the map \eqref{eq:picseq} $\Pic_n \to \Pic_n^{\alg}$ is an injection \cite[Proposition 7.5]{hms_pic}. However, it has been shown in many cases that $\kappa_n$ is non-trivial. The following is a conjecture of Hopkins.

\begin{conjecture}\label{conj:fin}
The group $\kappa_n$ is a finite $p$-group.
\end{conjecture}
In \cite[Theorem 4.4.1]{Heard}, Heard proves that for $p$ odd, $\kappa_n$ is a direct product of cyclic $p$-groups. Note also that a positive answer to Conjecture~\ref{conj:fingen} would imply Conjecture~\ref{conj:fin}.

In \cite{pstragowski_pic}, Pstr{\c a}gowski proves that for $2(p-1)>n^2+n$, $ \Pic_n  \cong \Pic_n^{\alg}$. The question of whether or not $ \Pic_n \to \Pic_n^{\alg}$ is surjective in general is open. Furthermore, the algebraic Picard group $\Pic_n^{\alg}$ is not known for any prime when $n > 2$. It is believed that $\Pic_n^{\alg}$ is finitely generated over $\Z_p$. In fact, the (folklore) expectation is that $\Pic_n^{\alg}$ is of rank two over $\Z_p$, with one summand generated by $L_{K(n)}S^1$ and the other by the spectrum $S\langle{\det}\rangle$ discussed below in Example~\ref{example:elementsinpic} (b).
The table of Figure~\ref{fig:pictable} summarizes the current state of the literature on these questions. The second author, Bobkova, Goerss and Henn have been working towards identifying $\Pic_2$ when $p=2$.

\begin{table}
\captionsetup{width=\textwidth}
\caption{The table below contains some known values of $\Pic_n$. Here, H--M--S stands for Hopkins--Mahowald--Sadofsky, G--M--H--R stands for Goerss--Henn--Mahowald--Rezk and K--S for Kamiya--Shimomura.}
\label{fig:pictable}
\centering
\small{
\begin{tabular}{|c | c | c | c | c | c | }
\hline
$n$ & $p$ & $\Pic_n$ & $\Pic_{n}^{\mathrm{alg}} $ & $\kappa_n$  &  Reference \\
\hline \hline
$1$ & $\geq 3$ & $\Z_p \times \Z/2(p-1)$ &  $\Z_p \times \Z/2(p-1)$  & $0$  & H--M--S \cite{hms_pic} \\
\hline
$1$ & $ 2$ & $\Z_2 \times \Z/2 \times \Z/4$ &  $\Z_2 \times (\Z/2)^2$  & $\Z/2$  & H--M--S  \cite{hms_pic}   \\
\hline
$2$ & $\geq 5$ & $\Z_p^2 \times \Z/2(p^2-1)$ &  $\Z_p^2 \times \Z/2(p^2-1)$  & $0$  &  \makecell{Due to Hopkins \\
See Lader \cite{lader}}\\
\hline
$2$ & $3$ & $\Z_3^2 \times \Z/16 \times (\Z/3)^2$ &  $\Z_3^2 \times \Z/16 $  & $(\Z/3)^2$ &  \makecell{Karamanov \cite{karamanov} \\ G--H--M--R \cite{ghmr_pic} \\ K--S \cite{kam_shim}}  \\
\hline
\end{tabular}}
\end{table}

\begin{remark}Analogously, there is a map from the Picard group $\Pic(\Sp_n)$ of the $E_n$-local category $\Sp_n$ to the category of $(E_n)_0E_n$-comodules. In \cite{hs_pic}, Hovey and Sadofsky use a variant of Theorem~\ref{thm:vcdim} to determine these Picard groups completely for large primes: They show that for $2(p-1) > n^2+n$, there is an isomorphism $\Pic(\Sp_n) \cong \Z$, generated by $L_nS^1$. 
\end{remark}

\begin{example}\label{example:elementsinpic}
We describe a few important elements in $\Pic_n$.
\begin{enumerate}[(a)]
\item The spheres $L_{K(n)}S^m$ for $m \in \Z$ are all invertible.
\item The determinant sphere $S\langle {\det} \rangle \in \Pic_n$. See \cite{detsphere} for a construction. It has the property that $(E_n)_*^{\vee}S\langle {\det} \rangle \cong (E_n)_*$ as $(E_n)_*$-modules, but with action of $\GG_n$ twisted by the determinant. Its image in $\Pic_n^{\alg,0} \cong H_c^1(\GG_n, (E_n)_0^{\times})$ is the homomorphism $\det \colon \GG_n  \to \Z_p^{\times} \subseteq  (E_n)_0^{\times}$ of Definition~\ref{defn:determinant}.
\item\label{item:fiber} Given $\lambda \in (\pi_0E_n^{h\GG_n^1})^{\times}$, one can define a element $S^{\lambda}$ via the fiber sequence
\[\xymatrix{S^{\lambda} \ar[r] & E_n^{h\GG_n^1} \ar[r]^-{\psi -\lambda } & E_n^{h\GG_n^1} .}\]
Some variation of this construction is discussed in Section 3.6 of \cite{westerland}. If the Adams--Novikov filtration of $\lambda$ is positive, one can show that $S^{\lambda}$ is exotic. At $p=3$, the subgroup of $(\pi_0E_2^{h\GG_2^1})^{\times}$ of positive Adams--Novikov filtration is isomorphic to $\Z/3$. The elements in one of the factors of $\Z/3$ in $\kappa_2 \cong \Z/3 \times \Z/3$ are of the form $S^{\lambda}$.

\item Some exotic elements cannot be constructed using (\ref{item:fiber}). One can instead use finite resolutions to construct them. The first example is at $p=2$ and $n=1$.
Recall that $K=E_1$ and $KO\simeq K^{hC_2}$. Since $K_*^{\vee}KO \cong K_*^{\vee}\Sigma^4KO$, rather than choosing $\EE_0 = \EE_1  = KO$ to topologically realize \eqref{eq:resheight1algK}, one can let $\EE_0 = \EE_1 = \Sigma^4KO$ to get a fiber sequence
\begin{align}\label{eq:realizationP1}
\xymatrix{P_1 \ar[r] &  \Sigma^4KO \ar[rr]^{\Sigma^{4}5^{2}\psi -1} & &  \Sigma^4KO}
\end{align}
where $\psi \in \GG_1 \cong \Z_2^{\times}$ is as in \eqref{eq:psichoice}. The
fiber $P_1$ is a generator of $\kappa_1$. By construction, $P_1\smsh KO \simeq \Sigma^4 KO$. See \cite[Example 5.1]{ghmr_pic} for more details.

Similarly, at $p=3$, there is an element $P_2 \in \kappa_2$ with the property that $P_2 \smsh E_2^{hG_{24}} \simeq \Sigma^{48} E_2^{hG_{24}}$; in fact,$P_2$ is a non-trivial exotic element which generates the other summand of $\Z/3 \subseteq \kappa_2$. See \cite[Theorem 5.5]{ghmr_pic}. It is constructed by modifying the realization of the duality resolution of Remark~\ref{rem:dualityfullresp3}. 
\end{enumerate}
\end{example}

We now discuss another element of $\Pic_n$ which plays an important role in $K(n)$-local homotopy theory and brings us to the topic of Gross--Hopkins duality. As an application of the period map mentioned in Remark~\ref{rem:grosshop}, Gross and Hopkins determine the dualizing complex of $\Sp_{K(n)}$, defined via a lift of Pontryagin duality for abelian groups. More precisely, the functor 
\[
\xymatrix{I_n^*(-) = \Hom(\pi_{-*}M_n(-), \Q_p/\Z_p)\colon \Sp_{K(n)}^{\mathrm{op}} \ar[r] & \Mod_{\Z_p}^{\mathrm{graded}}}
\]
is cohomological and thus representable by an object $I_n \in \Sp_{K(n)}$, the \emph{Gross--Hopkins dual of the sphere}. From an abstract point of view, it endows $\Sp_{K(n)}$ with a Serre duality functor, see~\cite{ars}.

By Theorem \ref{thm:descentss} and Corollary \ref{cor:hfpsslarge}, $I_n$ is determined by its Morava module $(E_n)_*^{\vee}(I_n) = \pi_*L_{K(n)}(E_n \smsh I_n)$ when $p$ is large with respect to the height $n$; otherwise, one might have to twist by an exotic element of the $K(n)$-local Picard group. The spectrum $I_n$ turns out to be invertible in $\Sp_{K(n)}$, and Gross and Hopkins use the period map \eqref{eq:permap} to show that there is an equivalence
\[
I_n  \simeq L_{K(n)}(S^{n^2-n} \smsh S\langle {\det} \rangle \smsh P_n),
\]
where $P_n \in \kappa_n$ and $S\langle {\det} \rangle $ is as in Example~\ref{example:elementsinpic} (b). This identification is also known as \emph{Gross--Hopkins duality}\index{Gross--Hopkins duality}. It turns out that $P_1$ for $p=2$ and $P_2$ for $p=3$ are the elements discussed in Example~\ref{example:elementsinpic} (d).

It has now been shown in many cases where $\GG_n$ has $p$-torsion that $P_n \not\simeq L_{K(n)}S^0$ by showing that $P_n \smsh E_n^{hF} \not\simeq E_n^{hF}$ for a suitable choice of subgroup $F \subseteq \GG_n$. 
 See \cite{BBS} and \cite{heard_li_shi}. These arguments rely on the intimate relationship between Gross--Hopkins duality and $K(n)$-local Spanier--Whitehead duality\index{Spanier--Whitehead duality}:
 For $X \in \Sp_{K(n)}$, let $I_nX =F(X,I_n)$. The invertibility of $I_n$ implies that $I_nX \simeq D_nX \smsh I_n$ so studying $I_nX$ amounts to understanding $P_n$ and $D_nX$.

We end this section with a few remarks on Spanier--Whitehead duality and, more specifically, on the problem of identifying $D_nE_n$. A first answer to this question due to Gross and Hopkins (see \cite[Proposition 16]{StrickGrossHop}) states that there is a weak equivalence $\Sigma^{-n^2}E_n \to D_nE_n $ which induces an isomorphism of $\GG_n$-modules on homotopy groups. This does not imply that $D_nE_n$ is equivalent to  $\Sigma^{-n^2}E_n$ as $\GG_n$-equivariant spectra. But it does suggest that there is a \emph{dualizing module}, that is, that $D_nE_n$ is self-dual up to a twist. 

In fact, the twist can be described as the $K(n)$-localization of a $p$-adic sphere. A first description of this sphere is given as the Spanier--Whitehead dual of 
\[S^{\GG_n} :=\left( \colim_{m, \mathrm{tr}} \Sigma_{+}^{\infty}B\mathbb{S}_n^{m}\right)^{\wedge}_p,\] 
where $\mathbb{S}_n^{m} \subseteq \mathbb{S}_n$ is the subgroup of elements congruent to $1$ modulo $\xi^{nm}$ (as in \eqref{eq:stabpres}). Here, the colimit is taken over transfers and one can show that $S^{\GG_n}$ has the homotopy type of a $p$-adic sphere of dimension $n^2$. This description is not practical for computations as the action of $\GG_n$ on $S^{\GG_n}$ induced by the conjugation actions on $B\mathbb{S}_n^{m}$ is mysterious. However, there is a conjectural description of the twist analogous to the identification of the dualizing object in the Wirthm{\"u}ller isomorphism for compact Lie groups~\cite[Chapter III]{lmsm_esht}. Let $\mathfrak{g}$ be the abelian group underlying $\mathcal{O}_{\D_n}$ (as defined in Section~\ref{sec:GGn}) endowed with the conjugation (or adjoint) action of $\GG_n$. Let 
 \[S^{\mathfrak{g}} = \left(\colim_{m, \mathrm{tr}}   \Sigma_{+}^{\infty} Bp^m \mathfrak{g} \right)^{\wedge}_p.\] 
 Again, $S^{\mathfrak{g}}$ has the homotopy type of a $p$-adic sphere of dimension $n^2$ and the action of $\GG_n$ on $Bp^m \mathfrak{g}$ induces an action on $S^{\mathfrak{g}}$.

\index{Linearization hypothesis}
\begin{linhypo}\label{conj:lin}
There is a $\GG_n$-equivariant equivalence 
\[S^{ \mathfrak{g}} \simeq S^{ \GG_n} .\] 
\end{linhypo}
\noindent

Inspired by Serre's definition of the dualizing module \cite[Chapter I, \S 3.5]{serre}, such a statement was first suggested to the experts by the strong connections in the $K(n)$-local category between Spanier--Whitehead duality, Brown--Comenetz duality, and Poincar\'e duality for the group $\GG_n$ (see Gross--Hopkins \cite{HopkinsGross} and Devinatz--Hopkins \cite[Section 5]{DH_action}). The hypothesis is stated in work of Clausen \cite[Section 6.4]{clausen_padicj}, not only for $\GG_n$, but for any $p$-adic analytic group. Clausen has recently announced a proof of the Linearization Hypothesis in this general form.

The linearization hypothesis leads to a $\GG_n$-equivariant equivalence 
\[D_nE_n \simeq  L_{K(n)}(E_n \smsh S^{- \mathfrak{g}}),\] 
where $S^{-\mathfrak{g}}=  F(S^{\mathfrak{g}}, S_p^0 )$, a description that lends itself well to applications.

\subsection{Compactifications and asymptotic algebraicity}\label{sec:assymptotics}

\index{asymptotic algebraicity}
We conclude this survey with a short overview of another recent direction in chromatic homotopy theory. As explained in Section~\ref{ssec:gn} and demonstrated in the examples above, chromatic homotopy theory at a fixed height $n$ simplifies when $p$ grows large, the essential transition occurring when $p-1>n$. This leads to the question of how to isolate those phenomena that hold generically, i.e., for all primes $p$ which are large with respect to the given height. The goal of this section is to outline a result of \cite{ultra1} that describes the \emph{compactification} of chromatic homotopy theory, which is based on the notion of ultraproducts\index{ultraproduct}. In particular, this provides a model of the limit of the $K(n)$-local categories when $p \to \infty$ that captures the generic behavior of these categories. A different approach using Goerss--Hopkins obstruction theory that gives an algebraic triangulated model for suitably large primes has recently appeared in the work of Pstr{\c a}gowski~\cite{pstragowski_alg}.

The Stone--\v{C}ech compactification of a topological space $X$ is the initial compact Hausdorff space $\beta X$ equipped with a continuous map $\iota$ from $X$. If $X$ is discrete, $\beta X$ can be modeled by the set of ultrafilters\index{ultrafilter} on $X$ endowed with the Stone topology; recall that an \emph{ultrafilter} on $X$ is a set $\cU$ of subsets of $X$ such that whenever $X$ is written as a disjoint union of finitely many subsets, then exactly one of them belongs to $\cU$. The structure map $\iota \colon X \to \beta X$ sends a point $x \in X$ to the principal ultrafilter at $x$, i.e., the set of subsets of $X$ that contain $x$. We denote the set of non-principal ultrafilters suggestively by $\partial \beta X = \beta X \setminus \iota(X)$. Assuming the axiom of choice, $X$ is infinite if and only if $\partial\beta X$ is non-empty. (In fact, the existence of non-principal ultrafilters is weaker than the axiom of choice, but the key point here is that there is no constructive way to find non-principal ultrafilters.)

Moreover, one can show that an open subset of $\beta X$ containing $\partial \beta X$ misses only finitely many points of $\iota X$. This may be thought of as a topological manifestation of the fundamental theorem of ultraproducts due to {\L}o{\'s}, which says that for a collection of models $(M_i)_{i\in I}$ of a first order theory, there is an equivalence for any formula $\phi$:
\begin{equation}\label{eq:los}
\begin{Bmatrix}
\phi \text{ hold in } M_i\\
\text{ for almost all } i \in I
\end{Bmatrix} 
\xymatrix@C=1.7pc{ \ar@{<=>}[r]^-{\sim} &}
\begin{Bmatrix}
\phi \text{ hold in } \prod_{\cU}M_i  \\
\text{ for all } \cU \in \partial\beta I
\end{Bmatrix},
\end{equation} 
where $\prod_{\cU}M_i$ denotes the \emph{ultraproduct} of the $M_i$ at $\cU$. While ultraproducts at non-principal ultrafilters thus capture generic information about the collection $(M_i)_{i\in I}$, they tend to also exhibit simplifying features. For example, for $\cU \in \partial\beta\bP$ a non-principal ultrafilter on the set of prime numbers $\bP$, the ultraproduct of $(\F_p)_{p \in \bP}$ turns out to be a rational field.
 
If $\cF$ is a presheaf on a topological space $X$ with values in a coefficient category $\cC$ that is closed under filtered colimits and products, then one may construct a naive completion of $\cF$ to be the presheaf on $\beta X$ given as the composite
\[
\xymatrix{\widehat{\cF}\colon \Open(\beta X)^{\mathrm{op}} \ar[r]^-{\iota^*} & \Open(X)^{\mathrm{op}} \ar[r]^-{\cF} & \cC,}
\]
where $\iota^*$ denotes the inverse image functor. Assuming that $X$ is discrete so that $\cF$ is a just a collection of stalks $(\cF_x)_{x\in X}$, the stalk of $\widehat{\cF}$ at an ultrafilter $\cU \in \beta X$ is given by
\begin{equation}\label{eq:ultraproduct}
\xymatrix{\widehat{\cF}_{\cU} \simeq \colim_{A \in \cU}\prod_{x \in A}\cF_x,}
\end{equation}
where the filtered colimit is taken over the projection maps induced by inclusions $A \subseteq A'$ in $\cU$. In particular, if $\cU = \cU_{x_0} \in \beta X$ is principal at a point $x_0 \in X$, then $\widehat{\cF}_{\cU_{x_0}} \simeq \widehat{\cF}_{x_0}$. The formula \eqref{eq:ultraproduct} exhibits the stalk $\widehat{\cF}_{\cU}$ as a categorical generalization of ultraproducts: indeed, if all the $\cF_x$ are non-empty and the coefficients are $\cC = \Set$, then this recovers the usual notion of ultraproducts mentioned above.

In Section \ref{sec:landscape}, we saw that the points of the spectrum $\Spc(\Sp)$ are in bijective correspondence with pairs $(p,n) \in \bP \times (\N \cup \{\infty\})$ with $(p,0) \sim (q,0)$ for all $p,q \in \bP$. From the point of view of tensor triangular geometry, one may think of the category of spectra as behaving like a bundle of categories over the space $\Spc(\Sp)$. Restricting to the discrete subspace $\bP \times \{n\} \subset \Spc(\Sp)$, this bundle should then be a disjoint union of the local categories $\Sp_{K(n)}$ for varying $p \in \bP$. Therefore, the above formalism yields a diagram
\[
\xymatrix{\coprod_{p \in \bP}\Sp_{K(n)} \ar[d] & \widehat{\coprod_{p \in \bP}\Sp_{K(n)}} \ar[d] & \prod_{\cU}^{\flat}\Sp_{K(n)} \ar[d] \\
\bP \times \{n\} \ar[r]^-{\iota} & \widehat{\bP\times \{n\}} \cong \beta\bP & \{\cU\}, \ar[l]_-{\supset}}
\]
in which the right vertical arrow exhibits $\prod_{\cU}^{\flat}\Sp_{K(n)}$ as the stalk of the compactification over $\cU \in \beta\bP$. (Here, the superscript $\flat$ indicates that the coefficient category is $\cC = \Cat_{\infty}^{\flat}$ for a suitable decoration $\flat$. For the details, we refer the interested reader to \cite{ultra1,ultra1.5}.)

The final ingredient in the formulation of the asymptotic algebraicity of chromatic homotopy theory is the algebraic model itself. Informally speaking, this is given by the stable $\infty$-category of periodized ind-coherent sheaves on the formal stack $\hat{\cH}(n)$ from Section \ref{sec:geomodel}. 
An algebraic avatar of this category has previously been studied by Franke \cite{franke}. Based on \cite{ultra1}, the main theorem of \cite{ultra1.5} can now be stated as follows:

\begin{theorem}\label{thm:ultra}
For any non-principal ultrafilter $\cU$ on $\bP$ there is a symmetric monoidal equivalence
\[
\xymatrix{\prod_{\cU}^{\flat}\Sp_{K(n)} \ar[r]^-{\sim} & \prod_{\cU}^{\flat}\IndCoh(\hat{\cH}(n))^{\mathrm{per}}}
\]
of $\Q$-linear stable $\infty$-categories.
\end{theorem}

The algebraic categories $\IndCoh(\hat{\cH}(n))^{\mathrm{per}}$ admit explicit algebraic models in terms of certain comodule categories. Therefore, the above result gives a precise formulation of the empirical observation that height $n$ chromatic homotopy theory becomes asymptotically algebraic when $p \to \infty$.

\bibliographystyle{alpha}
\bibliography{handbook-beaudrybarthel-bib}

\end{document}